 \documentclass[preprint,review,11pt]{elsarticle}

\usepackage{amssymb}
 \usepackage{amsthm,amsmath}
\usepackage{mathrsfs}

 \numberwithin{equation}{section}

\newtheorem{Theorem}{Theorem}[section]
\newtheorem{Lemma}{Lemma}[section]
\newtheorem{Remark}{Remark}[section]

\newcommand{\beq}{\begin{eqnarray}}
\newcommand{\eeq}{\end{eqnarray}}
\newcommand{\beqno}{\begin{eqnarray*}}
\newcommand{\eeqno}{\end{eqnarray*}}
\newenvironment {Proof} {\noindent {\bf Proof.}}{\quad $\square$\par\vspace{3mm}}
\newenvironment {Proof of Theorem $2.1$} {\noindent {\bf Proof of Theorem $2.1$.}}{\quad $\square$\par\vspace{3mm}}
\newenvironment {Proof of Theorem $2.2$} {\noindent {\bf Proof of Theorem $2.2$.}}{\quad $\square$\par\vspace{3mm}}
\newenvironment {Proof of Theorem $2.3$} {\noindent {\bf Proof of Theorem $2.3$.}}{\quad $\square$\par\vspace{3mm}}

\allowdisplaybreaks

\newcommand\E[1]{e^{i#1\Delta}}

\journal{Journal of Functional Analysis}

\begin{document}
\newcommand{\D}{\displaystyle}
\begin{frontmatter}



\title{{\LARGE Incompressible Limit of the Compressible Nematic Liquid Crystal Flow}\emph{}}


\author{Shijin Ding$^{1}$}
\author{Jinrui Huang$^{1,*}$\corref{cor1}}
\cortext[cor2]{Corresponding author. Email: huangjinrui1@163.com}
\author{Huanyao Wen$^1$}
\author{Ruizhao Zi$^2$}
\address{$^1$School of Mathematical Sciences,
South China Normal University, Guangzhou 510631, China}
\address{$^2$Department of Mathematics, Zhejiang University, Hangzhou 310027, China}

\begin{abstract} This paper is concerned with the incompressible limit of the compressible
hydrodynamic flow of liquid crystals with periodic boundary
conditions in $\mathbb{R}^N(N=2,3)$. The local and global existence of strong solutions for the incompressible system with small
initial data is rigorously proved via the incompressible limit. Furthermore, the convergence rates are obtained in some sense.
\end{abstract}


\end{frontmatter}

% \linenumbers
\newpage
\setcounter{section}{0} \setcounter{equation}{0}
\section{Introduction}
Liquid crystals exist in an intermediate state between solid and isotropic liquid which flow like fluid or viscous fluid and also have properties
of solid crystals such as certain optical properties. For the nematic liquid crystals, the axes of constituent molecules tend to align
parallel to each other along some preferred direction $n$, which is called the anisotropic axis. The hydrodynamic theory of liquid crystals was first proposed by
Ericksen \cite{Erickson1,Erickson2} and Leslie
\cite{leslie} in 1960s, see also the book by de Gennes \cite{Genn}. However, its rigorous mathematical analysis
did not appear until 1990s when Lin \cite{1} and Lin and Liu
\cite{2,4,3} made some important progress with the existence
of global weak solutions and the partial regularity of the
incompressible hydrodynamic flow equations of liquid crystals. The Ericksen-Leslie system is a macroscopic
continuum description for the time evolution of the materials under the influence of both the flow
velocity field $u$ and the direction field $n$ of
rod-like liquid crystals. It can  also be viewed as a multi-scales coupling system which describes the interaction between microstructure
(kinetic theory) and macrostructure (fluid), where the direction field $n$ is molecule-scale and the density, the fluid velocity field $u$ and the pressure are macroscopic.

Now we introduce a formal physical derivation of the compressible models for nematic liquid crystal flow through an appropriate energetic variational approach,
in which the least action principle gives the Hamiltonian parts (conservative force) of the hydrodynamic systems while the
maximum/minimum dissipation principle, namely, Onsager's principle, gives the dissipative parts
(dissipative force) of the systems. We refer the readers to \cite{sun,wxl} for the derivation of the incompressible nematic liquid crystal model
(a proper penalty approximation).

In the content of hydrodynamics, the basic variable is the flow map
(particle trajectory) $x(X,t)$, which is defined by the following ordinary
differential equation: \beq
    \label{l1}
    x_t(X,t)=u(x(X,t),t),\quad x(X,0)=X,
\eeq
where $X$ is the Lagrangian coordinate and $x$ is the Eulerian coordinate.

The deformation tensor $F(X,t)$ is defined as
\beq
    \label{l2}
    F(X,t)=\frac{\partial x}{\partial X}.
\eeq
Applying the chain rule, we see that $F(x,t)$ satisfies the following transport equation:
\beq
    \label{l3}
    F_t+u\cdot\nabla F=\frac{{\rm d}}{{\rm d}t}F=\frac{\partial x_t}{\partial X}=\frac{\partial u}{\partial X}=\frac{\partial u}{\partial x}\cdot\frac{\partial x}{\partial X}=\nabla uF,
\eeq
which stands for $F_{ij,t}+u_k\nabla_kF_{ij}=\nabla_ku_iF_{kj}$.

Define the density as
\beq
    \label{l4}
    \rho(x,t)=\frac{\rho_0(X)}{\det F}.
\eeq
By the identity of the variation for the determinant of a tensor
\beq
    \label{l5}
    \delta \det F=\det F {\rm tr}(F^{-1}\delta F),
\eeq
we have
\beq
    \label{l6}
    \rho_t+u\cdot\nabla\rho=\frac{{\rm d}}{{\rm d}t}\left(\frac{\rho_0(X)}{\det F}\right)=-\frac{\rho_0(X)}{(\det F)^2}\det F{\rm tr}\left(F^{-1}
    \frac{{\rm d}}{{\rm d}t}F\right)=-\rho\nabla\cdot u,
\eeq
and then we get the transport equation
\beq
    \label{l7}
    \rho_t+u\cdot\nabla\rho+\rho\nabla\cdot u=0,
\eeq
which can also be derived by the conservation law of mass.

The kinematic assumption (without dissipations) about the case of small molecules implies that $n$ is just transported by the flow trajectory, i.e.,
\beq
    \label{l11}
    \frac{D}{Dt}n=n_t+u\cdot\nabla n=0,
\eeq
and then we have
\beq
    \label{l12}
    n(x(X,t),t)=n_0(X),
\eeq
which also implies that the center of gravity of the molecules moves along the streamline of the velocity.
This is due to the fact that if the size of the molecules is small enough, then the directors are not affected by the stretching of the fluid (see \cite{zlz}).

We remark that for the case of big molecule (nematic), the transport of the orientation vector $n$ is governed by
\beq
    \frac{D}{Dt}n=n_t+(u\cdot\nabla) n+D_\beta (u)n,
\eeq
where $D_\beta(u)=-\frac{\nabla u-(\nabla u)^T}{2}-(-2\beta-1)\frac{\nabla u+(\nabla u)^T}{2}$, which is the stretching of the fluid on $n$
with $\beta$ depending only on the shape of the molecules (see \cite{zhang}).

The momentum equation of motion for hydrodynamic flow of liquid
crystals can be derived from the least action principle (Hamiltonian
principle). The action functional takes the form:
\beq
    \label{l8}
    \mathcal{A}(x)=\int_0^t\int_{\Omega_0}\left(\frac12\rho_0(X)|x_t|^2-W^\lambda(F)\right){\rm d}X{\rm d}t,\quad
\lambda>0,
\eeq
where $\Omega_0$ is the original domain occupied by
the material. The first part of $\mathcal{A}(x)$ denotes the kinetic
energy and the second one denotes the elastic energy.

By the definition of the density in (\ref{l4}), we have
\beq
    \label{l9}
    \int_0^t\int_{\Omega_0}\rho_0(X)|x_t|^2{\rm d}X{\rm d}t=\int_0^t\int_{\Omega_0}\frac{\rho_0(X)}{\det F}|x_t|^2\det F{\rm d}X{\rm d}t
    =\int_0^t\int_{\Omega}\rho(x,t)|u(x,t)|^2{\rm d}x{\rm d}t.
\eeq

In this paper, we consider the isotropic energy function $W^\lambda(F)$ of the form
\beq
    \label{l10}
    W^\lambda(F)=\left(\lambda^2\omega\left(\frac{\rho_0(X)}{\det F}\right)+\frac{\nu}{2}\left|F^{-1}\nabla_Xn_0(X)\right|^2\right)\det F,
\eeq
where $\omega$ is a $C^\infty$ function and denotes the energy density, $n_0(X):\Omega_0\rightarrow S^2$
is a unit-vector field which represents the molecular orientation of
the liquid crystal material. The first term on the
R.H.S. of (\ref{l10}) should be regarded as a penalization term which
drives the motion toward incompressibility in the limit as the
parameter $\lambda$ becomes large. And the second term on the R.H.S. of (\ref{l10}) $\frac{\nu}{2}\left|F^{-1}\nabla_Xn_0(X)\right|^2$
is equal to $\frac{\nu}{2}|\nabla n|^2$ in Eulerian coordinate by the kinematic assumption, where $\nu>0$ is viscosity of the fluid and denotes
the microscopic elastic relaxation time.

For any $y=(y_1,y_2,\ldots,y_N)\in C_c^1(\Omega\times[0,+\infty);\mathbb{R}^N)$,
let $x^\epsilon=x+\epsilon y(x)$,
$F^\epsilon=\frac{\partial x^\epsilon}{\partial X}$ and
$\rho^\epsilon=\rho(x^\epsilon(X,t),t)$. Then applying
$\delta=\frac{{\rm d}}{{\rm d}\epsilon}\big|_{\epsilon=0}$ to (\ref{l8}), we get
\beq
    \label{l13}
    0=\delta\mathcal{A}(x^\epsilon)
    =\int_0^t\int_{\Omega_0}\left(\rho_0(X)(x_t,y_t)-\delta
    W^\lambda(F)\right){\rm d}X{\rm d}t,
\eeq where $(f,g)$ is the inner product $f\cdot g$ for some $f,g:\Omega\rightarrow \mathbb{R}^N$.

Using the definition of $\rho(x,t)$, we have
\beq
    \label{l14}
    \nonumber\int_0^t\int_{\Omega_0}\rho_0(X)(x_t,y_t){\rm d}X{\rm d}t
    &=&-\int_0^t\int_{\Omega_0}\rho_0(X)(x_{tt},y){\rm d}X{\rm d}t=-\int_0^t\int_{\Omega_0}\frac{\rho_0(X)}{\det F}(x_{tt},y){\rm d}x{\rm d}t
    \\&=&
    -\int_0^t\int_{\Omega}\rho(x,t)(u_t+u\cdot\nabla u,y){\rm d}x{\rm d}t.
\eeq

Now we turn to show the calculation of the second term on the R.H.S of $(\ref{l13})$. On one hand, noting the facts that
\beq
    \label{l15}
    \nonumber
    \delta\rho^\epsilon&=&\delta\left(\frac{\rho_0(X)}{\det F^\epsilon}\right)
    =-\frac{\rho_0(X)}{(\det F)^2}\det F {\rm tr}\left(F^{-1}\delta F^\epsilon\right)
    \\&=&
    -\frac{\rho_0(X)}{\det F}{\rm tr}\left(\frac{\partial X}{\partial x}\frac{\partial y}{\partial X}\right)
    =-\frac{\rho_0(X)}{\det F}\nabla\cdot y,
\eeq
and
\beq
    \label{l16}
    \delta\left(\det F^\epsilon\right)=\det F{\rm tr}\left(F^{-1}\nabla_Xy\right)=\det F\nabla\cdot y,
\eeq
and then we use the definition of $\rho$ to give
\beq
    \label{l17}
    \nonumber &&-\lambda^2\int_0^t\int_{\Omega_0}\delta\left(\omega\left(\frac{\rho_0(X)}{\det F^\epsilon}\right)\det F^\epsilon\right){\rm d}X{\rm d}t
    \\&=& \nonumber
    -\lambda^2\int_0^t\int_{\Omega_0}\left(-\omega'\left(\frac{\rho_0(X)}{\det F}\right)\rho_0(X)\nabla\cdot y+\omega\left(\frac{\rho_0(X)}{\det F}\right)\det F\nabla\cdot y\right){\rm d}X{\rm d}t
    \\&=& \nonumber
    -\lambda^2\int_0^t\int_{\Omega}\left(-\omega'\left(\frac{\rho_0(X)}{\det F}\right)\rho_0(X)\nabla\cdot y+\omega\left(\frac{\rho_0(X)}{\det F}\right)\det F\nabla\cdot y\right)\frac{1}{\det F}{\rm d}x{\rm d}t
    \\&=& \nonumber
    -\lambda^2\int_0^t\int_{\Omega} \left(\omega(\rho)-\rho \omega^\prime(\rho)\right)\nabla\cdot y{\rm d}x{\rm d}t
    \\&=&
    -\lambda^2\int_0^t\int_{\Omega} (\nabla P(\rho), y){\rm d}x{\rm d}t,
\eeq
where $P(\rho)$ is the pressure function given by
\beq
    \label{ll8}
    P(\rho)=-\omega(\rho)+\rho\omega^\prime(\rho),
\eeq
which can also be obtained by the first law of thermodynamics. Here we show the proof briefly for completeness.

According to the first law of thermodynamics, we obtain
\beq
    \label{l19}
    {\rm d}W=-P{\rm d}V-S{\rm d}T,
\eeq where $W$, $V$, $S$ and $T$ denote the energy, volume, entropy
and temperature, respectively. Moreover,
\beq
    \label{l20}
    W=\omega V,\quad V=\frac{m}{\rho},\quad \omega=\omega(\rho),
\eeq
where $m$ is the total mass.

From (\ref{l19}), we have
\beq
    \label{l21}
    -P=\frac{\partial W}{\partial V}=\frac{\partial(\omega V)}{\partial V}=\omega+V\frac{\partial\omega}{\partial V}
    =\omega+V\frac{\partial\omega\left(\frac{m}{V}\right)}{\partial V}=\omega-\omega^\prime\frac{m}{V}=\omega-\rho\omega^\prime,
\eeq and then we get (\ref{ll8}).

On the other hand, by the Leibniz's formula, we have
\beq
    \label{l22}
    &&\nonumber
    -\frac{\nu}{2}\int_0^t\int_{\Omega_0}\delta\left(\left|(F^\epsilon)^{-1}\nabla_Xn_0(X)\right|^2\det F^\epsilon\right){\rm d}X{\rm d}t
    \\&=& \nonumber
    -\frac{\nu}{2}\int_0^t\int_{\Omega_0}|F^{-1}\nabla_Xn_0(X)|^2\delta\left(\det F^\epsilon\right){\rm d}X{\rm d}t
    \\&& \nonumber
    -\nu\int_0^t\int_{\Omega_0}\left(F^{-1}\nabla_Xn_0(X),\delta(F^\epsilon)^{-1}\nabla_Xn_0(X)\right)
    \det F{\rm d}X{\rm d}t,
    \\&=&
    V_1+V_2,
\eeq
where
\beq
    \label{l23}
    \nonumber V_1&=&-\frac{\nu}{2}\int_0^t\int_{\Omega_0}\left|F^{-1}\nabla_Xn_0(X)\right|^2\det F\nabla\cdot y{\rm d}X{\rm d}t
    \\&=&
    -\frac{\nu}{2}\int_0^t\int_{\Omega}|\nabla n|^2\nabla\cdot y{\rm d}x{\rm d}t
    =\frac{\nu}{2}\int_0^t\int_{\Omega}\left(\nabla\left(|\nabla n|^2\right),y\right){\rm d}x{\rm d}t,
\eeq
and
\beq
    \label{l24}
    \nonumber V_2&=&
    -\nu\int_0^t\int_{\Omega_0}
    \left(\partial_{X_j}n_0(X)\frac{\partial X_j}{\partial x_i},\partial_{X_r}n_0(X)\delta
    \left(\frac{\partial X_r}{\partial x^\epsilon_i}\right)\right)\det F{\rm d}X{\rm d}t
    \\&=& \nonumber
    -\nu\int_0^t\int_{\Omega_0}\left(\partial_{X_j}n_0(X)\frac{\partial X_j}{\partial x_i},\partial_{X_r}n_0(X)\left(-\frac{\partial X_r}{\partial x_k^\epsilon}\frac{{\rm d}}{{\rm d}\epsilon}\frac{\partial x_k^\epsilon}{\partial X_l}\frac{\partial X_l}{\partial x_i^\epsilon}\right)\bigg|_{\epsilon=0}\right)\det F{\rm d}X{\rm d}t
    \\&=& \nonumber
    \nu\int_0^t\int_{\Omega_0}\left(\partial_{X_j}n_0(X)\frac{\partial X_j}{\partial x_i},\partial_{X_r}n_0(X)\frac{\partial X_r}{\partial x_k}\frac{\partial y_k}{\partial X_l}\frac{\partial X_l}{\partial x_i}\right)\det F{\rm d}X{\rm d}t
    \\&=& \nonumber
    \nu\int_0^t\int_{\Omega}\left(\partial_{x_i}n,\partial_{x_k}n\frac{\partial y_k}{\partial x_i}\right){\rm d}x{\rm d}t
    \\&=&
    -\nu\int_0^t\int_{\Omega}\left(\nabla\cdot\left(\nabla n\otimes\nabla n\right),y\right){\rm d}x{\rm d}t,
\eeq
where we have taken the differential to $F^\epsilon
(F^\epsilon)^{-1}=I_N$ with respect to $\epsilon$ and then obtained
\beq
    \label{l25}
    \delta\left(\frac{\partial X_r}{\partial x^\epsilon_i}\right)=\left(-\frac{\partial X_r}{\partial x_k^\epsilon}\frac{{\rm d}}{{\rm d}\epsilon}\frac{\partial x_k^\epsilon}{\partial X_l}\frac{\partial X_l}{\partial x_i^\epsilon}\right)\bigg|_{\epsilon=0}.
\eeq

Particularly, in the case of the incompressible system, the
calculation above can be simplified by the incompressibility, i.e.,
$\det F=1$.

In summary,  we have derived from $\delta\mathcal{A}(x^\epsilon)=0$ that
\beq
    \label{nondissipation}
    \rho u_t+\rho u\cdot\nabla u+\lambda^2\nabla P(\rho)-\nu\left(\nabla n\otimes\nabla n-\frac{|\nabla n|^2}{2}I_N\right)=0.
\eeq

Finally, we are to deal with the dissipation. Firstly, take the internal dissipation into account together with the transport equation (\ref{l11}). For any $C^1$
map $ n:\Omega\times[0,+\infty)\rightarrow S^2$, let $\displaystyle
I(n)=\frac{\theta}{2}\int_{\Omega}|\nabla n|^2{\rm d}x$ be the energy
of $n$, where $\theta>0$ is the viscosity constant. The critical
points of the energy $I(n)$ are called the harmonic maps. The motion of $n$ satisfies the gradient flow (subject to the restraint $|n|=1$) that
\beq
    \label{D_t n}
    \frac{D}{Dt}n=-\frac{\delta I(n)}{\delta n},
\eeq
where $\frac{\delta I(n)}{\delta n}$ is the Frechet derivative.

Applying $\widetilde{\delta}=\frac{{\rm d}}{{\rm d}\tau}\big|_{\tau=0}$ to $I(m(\tau))$, where
\beq
    \label{l31}
    m(\tau)=\frac{n+\tau\varphi}{|n+\tau\varphi|},
\eeq
we have
\beq
    \label{l32}
    \frac{\delta I(n)}{\delta n}=-\theta(\Delta n+|\nabla n|^2n),
\eeq
which, combined with (\ref{D_t n}) (for the case of small molecule), yields the equation for the orientation field
\beq
    \label{n_t}
    n_t+(u\cdot\nabla)n=\theta(\Delta n+|\nabla n|^2n).
\eeq

Secondly, for the nematic liquid crystals with small rod-like molecules,
we have the following energy law (see \cite{hkl,zhang,Onsager1,Onsager2,Onsager3,wxl} for more details)
\beq
    \label{l26}
    \frac{{\rm d}}{{\rm d}t}E^{total}=-\triangle,
\eeq
where $E^{total}$ is the summation of kinetic energy and the elastic energy shown as follows
\beq
    \label{l27}
    E^{total}=\int_{\Omega}\left(\frac12\rho|u|^2+\lambda^2\omega(\rho)+\frac{\nu}{2}|\nabla n|^2\right){\rm d}x,
\eeq
and $\triangle$ denotes dissipation being a linear combination of the squares of various rate functions such as velocity,
rate of strain or the material derivative of internal variables, i.e.,
\beq
    \label{l28}
    \triangle=\int_{\Omega}\left(2\mu|\mathcal{D}u|^2+\kappa|\nabla\cdot u|^2+\frac{\nu}{\theta}|n_t+u\cdot\nabla n|^2\right){\rm d}x,
\eeq
where $\mathcal{D}u=\frac{\nabla u+\nabla^T u}{2}=\left(\frac{\nabla_iu_j+\nabla_ju_i}{2}\right)_{N\times N}$ is called the deformation tensor,
$\mu$ and $\kappa$ are the shear viscosity and the bulk viscosity coefficients satisfying $\mu>0$ and $2\mu+N\kappa\geq 0$.
The last term on the R.H.S. of (\ref{l28}) represents for the molecule-level microscopic dissipation.

Let $u^\tau=u+\tau
v$, where $v\in C_c^1(\Omega\times[0,+\infty);\mathbb{R}^N)$. Applying $\widetilde{\delta}=\frac{{\rm d}}{{\rm d}\tau}\big|_{\tau=0}$ to
$\triangle$, and integrating by parts, we have from (\ref{n_t}) and the fact $|n|=1$ that
\beq
    \label{l29}
    \nonumber
    \widetilde{\delta}\triangle&=&\nonumber
    \mu\int_{\Omega}\left(\nabla u+\nabla^T u,\nabla v+\nabla^T v\right){\rm d}x+2\kappa\int_{\Omega}(\nabla\cdot u,\nabla\cdot v){\rm d}x
    \\&& \nonumber
    +\frac{2\nu}{\theta}\int_\Omega (n_t+u\cdot\nabla n, v\cdot\nabla n){\rm d}x
    \\&=&\nonumber
    \mu\int_{\Omega}\left((\nabla u,\nabla v)+\left(\nabla^T u,\nabla^T v\right)\right){\rm d}x+\mu\int_{\Omega}\left(\left(\nabla u,\nabla^T v\right)+\left(\nabla^T u,\nabla v\right)\right){\rm d}x
    \\&& \nonumber
    +2\kappa\int_{\Omega}(\nabla\cdot u,\nabla\cdot v){\rm d}x+2\nu\int_\Omega(\Delta n+|\nabla n|^2n,v\cdot\nabla n){\rm d}x
    \\&=&
    2\int_{\Omega}\left(-\mu\Delta u-(\kappa+\mu)\nabla(\nabla\cdot u)+\nu\Delta n\cdot\nabla n,v\right){\rm d}x=0,
\eeq
where $\Delta n\cdot\nabla n=\sum_{i=1}^N\Delta n_i\nabla n_i$, implying that
\beq
    \label{dissipation}
    \mu\Delta u+(\kappa+\mu)\nabla(\nabla\cdot u)-\nu\Delta n\cdot\nabla n=0.
\eeq

The interesting fact from the above derivation is that the induced stress term $\nu\Delta n\cdot\nabla n=\nu\nabla\cdot
\left(\nabla n\otimes\nabla n-\frac{|\nabla n|^2}{2}I_N\right)$ can be derived either by the least action principle (contained in (\ref{nondissipation}))
or by the maximum dissipation principle (contained in (\ref{dissipation})).
Therefore, they can either be recognized as conservative or dissipative (see also \cite{hkl,wxl}). Clearly, if the system is
incompressible, the term $(\kappa+\mu)\nabla(\nabla\cdot u)$ in (\ref{dissipation}) should be zero.

Putting (\ref{nondissipation}) and (\ref{dissipation}) together, we have the momentum equation
\beq
    \label{l30}
    &&\nonumber\rho u_t+\rho (u\cdot\nabla)u+\lambda^2\nabla \left(P(\rho)\right)
    \\&=&\mu\Delta u+(\kappa+\mu)\nabla(\nabla\cdot u)-\nu\nabla\cdot\left(\nabla n\otimes\nabla n-\frac{|\nabla n|^2}{2}I_N\right).
\eeq

In conclusion, the compressible hydrodynamic flow equations of liquid crystals can be written as follows (the functions and the viscosity constants should depend on the value of the parameter $\lambda$):
\beq
\label{l38}
\begin{cases}
      \rho^\lambda_t+\nabla\cdot(\rho^\lambda u^\lambda)=0,\\
      (\rho^\lambda u^\lambda)_t+\nabla\cdot(\rho^\lambda u^\lambda\otimes{u^\lambda})+\lambda^2\nabla(P(\rho^\lambda))
      \\ \quad
      =\mu^\lambda\Delta{u^\lambda}+(\kappa^\lambda+\mu^\lambda)\nabla(\nabla\cdot u^\lambda)-\nu^\lambda\nabla\cdot\left(\nabla{n^\lambda}\odot{\nabla{n^\lambda}}-\frac{|\nabla n^\lambda|^2}{2}I_N\right),\\        {n^\lambda_t}+({u^\lambda}\cdot\nabla){n^\lambda}=\theta^\lambda(\Delta{n^\lambda}+|\nabla n^\lambda|^2{n^\lambda}).%\\[2mm]
\end{cases}\eeq
Here we consider $x\in T^N$, a torus in $\mathbb{R}^N$, $N=2$ or 3 and
$t>0$. The coefficients satisfy that $\mu^\lambda,\nu^\lambda,\theta^\lambda>0$,
and $2\mu^\lambda+N\kappa^\lambda\geq0$. The unknowns are the
density $\rho^\lambda:{T^N}\times[0,+\infty)\rightarrow \mathbb{R}^1$, the velocity field $u^\lambda:{T^N}\times[0,+\infty)\rightarrow \mathbb{R}^N$, and the molecular orientation field of the liquid crystal material
$n^\lambda: {T^N}\times[0,+\infty)\rightarrow S^2$, which
is a unit vector. $P(\rho)$ is the
smooth pressure-density function with $P'(\rho)>0$ for $\rho>0$. The symbol $\otimes$ represents the usual Kronecker multiplication, e.g., $u\otimes u=(u_i u_j)_{1\le i, j\le N}$, and $\nabla n\odot\nabla n$ represents the $N\times N$ matrix whose $(i,j)$-th entry is given by $\nabla_in\cdot \nabla_jn$ for $1\le i,j\le N$.

From mathematical point of view, it is reasonable to expect
that, as $\rho^\lambda\rightarrow 1$, the first equation in
(\ref{l38}) yields the incompressible condition $\nabla\cdot u=0$.
Suppose that the limits $u^\lambda\rightarrow u$ and
$n^\lambda\rightarrow n$ exist as $\lambda\rightarrow\infty$, and the viscosity coefficients satisfy that
\beq
    \label{l39}
    \mu^\lambda\rightarrow \mu>0,\ \ \kappa^\lambda\rightarrow \kappa,\ \ \nu^\lambda\rightarrow \nu>0,\ \ \theta^\lambda\rightarrow \theta>0,\ \ {\rm as} \ \ \lambda\rightarrow\infty.
\eeq Then (at least formally) we obtain
the following incompressible model by taking $\lambda\rightarrow
\infty$ that \beq \label{l40}
\begin{cases}\nabla\cdot u=0,\\
         u_t+(u\cdot \nabla) u+\nabla p=\mu
         \Delta{u}-
\nu\nabla\cdot(\nabla n\odot\nabla n),\\
        {n_t}+({u}\cdot\nabla){n}=\theta(\Delta{n}+|\nabla n|^2{n}),%\\[2mm]
\end{cases}\eeq
where $\nabla p$ is the limit of the term $\frac{\lambda^2}{\rho^\lambda}\nabla\left(P(\rho^\lambda)\right)-\frac{\nu^\lambda}{\rho^\lambda}\nabla\left(\frac{|\nabla n^\lambda|^2}{2}\right)$.
For simplicity, we assume in the following that $\kappa^\lambda\equiv\kappa$, $\mu^\lambda\equiv\mu$, $\mu^\lambda\equiv\mu$ and $\theta^\lambda\equiv\theta$
are constants independent of $\lambda$, which satisfy that $\mu,\nu,\theta>0$ and  $2\mu+N\kappa\geq0$.

In a series of papers, Lin \cite{1} and Lin and Liu \cite{2,4,3} addressed the existence and
partial regularity theory of suitable weak solution to the incompressible hydrodynamic flow of liquid crystals (\ref{l40}) of variable
length. More precisely, they considered the approximate equation of incompressible hydrodynamic flow of liquid
crystals $\big($$|\nabla n|^2n$ is replaced by $\frac{(1-|n|^2){n}}{\epsilon^2}$$\big)$, and proved \cite{2}, the local existence
of classical solutions and the global existence of weak solutions in dimension two and three.
For any fixed ${\epsilon}>0$, they also showed the existence and uniqueness of global classical solutions either in dimension two or
dimension three when the fluid viscosity $\mu$ is sufficiently large; Lin and Liu \cite{4} extended the classical theorem by
Caffarelli et al. \cite{cnk} on the Navier-Stokes equations that asserts the one dimensional parabolic Hausdorff
measure of the singular set of any {\it suitable} weak solution is
zero. See also \cite{Blanca,liu-walkington,Simon} for relevant results.  For the system (\ref{l40}), Lin et al. \cite{llw} proved that there exists
a global weak solution which is regular with the exception of at most finitely many time-slices in dimension two. For the density-dependent incompressible
flow of incompressible liquid crystals, on one hand,
Liu and Zhang \cite{liu-zhang} obtained the global weak solutions in dimension three with the initial density $\rho_0\in L^2$ for the model of variable length.
Jiang and Tan \cite{jiang-tan} improved the condition of $\rho_0$ to $\rho_0\in L^\gamma$, $\gamma>\frac32$. However, the estimates depend on $\epsilon$,
and thus they can not take the limit $\epsilon\rightarrow0$. On the other hand, considering the original term $|\nabla n|^2n$, Wen and Ding \cite{Wen-Ding}
proved the local existence and uniqueness of the strong solutions to the model in
a bounded domain in $\mathbb{R}^N$ ($N=2$ or $3$), provided that the initial density $\rho_0\geq 0$. Furthermore, they got the global existence and uniqueness
of the strong solutions with small initial data and $\inf\limits_{x\in\Omega}\rho_0>0$ in dimension two. Very recently, Li and Wang \cite{Li-Wang} proved the existence
and uniqueness of the local strong solutions with large initial data and the global strong solutions with small data in Besov space for the initial density away from
vacuum in dimension three. For an incompressible non-isothermal model, we refer to \cite{feireisl}.

The study for the compressible hydrodynamic flow (\ref{l38}) began in recent years.
Concerning the Dirichlet and Neumann boundary condition
for $(u,n)$ in dimension one, Ding et al. \cite{w1} obtained the existence and uniqueness of global classical solutions with the initial data
$(\rho_0,u_0,n_0)\in C^{1,\alpha}(I)\times C^{2,\alpha}(I)\times C^{2,\alpha}(I)$ and $\rho_0\geq c_0>0$, where $I=[0,1]$. They also addressed both the existence and
the uniqueness of global strong solutions for $0\leq\rho_0\in H^1(I)$ and $(u_0,n_0)\in H^1(I)\times H^2(I)$. Furthermore, Ding et al. \cite{w2} established the existence
of weak solution $(\rho,u,n)$ in dimension one with $0\leq\rho_0\in L^\gamma(I)$ for $\gamma>1$, $u_0\in L^2(I)$ and $n_0\in H^1(I)$. Huang et al. \cite{Huang-Wang-Wen1} studied the three dimensional Cauchy problem and other boundary problems. They obtained the local existence of unique strong solution
provided that the initial data $\rho_0,u_0,n_0$ are sufficiently regular and satisfy a natural compatibility condition, and then they proved a criterion for possible
breakdown of such a local strong solution at finite time in terms of the blow up of the quantities $\|\rho\|_{L_t^\infty L_x^\infty}$ and $\|\nabla n\|_{L_t^3L_x^\infty}$.
They obtained in \cite{Huang-Wang-Wen2} a blow up criterion in terms of $\|\mathcal{D}u\|_{L_t^1 L_x^\infty}$ and $\|\nabla n\|_{L_t^2L_x^\infty}$. Recently, Li et al. \cite{Lijing} established the global existence of classical solutions to the Cauchy problem with smooth initial data which are of small energy but possibly large oscillations with constant state as farfield condition which could be either vacuum or non-vacuum in dimension three. For the hydrodynamic
flow of liquid crystals of variable length, we refer to \cite{liu-liu-hao1,liu-liu-hao2,Wang-Yu2012}.

It is interesting to ask whether the solutions for the compressible flow of liquid crystals can converge to the solutions for the incompressible system.
For the flow of liquid crystals of variable length, this problem has been answered by recent papers \cite{liu-hao,Wang-Yu}. The present paper focuses on the system (\ref{l38})
and (\ref{l40}) and intends to answer such a problem.

This paper is organized as follows. In Section $2$, we state the main results. In Section $3$, we prove that the local strong solutions of (\ref{l38})
exist for sufficiently small disturbances from the general incompressible initial data, meanwhile, we get the uniform stability of the local solution family which
yields a lifespan of the system (\ref{l38}). In view of this fact, we show that the
local solutions for (\ref{l38}) converge to a local solution for the
limiting incompressible system (\ref{l40}) by means of compactness arguments. In Section $4$, the global existence of the strong solutions to the incompressible
system (\ref{l40}) is derived provided that the
initial data of the incompressible model are sufficiently small. In Section $5$, we obtain the
convergence rates about
$(\rho^\lambda,u^\lambda,n^\lambda)\rightarrow (1,u,n)$ in some
sense when $\lambda\rightarrow \infty$.   These results depend on
the techniques modified from \cite{Jiang 1,Klainerman,Lei,Lei2,Sideris}.

\setcounter{section}{1} \setcounter{equation}{0}
\section{Notations and Statements of Main Results}
First of all, we can rewrite (\ref{l38}) when the density is away from vacuum as follows:
\beq \label{l42}
\begin{cases}
      \rho^\lambda_t+\nabla\cdot(\rho^\lambda u^\lambda)=0,\\
      u^\lambda_t+(u^\lambda\cdot \nabla) u^\lambda+\frac{\lambda^2}{\rho^\lambda}\nabla P(\rho^\lambda)=\frac{\mu}{\rho^\lambda}\Delta{u^\lambda}+\frac{(\kappa+\mu)}{\rho^\lambda}\nabla(\nabla\cdot u^\lambda)-\frac{\nu}{\rho^\lambda}\Delta n^\lambda\cdot\nabla n^\lambda,\\        {n^\lambda_t}+({u^\lambda}\cdot\nabla){n^\lambda}=\theta(\Delta{n^\lambda}+|\nabla n^\lambda|^2{n^\lambda}).%\\[2mm]
\end{cases}\eeq

Throughout this paper, for convenience, let $\int_{T^N} f(x)=\int_{T^N} f(x) {\rm d}x$ and $\int_0^t g(s)=\int_0^t g(s){\rm d}s$, denote by $\|\cdot\|$,
$\|\cdot\|_s$ and $\|\cdot\|_\infty$ the norms in $L^2(T^N)$,
$H^s(T^N)$ and $L^\infty({T^N})$ respectively, and especially denote by $\|\cdot\|_{L^4}$ the norm in $L^4({T^N})$.

Define
\beq
\label{l43}
\begin{cases}
      E_s(U(t))=\frac12\sum\limits_{|\alpha|\leq s}\displaystyle\int_{T^N}\left(\lambda^2|\nabla^\alpha(\rho-1)|^2+|\nabla^\alpha u|^2+|\nabla\nabla^\alpha n|^2\right), \\
      \widetilde{E}_s(U(t))=\frac12\sum\limits_{|\alpha|\leq s}\displaystyle\int_{T^N}\left(\lambda^2\frac{P'(\rho)}{\rho}|\nabla^\alpha(\rho-1)|^2+\rho|\nabla^\alpha u|^2+|\nabla\nabla^\alpha n|^2\right),
\end{cases}\eeq where $U=(\rho,u,n)$.

It is easily seen that
\beq
    \label{l45}
    E_s(U(t))\sim\widetilde{E}_s(U(t)),
\eeq provided that $|\rho-1|$ is sufficiently
small (see \cite{Lei,Lei2,Sideris}).

Now we state the main results of this paper.
\begin{Theorem}\label{th.1.1} Consider the compressible model (\ref{l38})
with the following initial data
\beq
   \label{l46}
   \rho^\lambda(x,0)=1+\overline{\rho}_0^\lambda(x),\ \
   u^\lambda(x,0)=u_0(x)+\overline{u}_0^\lambda(x), \ \
   n^\lambda(x,0)=\frac{n_0(x)+\overline{n}_0^\lambda(x)}{\left|n_0(x)+\overline{n}_0^\lambda(x)\right|},
\eeq
where $u_0$, $n_0$ satisfy
\beq
    \label{l47}
    u_0(x)\in H^{s+1}({T^N}), \ \ \nabla\cdot u_0=0,\ \ n_0\in H^{s+2}({T^N},S^2)
\eeq
for any $s\geq \left[\frac{N}{2}\right]+2$. Moreover, for small positive constant $\delta_0$,
the functions $\overline{\rho}_0^\lambda(x)$, $\overline{u}_0^\lambda(x)$,
$\overline{n}_0^\lambda(x)$ are assumed to satisfy
\beq
    \label{l48}
    \|\overline{\rho}_0^\lambda\|_s\leq \lambda^{-2}\delta_0, \ \
    \|\overline{u}_0^\lambda(x)\|_{s+1}\leq \lambda^{-1} \delta_0,\ \
    \|\nabla n^\lambda(x,0)-\nabla n_0(x)\|_{s+1}\leq \lambda^{-1}\delta_0.
\eeq
Then the following statements hold.

Uniform stability: There exist constants $T_0$ and $C$
independent of $\lambda$ such that the unique strong solution
$(\rho^\lambda, u^\lambda, n^\lambda)$ of system (\ref{l38}) exists
for all large $\lambda$ on the time interval $[0,T_0]$ with properties:
\beq
\label{l49}
\begin{cases}
      E_s(U^\lambda(t))+\displaystyle\int_0^t\left(\mu\|\nabla u^\lambda\|_s^2+(\kappa+\mu)\|\nabla\cdot u^\lambda\|_s^2+\theta \|\nabla n^\lambda\|_{s+1}^2\right)\leq C,
      \\E_{s-1}(\partial_tU^\lambda(t))+\|n_t^\lambda\|^2+\displaystyle\int_0^t\left(\mu\|\nabla \partial_t             u^\lambda\|_{s-1}^2+(\kappa+\mu)\|\nabla\cdot \partial_tu^\lambda\|_{s-1}^2+\theta \|\nabla \partial_tn^\lambda\|_{s}^2\right)\leq C,
      \\ |n^\lambda|=1 \ \ {\rm in} \ \ Q_{T_0}=\overline{T^N}\times[0,T_0],%\\[2mm]
\end{cases}\eeq
where $E_{s-1}(\partial_tU(t)) =\frac12\sum\limits_{|\beta|\leq
s-1}\int_{T^N}\left(\lambda^2|\nabla^\beta\partial_t\rho|^2+|\nabla^\beta\partial_t
u|^2+|\nabla^\beta\partial_t\nabla n|^2\right)$.

Local existence of solutions for incompressible system: There exist
functions $u$ and $n$ such that \beq \label{l50}
\begin{cases}\rho^\lambda\rightarrow 1 \quad  {\rm in}\ \
             L^\infty([0,T_0];H^s)\cap {\rm Lip}([0,T_0];H^{s-1}),
             \\u^\lambda\rightharpoonup u  \quad {\rm weakly}^* \ {\rm in}\ \
             L^\infty([0,T_0];H^s)\cap {\rm Lip}([0,T_0];H^{s-1}),
             \\u^\lambda\rightarrow u \quad {\rm in}\ \
             C([0,T_0];H^{s'}),
             \\n^\lambda\rightharpoonup n  \quad {\rm weakly}^* \ {\rm in}\ \
             L^\infty([0,T_0];H^{s+1})\cap {\rm Lip}([0,T_0];H^{s}),
             \\n^\lambda\rightarrow n \quad {\rm in}\ \
             C([0,T_0];H^{s'+1})
             %\\[2mm]
\end{cases}\eeq for any $s'\in[0,s)$,
and the function pair $(u,n)$ is the unique strong solution of the
incompressible system of liquid crystal (\ref{l40}) with the initial
data \beq
    \label{l51}
    u(x,0)=u_0(x),\ \ n(x,0)=n_0(x),
\eeq
for some $p\in L^\infty([0,T_0];H^{s-1})\cap L^2([0,T_0];H^{s})$.

\end{Theorem}

\begin{Remark} By scaling we know that $\lambda^{-1}$ is the Mach number. Note that as $\lambda\rightarrow 0$, the initial data may be very large.
Theorem \ref{th.1.1} states that although the initial data depends on $\lambda$, the local existence time $T_0$ is independent of $\lambda$.
\end{Remark}

\begin{Theorem}\label{th.1.2} Consider the strong solutions $(\rho^\lambda,u^\lambda,n^\lambda)$ of system (\ref{l38}) obtained in Theorem \ref{th.1.1}.
Suppose in addition that the initial data satisfies
\beq
    \label{l52}
    \|u_0\|_s^2+\|\nabla n_0\|_{s}^2\leq \varepsilon_0,
\eeq
where $\varepsilon_0$ is a positive constant. If $\varepsilon_0$ is sufficiently small, then for any fixed $T>0$, the strong solution $(\rho^\lambda,u^\lambda,n^\lambda)$
satisfies the following estimates:
\beq
      \label{l53}
      E_s(U^\lambda(t))+\displaystyle\int_0^t\left(\mu\|\nabla u^\lambda\|_s^2+(\kappa+\mu)\|\nabla\cdot u^\lambda\|_s^2+\theta \|\nabla^2 n^\lambda\|_{s}^2\right)\leq 4(\varepsilon_0+\lambda^{-2}\delta_0^2)
\eeq
for $t\in [0,T^\lambda)$, and
\beq
     \label{l54}
     &&\nonumber
     E_{s-1}(\partial_tU^\lambda(t))+\|n^\lambda_t\|^2+\displaystyle\int_0^t\left(\mu\|\nabla \partial_t             u^\lambda\|_{s-1}^2+(\kappa+\mu)\|\nabla\cdot \partial_tu^\lambda\|_{s-1}^2+\theta \|\nabla \partial_tn^\lambda\|_{s}^2\right)
     \\&\leq& C\exp{Ct},
\eeq for $t\in[0,T]$, where $T^\lambda>T$ and $T^\lambda\rightarrow\infty$ as
$\lambda\rightarrow\infty$. Furthermore, as
$\lambda\rightarrow\infty$, $(\rho^\lambda,u^\lambda,n^\lambda)$
converges to the unique global strong solution $(1,u,n)$ of the
incompressible system of liquid crystals (\ref{l40}), and
\beq
    \label{l55}
    \|u\|_s^2+\|\nabla n\|_{s}^2+\int_0^t\left(\mu\|\nabla u\|_s^2+\theta \|\nabla^2 n\|_{s}^2\right)\leq C_1\varepsilon_0,
\eeq
for any $t>0$, where $C_1$ is a uniform constant independent of $\varepsilon_0$ and $t$.
\end{Theorem}

\begin{Remark} Since this paper is mainly concerned with the limit as $\lambda\rightarrow\infty$, we do not state the global existence of the compressible model
for fixed $\lambda>0$ and small initial data. The global existence of the solutions to the slightly compressible model is considered in  \cite{Ding2}.
\end{Remark}

\begin{Theorem}\label{th.1.3} Under the assumptions of Theorem \ref{th.1.1}, the
convergence rates of $\rho^\lambda$, $u^\lambda$ and $n^\lambda$
$(\lambda\rightarrow \infty)$ are deduced as
\beq
    \label{l56}
    \lambda\|\rho^\lambda-1\|_s^2+\|u^\lambda-u\|^2+\|n^\lambda-n\|_2^2+\int_0^t\left(\|u^\lambda-u\|_1^2+\|n^\lambda-n\|_3^2\right)\leq C\lambda^{-1},
\eeq for $t\in [0,T_0]$. Furthermore, we have
\beq
    \label{l57}
    \|\nabla(\rho^\lambda-1)\|_{s-2}^2\leq C\lambda^{-4}, \ \forall \ t\in [0,T_0].
\eeq

The statement also holds for the strong solution given in Theorem \ref{th.1.2} for $t\in [0,T^\lambda)$.
\end{Theorem}
\begin{Remark} In Theorem \ref{th.1.3}, we do not give the
convergence rates for the higher order derivatives of $u^\lambda$
and $n^\lambda$ since we know little about the convergence rate of
the pressure.
\end{Remark}

\setcounter{section}{2} \setcounter{equation}{0}
\section{Local Existence and Uniform Stability}
In this section, we will give the uniform estimates for our results and then prove Theorem \ref{th.1.1}.
Let $U_0=\left(1+\overline{\rho}^\lambda_0,u_0+\overline{u}^\lambda_0,\frac{n_0+\overline{n}^\lambda_0}{|n_0+\overline{n}^\lambda_0|}\right)$.
We consider a set of functions $B^\lambda_{T_0}(U_0)$ contained in $\left\{(\rho,u,n):(\rho,u,\nabla n)\in L^\infty([0,T_0];H^s)\cap {\rm Lip}([0,T_0];H^{s-1})\right\}$
with $s\geq 3$ and defined by
\beq \label{s1}
\begin{cases}|\lambda(\rho-1)|+|u-u_0|+|\nabla n-\nabla n_0|<\delta,
             \\E_s(U(t))+\|n\|^2+\displaystyle\int_0^t\left(\mu\|\nabla
             u\|_s^2+(\kappa+\mu)\|\nabla\cdot u\|_s^2+\theta \|\nabla n\|_{s+1}^2\right)\leq
             K_1,
             \\E_{s-1}(\partial_tU(t))+\|n_t\|^2+\displaystyle\int_0^t\left(\mu\|\nabla
             \partial_t
             u\|_{s-1}^2+(\kappa+\mu)\|\nabla\cdot \partial_tu\|_{s-1}^2+\theta \|\nabla \partial_tn\|_{s}^2\right)\leq
             K_2.
             %\\[2mm]
\end{cases}\eeq

For any $V=(\xi^\lambda,\upsilon^\lambda,m^\lambda)\in B^\lambda_{T_0}(U_0)$, define $U=(\rho^\lambda,u^\lambda,n^\lambda)=\Lambda(V)$ as the unique
solution of the following ``linearized" problem
\beq \label{s2}
\begin{cases}\rho^\lambda_t+(v^\lambda\cdot\nabla)\rho^\lambda+\xi^\lambda\nabla\cdot u^\lambda=0,\\
         u^\lambda_t+(v^\lambda\cdot \nabla) u^\lambda+\lambda^2\frac{P'(\xi^\lambda)}{\xi^\lambda}\nabla\rho^\lambda=\frac{\mu}{\xi^\lambda}
         \Delta{u^\lambda}+\frac{(\kappa+\mu)}{\xi^\lambda}\nabla(\nabla\cdot u^\lambda)-\frac{\nu}{\xi^\lambda}\Delta n^\lambda\cdot\nabla n^\lambda,\\
        {n^\lambda_t}-\theta\Delta{n}^\lambda=-({v}^\lambda\cdot\nabla){m}^\lambda+\theta|\nabla m^\lambda|^2{m}^\lambda, %\\[2mm]
\end{cases}\eeq
for which the existence and uniqueness of the solutions is guaranteed by the standard theory of parabolic equations and Navier--Stokes equations.
Now we are to show that for appropriate choices of $T_0$, $\delta$,
$K_1$ and $K_2$ independent of $\lambda$, $\Lambda$ maps
$B^\lambda_{T_0}(U_0)$ into itself and it is a contraction in
certain function spaces. We emphasize that the solutions will depend
on the value of the parameter $\lambda$, but for convenience, the
dependence will not always be displayed in this section.

We need the following lemma for the proofs.
\begin{Lemma} \cite{Hagstrom2,Klainerman,Lei2} \label{le:2.1} Let $s>\frac{N}{2}$. For any functions $f$, $g$ $($possibly vector-valued in $\mathbb{R}^n$$)$ in Sobolev space
$H^s(T^N,\mathbb{R}^n)$, $\Phi\in C^s(\mathbb{R}^n)$, $\|\nabla^j\Phi\|_\infty<+\infty$ $(j=1,2,\ldots,s)$, and multi-index $\alpha$ satisfying $|\alpha|\leq s$, we have the Sobolev inequality:
\beq
    \|f\|_\infty\leq \|f\|_s,
\eeq
the estimate based on the chain rule:
\beq
    \label{4.7}
    \|\nabla^r(\Phi\circ f)\|\leq C(1+\|f\|_\infty^{s-1})\|\nabla f\|_{s-1}, \ {\rm for} \  1\leq r\leq s,
\eeq
and the estimates based on Leibniz' rule:
\beq
    \label{4.8}
    \|\nabla^\alpha(fg)\|&\leq& C(\|f\|_\infty\|\nabla^\alpha g\|+\|g\|_\infty\|\nabla^\alpha f\|),
    \\
    \|\nabla^\alpha(fg)-f\nabla^\alpha g\|&\leq& C(\|\nabla f\|_\infty\|g\|_{s-1}+\|g\|_\infty\|\nabla f\|_{s-1}),
\eeq
where the constants $C$ are independent of $f$, $g$, but may depend on $|\alpha|$ and $\Phi$.
\end{Lemma}

Before proceeding any further, we apply $D^{\alpha_1}$ to the first and the second equations of (\ref{s2}) and $D^{\alpha_2}$ to the third one respectively, and then we get
\beq \label{s3}
\begin{cases}
         \partial_t D^{\alpha_1} \rho+(v\cdot\nabla) D^{\alpha_1}\rho+\xi\nabla\cdot D^{\alpha_1} u=\Pi_1,\\
         \partial_tD^{\alpha_1} u+(v\cdot \nabla) D^{\alpha_1} u+\lambda^2\frac{P'(\xi)}{\xi}\nabla D^{\alpha_1}
         \rho=\Pi_2,\\
        \partial_t D^{\alpha_2} n-\theta\Delta D^{\alpha_2} n=-D^{\alpha_2}\left(v\cdot\nabla m\right)
        +\theta D^{\alpha_2}(|\nabla m|^2m),%\\[2mm]
\end{cases}\eeq where
\beqno
      \Pi_1&=&-[D^{\alpha_1}(v\cdot\nabla\rho)-(v\cdot\nabla) D^{\alpha_1}\rho]-[D^{\alpha_1}(\xi\nabla\cdot u)-\xi\nabla\cdot D^{\alpha_1} u],
      \quad\quad\quad\quad\quad\quad\quad\quad\quad\quad\quad
      \\
      \Pi_2&=&\frac{\mu}{\xi}\Delta D^{\alpha_1} u
      +\frac{\kappa+\mu}{\xi}\nabla D^{\alpha_1} (\nabla\cdot u)-\frac{\nu}{\xi}D^{\alpha_1} (\Delta n\cdot\nabla n)
      \\&&
      -[D^{\alpha_1}(v\cdot\nabla u)-(v\cdot\nabla) D^{\alpha_1} u]-\lambda^2\left[D^{\alpha_1}\left(\frac{P'(\xi)}{\xi}\nabla\rho\right)-\frac{P'(\xi)}{\xi}\nabla D^{\alpha_1}\rho\right]\\&&+\left[D^{\alpha_1}\left(\frac{\mu}{\xi}\Delta u\right)-\frac{\mu}{\xi}\Delta D^{\alpha_1} u\right]+\left\{D^{\alpha_1}\left[\frac{\kappa+\mu}{\xi}\nabla(\nabla\cdot u)\right]-\frac{\kappa+\mu}{\xi}\nabla D^{\alpha_1}(\nabla\cdot u)\right\}
      \\&&
      -\left\{D^{\alpha_1}\left[\frac{\nu}{\xi}(\Delta n\cdot\nabla n)\right]-\frac{\nu}{\xi}D^{\alpha_1}(\Delta n\cdot\nabla n)\right\}.
\eeqno

We will prove that $\Lambda$ maps $B^\lambda_{T_0}(U_0)$ into itself by
two steps and denote by $C$ the constants
independent of $\lambda$, $K_1$ and $K_2$ in these two steps. Without loss of generality, we assume that $T_0^{-1}$, $\lambda$, $K_1$, $K_2>1$.

\textbf{Step One: Estimates of $n$.}

Taking the $L^2$ inner product of $(\ref{s3})_3$ with $D^{\alpha_2} n$,
and then integrating by parts, we have
 \beq
    \label{s4}
    \nonumber
    &&\frac12\frac{{\rm d}}{{\rm d}t}\int_{T^N} |D^{\alpha_2} n|^2+\theta\int_{T^N}|D^{\alpha_2}\nabla n|^2
    \\&=& \nonumber
    -\int_{T^N}D^{\alpha_2} (v\cdot\nabla m)\cdot D^{\alpha_2} n+\theta \int_{T^N} D^{\alpha_2}(|\nabla m|^2m)\cdot D^{\alpha_2} n
    \\&=&
    N_1+N_2.
\eeq
Next, we will give the estimates of $N_1$ and $N_2$ in the following two cases respectively.

Case 1: $D^{\alpha_2}=\nabla^{\alpha_2}$ with $|\alpha_2|\leq s+1$.

It follows from Lemma \ref{le:2.1}, the Sobolev embedding $H^2(T^N)\hookrightarrow L^\infty(T^N)$ for $N=2,3$, and the Cauchy inequality that
\beq
    \label{s9}
    \nonumber
    |N_1|&\leq& \nonumber
    \|\nabla^{\alpha_2}(v\cdot\nabla m)\|\|\nabla^{\alpha_2} n\|
    \\&\leq& \nonumber
    C(\|v\|_\infty\|\nabla m\|_{s+1}+\|\nabla m\|_\infty\|v\|_{s+1})\|\nabla^{\alpha_2} n\|
    \\&\leq& \nonumber
    (K_1^{\frac12}\|\nabla m\|_{s+1}+K_1+K_1^{\frac12}\|\nabla v\|_{s})\|\nabla^{\alpha_2} n\|
    \\&\leq&
    K_1^2\|\nabla^{\alpha_2} n\|^2+CK_1^{-1}\left(\|\nabla m\|_{s+1}^2+\|\nabla v\|_s^2\right)+C,
    \\ \nonumber
    \label{s10}
    |N_2|&\leq& \nonumber
    \|\nabla^{\alpha_2}(|\nabla m|^2m)\|\|\nabla^{\alpha_2} n\|
    \\&\leq& \nonumber
    C\left(\|\nabla m\|_\infty^2\|m\|_{s+1}+\|m\|_\infty\|\nabla m\|_\infty\|\nabla m\|_{s+1}\right)\|\nabla^{\alpha_2} n\|
    \\&\leq& \nonumber
    CK_1\|\nabla m\|_{s+1}\|\nabla^{\alpha_2} n\|
    \\&\leq&
    K_1^3\|\nabla^{\alpha_2}n\|^2+CK_1^{-1}\|\nabla m\|^2_{s+1}.
\eeq

Case 2: $D^{\alpha_2}=\partial_t\nabla^{\alpha_1}$ with $|\alpha_1|\leq s$.

Similarly, we have
\beq
    \label{s11}
    \nonumber
    |N_1|&\leq& \nonumber
    \left(\|\nabla^{\alpha_1}(v_t\cdot\nabla m)\|+\|\nabla^{\alpha_1}(v\cdot\nabla m_t)\|\right)\|\nabla^{\alpha_1}n_t\|
    \\&\leq& \nonumber
    C\left(\|v_t\|_\infty\|\nabla m\|_s+\|\nabla m\|_\infty\|v_t\|_{s}+\|v\|_\infty\|\nabla m_t\|_s+\|\nabla m_t\|_\infty\|v\|_s\right)\|\nabla^{\alpha_1}n_t\|
    \\&\leq& \nonumber
    C\left(K_1^{\frac12}K_2^{\frac12}+K_1^{\frac12}\|\nabla v_t\|_{s-1}+K_1^{\frac12}\|\nabla m_t\|_s\right)\|\nabla^{\alpha_1}n_t\|
    \\&\leq&
    K_1K_2\|\nabla^{\alpha_1}n_t\|^2+CK_2^{-1}\left(\|\nabla v_t\|_{s-1}^2+\|\nabla m_t\|_s^2\right).
    \\ \nonumber
    \label{s12}
    |N_2|&\leq& \left(\|\nabla^{\alpha_1}\left((\nabla m:\nabla m_t)m\right)\|+\|\nabla^{\alpha_1}\left(|\nabla m|^2m_t\right)\|\right)\|\nabla^{\alpha_1}n_t\|
    \\&\leq& \nonumber
    C\left(\|\nabla m\|_\infty\|\nabla m_t\|_\infty\|m\|_s+\|m\|_\infty\left(\|\nabla m\|_\infty\|\nabla m_t\|_s+\|\nabla m_t\|_\infty\|\nabla m\|_s\right)\right.
    \\&& \nonumber
    \left.+\|\nabla m\|_\infty^2\|m_t\|_s+\|m_t\|_\infty\|\nabla m\|_\infty\|\nabla m\|_s\right)\|\nabla^{\alpha_1} n_t\|
    \\&\leq& \nonumber
    C\left(K_1K_2^{\frac12}+K_1\|\nabla m_t\|_s\right)\|\nabla^{\alpha_1} n_t\|
    \\&\leq&
    K_1^2K_2\|\nabla^{\alpha_1}n_t\|^2+CK_2^{-1}\|\nabla m_t\|_s^2.
\eeq

On one hand, substituting the estimates (\ref{s9}) and (\ref{s10}) into (\ref{s4}), and then summing over $\alpha_2$, we obtain
\beq
    \label{s16}
    &&\nonumber
    \frac{{\rm d}}{{\rm d}t}\sum\limits_{|\alpha_2|\leq s+1}\|\nabla^{\alpha_2} n\|^2+\theta\sum\limits_{|\alpha_2|\leq s}\|\nabla^{\alpha_2}\nabla n\|^2
    \\&\leq&
     CK_1^3\|n\|_{s+1}^2+CK_1^{-1}\left(\|\nabla m\|_{s+1}^2+\|\nabla v\|_s^2\right)+C.
\eeq

Then integrating the above inequality over $[0,t]\subseteq [0,T_0]$, we find that
\beq
    \label{jf3.14}
    \|n\|_{s+1}^2(t)+\theta\int_0^t\|\nabla n\|_{s+1}^2\leq CK_1^3\int_0^t\|n\|_{s+1}^2+\|n(x,0)\|_{s+1}^2+C
\eeq
which directly yields from the Gronwall inequality and (\ref{l46})--(\ref{l48}) that
\beq
    \label{s17}
    \nonumber
    \|n\|_{s+1}^2(t)&\leq& \exp{(CK_1^3T_0)}\left(\|n(x,0)\|_{s+1}^2+C\right)
    \\&\leq&
    \exp{(CK_1^3T_0)}\left(\|n_0\|_{s+1}^2+\lambda^{-2}\delta_0+C\right)\leq C
\eeq for $t\in[0,T_0]$, provided that $T_0<T_1=K_1^{-3}$. It follows from (\ref{jf3.14}) that $\theta\int_0^t\|\nabla n\|_{s+1}^2\leq C$ for $t\in[0,T_0]$,
provided that $T_0<T_1$.

On the other hand, substituting the estimates (\ref{s11}) and (\ref{s12}) into (\ref{s4}), and then summing over $\alpha_1$, we obtain
\beq
    \label{xin22}
    &&\nonumber
    \frac{{\rm d}}{{\rm d}t}\sum\limits_{|\alpha_1|\leq s}\int_{T^N} |\nabla^{\alpha_1} n_t|^2+\theta\sum\limits_{|\alpha_1|\leq s}\int_{T^N}|\nabla^{\alpha_1}\nabla n_t|^2
    \\&\leq&
    CK_1^2K_2\|n_t\|_{s}^2+CK_2^{-1}\left(\|\nabla v_t\|_{s-1}^2+\|\nabla m_t\|_s^2\right).
\eeq

Recalling the constraints of the initial data and
(\ref{l38})$_3$, and then using Lemma \ref{le:2.1}, we get
\beq
    \label{xin23}
    \nonumber
    \|n_t(x,0)\|_{s}&\leq&
    C\left(\|(u(x,0)\cdot\nabla) n(x,0)\|_{s}+\|\Delta n(x,0)+|\nabla n(x,0)|^2n(x,0)\|_{s}\right)
    \\&\leq&\nonumber
    C\left(\|u_0\|_{s}+\lambda^{-1}\delta_0\right)\left(\|\nabla n_0\|_{s}+\lambda^{-1}\delta_0\right)+C\left(\|\Delta n_0\|_{s}+\lambda^{-1}\delta_0\right)
    \\&&+\nonumber
    C\left(\|\nabla n_0\|_{s}^2+\lambda^{-2}\delta_0^2\right)\left(\|n_0\|_{s}+\lambda^{-1}\delta_0\right)+C\left(\|\nabla n_0\|_{s}^2+\lambda^{-2}\delta_0^2\right)
    \\&\leq&
    C.
\eeq

Then integrating (\ref{xin22}) over $[0,t]\subseteq [0,T_0]$ and using the Gronwall inequality, we have
\beq
    \label{xin24}
    \|n_t\|_{s}^2+\theta\int_0^t\|\nabla n_t\|_{s}^2ds\leq C
\eeq for $t\in[0,T_0]$, provided that $T_2\doteq\{T_1,K_1^{-2}K_2^{-1}\}$.

\textbf{Step Two: Estimates of $\rho$ and $u$.}

Taking the $L^2$ inner product of $(\ref{s3})_1$ and $(\ref{s3})_2$ with $\lambda^2\frac{P'(\xi)}{\xi}D^{\alpha_1} (\rho-1)$
and $\xi D^{\alpha_1} u$ respectively, and then integrating by parts, we arrive at
\beq
    \label{s23}
    \nonumber&&\frac12\frac{{\rm d}}{{\rm d}t}\int_{T^N}\left(\frac{P'(\xi)}{\xi}|\lambda D^{\alpha_1}(\rho-1)|^2+\xi|D^{\alpha_1} u|^2\right)
    \\&&
    +\mu\int_{T^N}|D^{\alpha_1}\nabla u|^2+(\kappa+\mu)\int_{T^N}|D^{\alpha_1}(\nabla\cdot u)|^2=
    \sum_{k=1}^9I_k,
\eeq
where
\beqno
    I_1&=&\frac{1}{2}\int_{T^N}\left[|\lambda D^{\alpha_1}(\rho-1)|^2\partial_t\left(\frac{P'(\xi)}{\xi}\right)+\xi_t|D^{\alpha_1} u|^2\right],
    \\
    I_2&=&\frac{1}{2}\int_{T^N}\left[|\lambda D^{\alpha_1}(\rho-1)|^2\nabla\cdot\left(\frac{P'(\xi)}{\xi}v\right)+\nabla\cdot(\xi v)|D^{\alpha_1} u|^2\right],
    \\
    I_3&=&\lambda^2\int_{T^N} P''(\xi)D^{\alpha_1}(\rho-1)D^{\alpha_1} u\cdot\nabla\xi,
    \\
    I_4&=&-\lambda^2\int_{T^N}\frac{P'(\xi)}{\xi}D^{\alpha_1}(\rho-1)\biggr\{[D^{\alpha_1}(v\cdot\nabla\rho)-v\cdot\nabla D^{\alpha_1}\rho]+[D^{\alpha_1}(\xi\nabla\cdot u)-\xi\nabla \cdot D^{\alpha_1} u]\biggr\},
    \\
    I_5&=&-\nu\int_{T^N}D^{\alpha_1}(\Delta n\cdot\nabla n)\cdot D^{\alpha_1} u,
    \\
    I_6&=&-\int_{T^N}\xi\left[D^{\alpha_1}(v\cdot\nabla u)-v\cdot\nabla D^{\alpha_1} u\right]\cdot D^{\alpha_1} u,
    \\
    I_7&=&-\lambda^2\int_{T^N}\xi\left[D^{\alpha_1}\left(\frac{P'(\xi)}{\xi}\nabla\rho\right)-\frac{P'(\xi)}{\xi}\nabla D^{\alpha_1} \rho\right]\cdot D^{\alpha_1} u,
    \\
    I_8&=&\int_{T^N}\xi\left[D^{\alpha_1}\left(\frac{\mu}{\xi}\Delta u\right)-\frac{\mu}{\xi}\Delta D^{\alpha_1} u\right]\cdot D^{\alpha_1} u
    \\&&
    +\int_{T^N}\xi\left\{D^{\alpha_1}\left[\frac{\kappa+\mu}{\xi}\nabla(\nabla\cdot u)\right]-\frac{\kappa+\mu}{\xi}\nabla D^{\alpha_1}(\nabla\cdot u)\right\}\cdot D^{\alpha_1} u,
    \\
    I_9&=&-\int_{T^N}\xi\left\{D^{\alpha_1}\left[\frac{\nu}{\xi}(\Delta n\cdot\nabla n)\right]-\frac{\nu}{\xi}D^{\alpha_1}(\Delta n\cdot\nabla n)\right\}\cdot D^{\alpha_1} u.\quad\quad\quad\quad\quad\quad\quad\quad\quad\quad
\eeqno

First of all, we give the estimates of $I_1$--$I_3$. Choose $\delta$ small enough so that (\ref{l45}) and $|\xi-1|\leq \frac12$ hold.
Then by the smoothness of the pressure $P(\cdot)$, we have
\beq
    \label{xin26}
    \nonumber
    |I_1|&\leq& \left\|\frac{\xi P''(\xi)-P'(\xi)}{\xi^2}\right\|_\infty\|\xi_t\|_\infty\|\lambda D^{\alpha_1}(\rho-1)\|^2+\|\xi_t\|_\infty\|D^{\alpha_1} u\|^2
    \\&\leq&\nonumber
    C\|\xi_t\|_2\left(\|\lambda D^{\alpha_1}(\rho-1)\|^2+\|D^{\alpha_1} u\|^2\right)
    \\&\leq&
    C\lambda^{-1} K_2^{\frac12}\left(\|\lambda D^{\alpha_1}(\rho-1)\|^2+\|D^{\alpha_1} u\|^2\right),
    \\ \label{xin27}
    |I_2|&\leq&\nonumber \left(\left\|\frac{\xi P''(\xi)-P'(\xi)}{\xi^2}\right\|_\infty\|\nabla\xi\|_\infty\|v\|_\infty+\left\|\frac{P'(\xi)}{\xi}\right\|_\infty\|\nabla\cdot v\|_\infty\right)\|\lambda D^{\alpha_1}(\rho-1)\|^2
    \\& &\nonumber
    +\left(\|\nabla \xi\|_\infty\|v\|_\infty+\|\xi\|_\infty\|\nabla\cdot v\|_\infty\right)\|D^{\alpha_1} u\|^2
    \\&\leq&\nonumber
    C\left(\lambda^{-1}\|\lambda\nabla\xi\|_2\|v\|_2+\|v\|_3)(\|\lambda D^{\alpha_1}(\rho-1)\|^2+\|D^{\alpha_1} u\|^2\right)
    \\&\leq&
    C\left(\lambda^{-1} K_1+K_1^{\frac12}\right)\left(\|\lambda D^{\alpha_1}(\rho-1)\|^2+\|D^{\alpha_1} u\|^2\right),
    \\\label{xin28}
    |I_3|&\leq& \nonumber
    \lambda\|P''(\xi)\|_\infty\|\nabla \xi\|_\infty\left(\|\lambda D^{\alpha_1}(\rho-1)\|^2+\|D^{\alpha_1} u\|^2\right)
    \\&\leq&\nonumber
    C\|\lambda\nabla\xi\|_2\left(\|\lambda D^{\alpha_1}(\rho-1)\|^2+\|D^{\alpha_1} u\|^2\right)
    \\&\leq&
    CK_1^{\frac12}\left(\|\lambda D^{\alpha_1}(\rho-1)\|^2+\|D^{\alpha_1} u\|^2\right).
\eeq

Secondly, we give the estimates of $I_4$--$I_9$ in the following three cases.

Case 1: $|\alpha_1|=0$ or $D^{\alpha_1}=\partial_t$.

If $|\alpha_1|=0$, then the quantities $I_j$, $
4\leq j\leq9$, are all equal to zero except for $I_5$. Thus it suffices to estimate $I_5$. Using (\ref{s17}), we have
\beq
    \label{s37}
    |I_5|\leq C\|u\|\|\nabla n\|_\infty\|\Delta n\|\leq C\|u\|^2+C.
\eeq

If $D^{\alpha_1}=\partial_t$, then by using (\ref{s17}), (\ref{xin24}) and the Cauchy inequality, we have
\beq
    \nonumber
    \label{xin30}
    |I_4|&\leq& C\|\lambda\rho_t\|\left(\|v_t\|_\infty\|\lambda\nabla\rho\|+\|\lambda\xi_t\|_\infty\|\nabla\cdot
    u\|\right)
    \\&\leq&
    CK_2^{\frac12}\left(\|\lambda\nabla\rho\|^2+\|\lambda\rho_t\|^2+\|\nabla\cdot u\|^2\right), \\ \label{s38}
    |I_5|&\leq& C\left(\|\Delta n_t\|\|\nabla n\|_\infty+\|\Delta n\|_\infty\|\nabla n_t\|\right)\|u_t\|
    \\&\leq&
    C\left(\|\Delta n_t\|\|\nabla n\|_2+\|\Delta n\|_2\|\nabla n_t\|\right)\|u_t\|
    \leq \label{s39}
    C\|u_t\|^2+C,
    \\ \label{s40}
    |I_6|&\leq& C\|v_t\|_\infty\|\nabla u\|\|u_t\|\leq CK_2^{\frac12}(\|\nabla u\|^2+\|u_t\|^2),
    \\|I_7|&\leq& C\|\lambda\xi_t\|_\infty\|\lambda\nabla\rho\|\|u_t\|\leq CK_2^{\frac12}\left(\|\lambda\nabla\rho\|^2+\|u_t\|^2\right), \quad\  \label{s41}
    \\
    |I_8|&\leq&
    C\lambda^{-1}\|\lambda\xi_t\|_\infty\left(\|\Delta u\|+\|\nabla(\nabla\cdot u)\|\right)\|u_t\|
    \leq C\lambda^{-1}K_2^{\frac12}(\|u\|_2^2+\|u_t\|^2),\quad\label{s42}
    \\ \label{xin36}
    |I_9|&\leq&C\lambda^{-1}\|\lambda\xi_t\|_\infty\|\Delta n\|\|\nabla n\|_\infty\|u_t\|
    \leq
    C\lambda^{-2}K_1\|u_t\|^2+C.
\eeq

Case 2: $D^{\alpha_1}=\nabla^{\alpha_1}$ for $1\leq|\alpha_1|\leq s$, where $s\geq3$.

By Lemma \ref{le:2.1}, we have
\beq
    \nonumber |I_4|&\leq& C\left(\|\nabla v\|_\infty\|\lambda\nabla\rho\|_{s-1}+\|\lambda\nabla\rho\|_\infty\|v\|_s+\|\lambda\nabla\xi\|_\infty\|\nabla\cdot u\|_{s-1}+\|\nabla\cdot u\|_\infty\|\lambda\nabla \xi\|_{s-1}\right)
    \\& &\nonumber
    \times\|\lambda\nabla^{\alpha_1}(\rho-1)\|
    \\&\leq&
    CK_1^{\frac12}\left(\|\lambda(\rho-1)\|_s^2+\|u\|_s^2\right),\label{s44}
    \\ \nonumber
    |I_5|&=&\left|\nu\int_{T^N}\nabla^{\alpha_1}\left(\nabla n\odot\nabla n-\frac{|\nabla n|^2}{2}I_N\right)\cdot\nabla\nabla^{\alpha_1} u\right|
    \\&\leq&
    C\|\nabla u\|_s\|\nabla n\|_\infty\|\nabla n\|_{s}
    \leq
    \tau\|\nabla u\|_s^2+C(\tau),\label{s45}
\eeq
where we have used integration by parts and (\ref{s17}) in the estimate of $I_5$.

Similarly as in (\ref{s44}), one obtains
\beq
    \label{s46}
    |I_6|&\leq& C\|u\|_s(\|\nabla v\|_\infty\|\nabla u\|_{s-1}+\|\nabla u\|_\infty\|v\|_s)\leq CK_1^{\frac12}\|u\|_s^2.
\eeq

On the other hand, by Lemma \ref{le:2.1}, we have
\beq
    \nonumber
    |I_7|&\leq&\nonumber C\lambda^2\left(\left\|\frac{\xi P''(\xi)-P'(\xi)}{\xi^2}\nabla\xi\right\|_\infty\|\nabla\rho\|_{s-1} +\|\nabla\rho\|_\infty\left\|\nabla\left(\frac{P'(\xi)}{\xi}\right)\right\|_{s-1}\right)\|u\|_s
    \\&\leq&\nonumber
    C\left(\|\lambda\nabla\xi\|_2\|\lambda\nabla\rho\|_{s-1}+\|\lambda\nabla\rho\|_2\|\lambda\nabla\xi\|_{s-1}\right)\|u\|_s
    \\&\leq&
    CK_1^{\frac12}\left(\|\lambda(\rho-1)\|_s^2+\|u\|_s^2\right),\label{s48}
    \\ \nonumber
    |I_8|&\leq&\nonumber C\|u\|_s\left(\left\|\nabla\left(\frac{1}{\xi}\right)\right\|_\infty\|\Delta u\|_{s-1}+\|\Delta u\|_\infty\left\|\nabla\left(\frac{1}{\xi}\right)\right\|_{s-1}\right.
    \\&&\nonumber
    \left.+\left\|\nabla\left(\frac{1}{\xi}\right)\right\|_\infty\|\nabla(\nabla\cdot u)\|_{s-1}+\|\nabla(\nabla\cdot u)\|_\infty\left\|\nabla\left(\frac{1}{\xi}\right)\right\|_{s-1}\right)
    \\&\leq&\nonumber
    C\lambda^{-1} K_1^{\frac12}\|u\|_s(\|\nabla u\|_s+\|\nabla\cdot u\|_s)
    \\&\leq&
    \tau\|\nabla u\|_s^2+\tau\|\nabla\cdot u\|_s^2+C(\tau)\lambda^{-2} K_1\|u\|_s^2,\label{s49}
    \\ \nonumber
    |I_9|&\leq&\nonumber C\|u\|_s\left(\left\|\nabla\left(\frac{1}{\xi}\right)\right\|_\infty\|\Delta n\cdot\nabla n\|_{s-1}+\|\Delta n\cdot\nabla n\|_\infty\left\|\nabla\left(\frac{1}{\xi}\right)\right\|_{s-1}\right)
    \\&\leq&\nonumber
    C\lambda^{-1}\|u\|_s\left(\|\lambda\nabla\xi\|_\infty(\|\Delta n\|_\infty\|\nabla n\|_{s-1}+\|\nabla n\|_\infty\|\Delta n\|_{s-1})
    +\|\Delta n\|_2\|\nabla n\|_2\|\lambda\nabla\xi\|_{s-1}\right)
    \\&\leq&
    C\lambda^{-2}K_1\|u\|_s^2+C, \label{s50}
\eeq
where we have used (\ref{s17}) in the estimate of $I_9$.

Putting the estimates (\ref{xin26})--(\ref{s37}) and (\ref{s44})--(\ref{s50}) together, summing over $\alpha_1$, and then choosing $\tau$ small enough, we get
\beq
    \label{s51}
    \nonumber &&\frac{{\rm d}}{{\rm d}t}\sum\limits_{|\alpha_1|\leq s}\int_{T^N} \left(\frac{P'(\xi)}{\xi}\lambda^2|\nabla^{\alpha_1}(\rho-1)|^2+\xi|\nabla^{\alpha_1} u|^2\right)
    \\&& \nonumber
    +\sum\limits_{|\alpha_1|\leq s}\mu\int_{T^N}|\nabla^{\alpha_1}\nabla u|^2+\sum\limits_{|\alpha_1|\leq s}(\kappa+\mu)\int_{T^N}|\nabla^{\alpha_1}(\nabla\cdot u)|^2
    \\&\leq&
    C\left(K_1+K_2\right)\left(\|\lambda(\rho-1)\|_s^2+\|u\|_s^2\right)+C.
\eeq

Recalling the constraints of the initial data (\ref{l47}) and (\ref{l48}), we have
\beq
    \label{s52}
    \|\lambda\overline{\rho}_0\|_s^2+\|u_0+\overline{u}_0\|_s^2\leq \left(3\lambda^{-2}\delta_0^2+2\|u_0\|_s^2\right).
\eeq

Then by (\ref{l45}), (\ref{s51}), (\ref{s52}) and the Gronwall inequality, we get
\beq
    \label{s53}
    \|\lambda(\rho-1)\|_{s}^2(t)+\|u\|_s^2(t)
    \leq
    \exp{(C\left(K_1+K_2\right)T_0)}\left(3\lambda^{-2}\delta_0^2+2\|u_0\|_s^2+CT_0\right)
    \leq C
\eeq
for $t\in[0,T_0]$, provided that $T_0<T_2$. Furthermore, we integrate (\ref{s51}) over
$[0,t]\subseteq [0,T_0]$ and then get
\beq
    \label{s54}
    \mu\int_0^t\|\nabla u\|_{s}^2+(\kappa+\mu)\int_0^t\|\nabla\cdot u\|_s^2
    \leq C  \ {\rm for} \ t\in[0,T_0].
\eeq

Case 3: $D^{\alpha_1}=\nabla^\beta\partial_t$ for $1\leq|\beta|\leq s-1$, where $s\geq3$.

On one hand, for
$I_5$, $I_6$ and $I_9$, by (\ref{s17}), (\ref{xin24}), (\ref{s53}) and the Cauchy inequality, we are led to
\beq
    \nonumber
    |I_5|&\leq& \nonumber
    C\left(\|\Delta n_t\|_\infty\|\nabla n\|_{s-1}+\|\nabla n\|_\infty\|\Delta n_t\|_{s-1}\right)\|\nabla^\beta u_t\|
    \\&& \nonumber
    +C\left(\|\Delta n\|_\infty\|\nabla n_t\|_{s-1}+\|\nabla n_t\|_\infty\|\Delta n\|_{s-1}\right)\|\nabla^\beta u_t\|
    \\&\leq&
    C\|u_t\|_{s-1}^2+C(\|\Delta n_t\|_{s-1}^2+1),\label{s57}
    \\ \nonumber
    |I_6|&\leq& \nonumber
    C\|u_t\|_{s-1}\left(\|\nabla v\|_\infty\|\nabla u_t\|_{s-2}+\|\nabla u_t\|_\infty\|v\|_{s-1}+\|v_t\|_\infty\|\nabla u\|_{s-1}+\|\nabla u\|_\infty\|v_t\|_{s-1}\right)
%    \\&\leq& \nonumber
%    C\|u_t\|_{s-1}\left(K_1^{\frac12}\|u_t\|_{s-1}+K_1^{\frac12}\|\nabla u_t\|_2+K_2^{\frac12}\|u\|_s\right)
    \\&\leq&
    \tau\|\nabla u_t\|_2^2+C(\tau)\left(K_1+K_2\right)\|u_t\|_{s-1}^2+C,\label{xin47}
    \\ \nonumber
    |I_9|&\leq&\nonumber
    C\|\nabla^\beta u_t\|\bigg(\left\|\left(\frac{1}{\xi}\right)_t\right\|_\infty\|\Delta n\cdot\nabla n\|_{s-1}+\|\Delta n\cdot\nabla n\|_\infty\left\|\left(\frac{1}{\xi}\right)_t\right\|_{s-1}
    \\&&+\nonumber
    \left\|\nabla\left(\frac{1}{\xi}\right)\right\|_\infty\|(\Delta n\cdot\nabla n)_t\|_{s-2}+\|(\Delta n\cdot\nabla n)_t\|_\infty\left\|\nabla\left(\frac{1}{\xi}\right)\right\|_{s-2}\bigg)
    \\&\leq&\nonumber
    C\|u_t\|_{s-1}\big(\left\|\xi_t\right\|_\infty\|\Delta n\|_\infty\|\nabla n\|_{s-1}+\left\|\xi_t\right\|_\infty\|\nabla n\|_\infty\|\Delta n\|_{s-1}
    \\&&\nonumber
    +\|\Delta n\|_\infty\|\nabla n\|_\infty\|\xi_t\|_{s-1}(1+\|\nabla\xi\|_{s-2})+\left\|\nabla\xi\right\|_\infty\|\Delta n_t\|_\infty\|\nabla n\|_{s-2}
    \\&& \nonumber
    +\left\|\nabla\xi\right\|_\infty\|\nabla n\|_\infty\|\Delta n_t\|_{s-2}+\left\|\nabla\xi\right\|_\infty\|\Delta n\|_\infty\|\nabla n_t\|_{s-2}
    +\left\|\nabla\xi\right\|_\infty\|\nabla n_t\|_\infty\|\Delta n\|_{s-2}\\&& \nonumber+\|\Delta n\|_\infty\|\nabla n_t\|_\infty\left\|\nabla\xi\right\|_{s-2}
    +\|\Delta n_t\|_\infty\|\nabla n\|_\infty\left\|\nabla\xi\right\|_{s-2}\big)
    \\&\leq&
    C\left(K_1^2+K_2^2\right)\|u_t\|_{s-1}^2+C(\|\Delta n_t\|_2^2+1).\label{xin53}
\eeq
where we have used in (\ref{xin53}) the following fact
\beq
    \label{s59}
    \nonumber
    \left\|\left(\frac{1}{\xi}\right)_t\right\|_{s-1}&\leq& C\left\|\frac{1}{\xi^2}\xi_t\right\|_{s-1}
    \leq
    C\left\|\frac{1}{\xi^2}\right\|_\infty\|\xi_t\|_{s-1} +C\|\xi_t\|_\infty\left\|\frac{1}{\xi^2}\right\|_{s-1}
    \\&\leq&
    C\|\xi_t\|_{s-1}(1+\|\nabla\xi\|_{s-2}).
\eeq

On the other hand, for the terms $I_4$, $I_7$ and $I_8$, some extra discussions have to be given since the methods to estimate these three terms are different between
the case $s=3$ and $s\geq 4$. Indeed, if $s\geq4$, then we use Lemma \ref{le:2.1}, the Sobolev imbedding, the Cauchy inequality, (\ref{s17}), (\ref{xin24}), (\ref{s53})
and a similar calculation as in (\ref{s59}) to get
\beq
    |I_4|&\leq&\nonumber C\|\lambda\rho_t\|_{s-1}\left(\|\nabla v\|_\infty\|\lambda\nabla\rho_t\|_{s-2}+\|\lambda\nabla\rho_t\|_\infty\|v\|_{s-1} +\|v_t\|_\infty\|\lambda\nabla\rho\|_{s-1}\right.
    \\&&\nonumber
    \left.+\|\lambda\nabla\rho\|_\infty\|v_t\|_{s-1}+\|\lambda\nabla\xi\|_\infty\|\nabla\cdot u_t\|_{s-2}+\|\nabla\cdot u_t\|_\infty\|\lambda\nabla\xi\|_{s-2}\right.
    \\&&\nonumber
    \left.+\|\lambda\xi_t\|_\infty\|\nabla\cdot u\|_{s-1}+\|\nabla\cdot u\|_\infty\|\lambda\xi_t\|_{s-1}\right)
%    \\&\leq&\nonumber
%    C\|\lambda\rho_t\|_{s-1}\left(K_1^{\frac12}\|\lambda\rho_t\|_{s-1}+K_2^{\frac12}\|\lambda(\rho-1)\|_{s} +K_1^{\frac12}\|u_t\|_{s-1}+K_2^{\frac12}\|u\|_s\right)
    \\&\leq&
    C(K_1+K_2)\left(\|\lambda\rho_t\|_{s-1}^2+\|u_t\|_{s-1}^2\right)+C,\label{s56}
    \\
    \nonumber
    |I_7|&\leq&
    C\lambda^2\|u_t\|_{s-1}\left(\left\|\nabla\xi\right\|_\infty\|\nabla\rho_t\|_{s-2}
    +\|\nabla\rho_t\|_\infty\left\|\nabla\xi\right\|_{s-2}+\left\|\xi_t\right\|_\infty\|\nabla\rho\|_{s-1}\right.
    \\&&\nonumber\left.
    +\|\nabla\rho\|_\infty\left\|\xi_t\right\|_{s-1}(1+\|\nabla\xi\|_{s-2})\right)
%    \\&\leq&\nonumber
%    C\|u_t\|_{s-1}\left(K_1^{\frac12}\|\lambda\rho_t\|_{s-1} +K_2^{\frac12}\left(1+\lambda^{-1}K_1^{\frac12}\right)\|\lambda\nabla\rho\|_{s-1}\right)
    \\&\leq&
    C\left(K_1^2+K_2^2\right)\left(\|\lambda\rho_t\|_{s-1}^2+\|u_t\|_{s-1}^2\right)+C,\label{s60}
    \\ \label{s61}
    |I_8|&\leq&\nonumber
    C\|u_t\|_{s-1}\left(\left\|\nabla\xi\right\|_\infty\|\Delta u_t\|_{s-2}+\|\Delta u_t\|_\infty\left\|\nabla\xi\right\|_{s-2} +\left\|\xi_t\right\|_\infty\|\Delta u\|_{s-1}\right.
    \\&&\nonumber
    \left.+\|\Delta u\|_\infty\|\xi_t\|_{s-1}(1+\|\nabla\xi\|_{s-2}) +\left\|\nabla\xi\right\|_\infty\|\nabla(\nabla\cdot u_t)\|_{s-2}+\|\nabla(\nabla\cdot u_t)\|_\infty\left\|\nabla\xi\right\|_{s-2}\right.
    \\&& \nonumber
    \left.+\left\|\xi_t\right\|_\infty\|\nabla(\nabla\cdot u)\|_{s-1}+\|\nabla(\nabla\cdot u)\|_\infty\|\xi_t\|_{s-1}(1+\|\nabla\xi\|_{s-2})\right)
%    \\&\leq&\nonumber
%    C\lambda^{-1}\|u_t\|_{s-1}\bigg[K_1^{\frac12}(\|\nabla u_t\|_{s-1}+\|\nabla\cdot u_t\|_{s-1})
%    +K_2^{\frac12}\left(\lambda^{-1}K_1^{\frac12}+1\right)(\|\nabla u\|_s+\|\nabla\cdot u\|_s)\bigg]
    \\&\leq&
    \tau\|\nabla u_t\|_{s-1}^2+\tau\|\nabla\cdot u_t\|_{s-1}^2+C(\tau)\left(K_1^2+K_2^2\right)\|u_t\|_{s-1}^2+C\left(\|\nabla u\|_s^2+\|\nabla\cdot u\|_s^2\right).\quad\quad\
\eeq

When $s=3$, more refined estimates for the terms $I_4$, $I_7$, and $I_8$ are needed. Here we only deal with the terms with $|\beta|=s-1=2\
(D^{\alpha_1}=\nabla_i\nabla_j\partial_t)$ since for the case $|\beta|=1$, one can get the estimate in a similar manner (actually more easily),
or one can estimate the term $\int_{T^N}\left(|\lambda\nabla^\beta\rho_t|^2+|\nabla^\beta u_t|^2\right)$ by the interpolation since we have dealt with the cases
$D^{\alpha_1}=\partial_t$ and $D^{\alpha_1}=\nabla_i\nabla_j\partial_t$. We use Lemma \ref{le:2.1}, the Sobolev embedding
$H^1(T^N)\hookrightarrow L^4(T^N)$ and $H^2(T^N)\hookrightarrow L^\infty(T^N)$, (\ref{s53}) and the Cauchy inequality to give
\beq
    \label{s63}
    \nonumber
    |I_4|&\leq& \lambda\|\lambda\nabla_i\nabla_j\rho_t\|\left(\|\nabla_i\nabla_j(v\cdot\nabla \rho_t)-v\cdot\nabla\nabla_i\nabla_j\rho_t\|+\|\nabla_i\nabla_j(v_t\cdot\nabla\rho)\|\right.
    \quad\quad\quad\quad\quad\quad\quad\quad
    \\&&\nonumber
    \left.+\|\nabla_i\nabla_j(\xi\nabla\cdot u_t)-\xi\nabla\cdot\nabla_i\nabla_j u_t\|+\|\nabla_i\nabla_j(\xi_t\nabla\cdot u)\|\right)
    \\&\leq& \nonumber
    \lambda\|\lambda\rho_t\|_2\left(\|\nabla_i\nabla_j v\cdot\nabla\rho_t+\nabla_j v\cdot\nabla\nabla_i\rho_t+\nabla_i v\cdot\nabla\nabla_j\rho_t\|+\|\nabla^2(v_t\cdot\nabla\rho)\|\right.
    \\&&\nonumber
    \left.+\|\nabla_i\nabla_j\xi(\nabla\cdot u_t)+\nabla_j\xi\nabla_i(\nabla\cdot u_t)+\nabla_i\xi\nabla_j(\nabla\cdot u_t)\|+\|\nabla^2(\xi_t\cdot\nabla u)\|\right)
    \\&\leq& \nonumber
    C\lambda\|\lambda\rho_t\|_2\left(\|\nabla^2 v\|_{L^4}\|\nabla\rho_t\|_{L^4}+\|\nabla v\|_{\infty}\|\nabla^2\rho_t\|+\|v_t\|_{\infty}\|\nabla \rho\|_2+\|\nabla\rho\|_{\infty}\|v_t\|_2\right.
    \\&& \nonumber
    \left.+\|\nabla^2 \xi\|_{L^4}\|\nabla\cdot u_t\|_{L^4}+\|\nabla\xi\|_\infty\|\nabla(\nabla\cdot u_t)\|+\|\xi_t\|_\infty\|\nabla u\|_2+\|\nabla u\|_\infty\|\xi_t\|_2\right)
    \\&\leq& \nonumber
    C\lambda\|\lambda\rho_t\|_2\left(\|\nabla^2 v\|_1\|\nabla\rho_t\|_1+\|\nabla v\|_{2}\|\nabla^2 \rho_t\|+\|v_t\|_{2}\|\nabla\rho\|_2+\|\nabla \rho\|_{2}\|v_t\|_2\right.
    \\&& \nonumber
    \left.+\|\nabla^2 \xi\|_{1}\|\nabla\cdot u_t\|_{1}+\|\nabla\xi\|_2\|\nabla(\nabla\cdot u_t)\|+\|\xi_t\|_2\|\nabla u\|_2+\|\nabla u\|_2\|\xi_t\|_2\right)
%    \\&\leq& \nonumber
%    CK_1^{\frac12}\|\lambda\rho_t\|_2^2+CK_2^{\frac12}\|\lambda\rho_t\|_2+CK_1^{\frac12}\|\lambda\rho_t\|_2\|u_t\|_2 +CK_2^{\frac12}\|\lambda\rho_t\|_2
    \\&\leq&
    C(K_1+K_2)\left(\|\lambda\rho_t\|_2^2+\|u_t\|_2^2\right)+C.
\eeq

Clearly, we have
\beq
    \label{s66}
    \nonumber
    I_7&=& -\lambda^2\int_{T^N}\xi\left\{\left[\nabla_i\nabla_j\left(\frac{P'(\xi)}{\xi}\nabla\rho_t\right)-\frac{P'(\xi)}{\xi}\nabla\nabla_i\nabla_j \rho_t\right]+\nabla_i\nabla_j\left[\left(\frac{P'(\xi)}{\xi}\right)_t\nabla\rho\right]\right\}\cdot\nabla_i\nabla_j u_t
    \\&=&
    I_{7,1}+I_{7,2},
\eeq
and then we get
\beq
    \label{s67}
    \nonumber
    |I_{7,1}|&\leq& C\lambda^2\|u_t\|_2\left(\|\nabla\xi\|_\infty^2\|\nabla\rho_t\|+\|\nabla^2\xi\|_{L^4}\|\nabla\rho_t\|_{L^4}+\|\nabla\xi\|_\infty\|\nabla^2\rho_t\|\right)
    \\&\leq& \nonumber
    C\|u_t\|_2\left(\|\lambda\nabla\xi\|_2^2\|\lambda\nabla\rho_t\|+\|\lambda\nabla^2\xi\|_1\|\lambda\nabla\rho_t\|_1 +\|\lambda\nabla\xi\|_2\|\lambda\nabla^2\rho_t\|\right)
    \\&\leq&
    CK_1\left(\|\lambda\rho_t\|_2^2+\|u_t\|_2^2\right),
\eeq
where we have used the fact that
\beq
    \label{s68}
    \nonumber
    &&\nabla_i\nabla_j\left(\frac{P'(\xi)}{\xi}\nabla\rho_t\right)-\frac{P'(\xi)}{\xi}\nabla\nabla_i\nabla_j\rho_t
    \\&=& \nonumber
    \left(\frac{P'''(\xi)}{\xi}-\frac{2P''(\xi)}{\xi^2}+\frac{2P'(\xi)}{\xi^3}\right)\nabla_i\xi\nabla_j\xi\nabla\rho_t+\frac{\xi P''(\xi)-P'(\xi)}{\xi^2}\nabla_i\nabla_j\xi\nabla\rho_t
    \\&&
    +\frac{\xi P''(\xi)-P'(\xi)}{\xi^2}\nabla_j\xi\nabla\nabla_i\rho_t+\frac{\xi P''(\xi)-P'(\xi)}{\xi^2}\nabla_i\xi\nabla\nabla_j\rho_t.
\eeq

Moreover, we have from Lemma \ref{le:2.1}, (\ref{s59}) and (\ref{s66}) that
\beq
    \label{s69}
    \nonumber
    |I_{7,2}|&\leq& C\|u_t\|_2\left(\left\|\lambda\xi_t\right\|_\infty\|\lambda\nabla\rho\|_2 +\|\lambda\nabla\rho\|_\infty\left\|\lambda\xi_t\right\|_2(1+\|\nabla\xi\|_1)\right)
%    \\&\leq&\nonumber
%    C\|u_t\|_2\|\lambda\nabla\rho\|_2\left(\|\lambda\xi_t\|_2+\lambda^{-1}\|\lambda\xi_t\|_2\|\lambda\nabla\xi\|_1\right)
    \\&\leq&
    C\left(K_1^2+K_2^2\right)\|u_t\|_2^2+C.
\eeq

By combining (\ref{s66}), (\ref{s67}) and (\ref{s69}), we have
\beq
    \label{s70}
    |I_7|\leq C\left(K_1^2+K_2^2\right)\left(\|\lambda\rho_t\|_2^2+\|u_t\|_2^2\right)+C.
\eeq

Finally, we turn to give the estimate of $I_8$.
\beq
    \label{s71}
    \nonumber
    |I_8|&\leq&\nonumber
    C\|u_t\|_2\left(\|\nabla\xi\|_\infty^2\|\Delta u_t\|+\|\nabla^2\xi\|_{L^4}\|\Delta u_t\|_{L^4}+\|\nabla\xi\|_\infty\|\nabla\Delta u_t\|+\left\|\xi_t\right\|_\infty\|\Delta u\|_2\right.
    \\&&\nonumber
    \left.+\|\Delta u\|_\infty\|\xi_t\|_2(1+\|\nabla\xi\|_1)+\|\nabla\xi\|_\infty^2\|\nabla(\nabla\cdot u_t)\|+\|\nabla^2\xi\|_{L^4}\|\nabla(\nabla\cdot u_t)\|_{L^4}\right.
    \\&&\nonumber
    \left.+\|\nabla\xi\|_\infty\|\nabla^2(\nabla\cdot u_t)\|
    +\left\|\xi_t\right\|_\infty\|\nabla(\nabla\cdot u)\|_2
    +\|\nabla(\nabla\cdot u)\|_\infty\|\xi_t\|_2(1+\|\nabla\xi\|_1)\right.)
    \\&\leq&\nonumber
    C\|u_t\|_2\left(\lambda^{-2}\|\lambda\nabla\xi\|_2^2\left(\|\nabla u_t\|_1+\|\nabla\cdot u_t\|_1\right)
    +\lambda^{-1}\|\lambda\nabla\xi\|_2\left(\|\nabla u_t\|_2+\|\nabla\cdot u_t\|_2\right)\right.
    \\&& \nonumber
    \left.
    +\lambda^{-1}\|\lambda\xi_t\|_2(1+\lambda^{-1}\|\lambda\nabla\xi\|_1)\left(\|\nabla u\|_3+\|\nabla\cdot u\|_3\right)\right)
    \\&\leq&
    \tau\|\nabla u_t\|_2^2+\tau\|\nabla\cdot u_t\|_2^2+C(\tau)\left(K_1^2+K_2^2\right)\|u_t\|_2^2
    +C\left(\|\nabla u\|_3^2+\|\nabla\cdot u\|_3^2\right),\quad\quad\quad
\eeq
where we have used the fact that
\beq
    \label{s72}
    \nonumber &&\nabla_i\nabla_j\left(\frac{1}{\xi}\Delta u_t\right)-\frac{1}{\xi}\Delta \nabla_i\nabla_j u_t
    \\&=& \frac{2}{\xi^3}\nabla_i\xi\nabla_j\xi\Delta u_t-\frac{1}{\xi^2}\nabla_i
    \nabla_j\xi\Delta u_t-\frac{1}{\xi^2}\nabla_j\xi\nabla_i\Delta u_t-\frac{1}{\xi^2}\nabla_i\xi\nabla_j\Delta u_t,
\eeq and similar calculations on the term
$\nabla_i\nabla_j\left(\frac{1}{\xi}\nabla\left(\nabla\cdot
u_t\right)\right)-\frac{1}{\xi}\nabla_i\nabla_j\nabla(\nabla\cdot
u_t)$ as in (\ref{s72}).

Now putting the estimates (\ref{xin26})--(\ref{xin28}), (\ref{xin30})--(\ref{xin36}), (\ref{s57})--(\ref{s71}) together, summing over $\beta$, and then choosing $\tau$
small enough, we have
\beq
    \label{s77}
    \nonumber&&\frac{{\rm d}}{{\rm d}t}\sum\limits_{|\beta|\leq s-1}\int_{T^N}\left(\frac{P'(\xi)}{\xi}|\lambda\nabla^\beta \rho_t|^2+\xi|\nabla^\beta u_t|^2\right)
    \\&&+\nonumber
    \sum\limits_{|\beta|\leq s-1}\mu\int_{T^N}|\nabla^\beta\nabla u_t|^2+\sum\limits_{|\beta|\leq s-1}(\kappa+\mu)\int_{T^N}|\nabla^\beta(\nabla\cdot u_t)|^2
    \\&\leq&
    C\left(K_1^2+K_2^2\right)\left(\|\lambda\rho_t\|_{s-1}^2+\|u_t\|_{s-1}^2\right)
    +C\left(\|\nabla n_t\|_s^2+\|\nabla u\|_s^2+\|\nabla\cdot u\|_s^2+1\right).
\eeq

Recalling the constraints of the initial data and (\ref{l42}), we have
\beq
    \label{s78}
    \nonumber &&\|\lambda\nabla^\beta\rho_t(x,0)\|^2+\|\nabla^\beta u_t(x,0)\|^2
    \\&\leq& \nonumber
    C\left(\|\lambda((u_0+\overline{u}_0^\lambda)\cdot\nabla)\overline{\rho}_0^\lambda\|_{s-1}^2 +\|\lambda(\overline{\rho}_0^\lambda+1)\nabla\cdot\overline{u}_0^\lambda\|_{s-1}^2 +\|(1+\overline{\rho}_0^\lambda)^{-1}\lambda^2\nabla\overline{\rho}_0^\lambda\|_{s-1}^2\right.
    \\&&+ \nonumber
    \left.\|((u_0+\overline{u}_0^\lambda)\cdot\nabla)(u_0+\overline{u}_0^\lambda)\|_{s-1}^2+\|(1+\overline{\rho}_0^\lambda)^{-1}(\Delta (u_0+\overline{u}_0^\lambda)+\nabla(\nabla\cdot\overline{u}_0^\lambda))\|_{s-1}^2\right.
    \\&&\nonumber+
    \left.\|(1+\overline{\rho}_0^\lambda)^{-1}\Delta(n_0+\overline{n}_0^\lambda)\cdot\nabla(n_0+\overline{n}_0^\lambda)\|_{s-1}^2\right)
    \\&\leq& C.
\eeq

Then by the Gronwall inequality, (\ref{xin24}) and (\ref{s54}), we get
\beq
    \label{s80}
    \|\lambda\rho_t\|_{s-1}^2(t)+\|u_t\|_{s-1}^2(t)+\mu\int_0^t\|\nabla u_t\|_{s-1}^2+(\kappa+\mu)\int_0^t\|\nabla\cdot u_t\|_{s-1}^2\leq C
\eeq
for $t\in[0,T_0]$, provided that $T_0$ is small enough such that $T_0<T_3\doteq\min\left\{T_2,(K_1^2+K_2^2)^{-1}\right\}$.

It remains to show the first inequality of (\ref{s1}). It suffices to show $\|\lambda(\rho-1)\|_s+\|u-u_0\|_s+\|\nabla
n-\nabla n_0\|_{s}<c_0^{-1}\delta$ by the Sobolev inequality, where $c_0$ is the Sobolev constant. Let
$\overline{\rho}=\rho-1$, $\overline{u}=u-u_0$ and
$\nabla\overline{n}=\nabla n-\nabla n_0$.

By the Cauchy inequality with parameter $\tau$, we proceed as the proof of (\ref{s4}),
(\ref{s9}) and (\ref{s10}) for the case $D^{\alpha_2}=\nabla_i\nabla^{\alpha_1}$ $(|\alpha_1|\leq s)$ to get
\beq
    \label{s81}
    \nonumber
    &&\frac{{\rm d}}{{\rm d}t}\sum\limits_{|\alpha_1|\leq s}\|\nabla\nabla^{\alpha_1}\overline{n}\|^2+\theta\sum\limits_{|\alpha_1|\leq s}\|\nabla^2\nabla^{\alpha_1} \overline{n}\|^2
    \\&\leq&
    C(\tau)K_1^3\|\nabla\overline{n}\|_s^2+\tau K_1^{-1}\left(\|\nabla m\|_{s+1}^2+\|\nabla v\|_s^2\right)+C\|\Delta n_0\|_{s}^2.
\eeq
%where we have dealt with the term such as $\int_{T^N}\nabla_i\nabla^\beta(v\cdot \nabla
%n_0)\cdot\nabla_i\nabla^\beta\Delta\overline{n}$ by integration by parts and the Cauchy inequality
%as follows:
%\beq
%    \label{s82}
%    \nonumber
%    \left|\int_{T^N}\nabla_i\nabla^\beta(v\cdot \nabla n_0)\cdot\nabla_i\nabla^\beta\Delta \overline{n}\right|&=&\left|-\int_{T^N}\nabla^\beta(v\cdot \nabla n_0)\cdot\Delta\nabla^\beta\Delta \overline{n}\right|
%    \\&\leq&\frac{\theta}{12}\|\Delta\nabla^\beta\overline{n}\|^2+C\|\nabla^\beta(v\cdot\nabla n_0)\|^2,
%\eeq
%and the terms $\int_{T^N}\nabla_i\nabla^\beta(\Delta
%n_0)\cdot\nabla_i\nabla^\beta\Delta\overline{n}$ and
%$\int_{T^N}\nabla_i\nabla^\beta(|\nabla
%m|^2n_0)\cdot\nabla_i\nabla^\beta\Delta\overline{n}$ were
%treated by the same method shown in (\ref{s82}).

Since (\ref{l48}) implies that $\|\nabla\overline{n}(x,0)\|_s^2\leq \lambda^{-2}\delta_0^2$, we conclude from the Gronwall inequality that
\beq
    \label{s83}
    \|\nabla\overline{n}\|_{s}^2\leq \exp{\left(C(\tau)K_1^3T_0\right)}\left(\lambda^{-2}\delta_0^2+C\tau+CT_0\right)<c_0^{-1}\delta,
\eeq
where we have chosen $\lambda^{-1}$, $\tau$ and $T_0(<T_3)$  sufficiently small such that (\ref{s83}) holds.

Similarly as in the proof of (\ref{s51}) and (\ref{s81}), we have
\beq
    \label{s85}
    \nonumber
    &&\frac{{\rm d}}{{\rm d}t}\sum_{|\alpha_1|\leq s}\int_{T^N}\left(\frac{P'(\xi)}{\xi}|\lambda\nabla^{\alpha_1} \overline{\rho}|^2+\xi|\nabla^{\alpha_1}\overline{u}|^2\right)
    \\&& \nonumber
    +\sum\limits_{|\alpha_1|\leq s}\mu\int_{T^N}|\nabla\nabla^{\alpha_1}\overline{u}|^2
    +\sum\limits_{|\alpha_1|\leq s}(\kappa+\mu)\int_{T^N}|\nabla^{\alpha_1} (\nabla\cdot\overline{u})|^2
    \\&\leq&\nonumber
    C\left(K_1+K_2\right)\left(\|\lambda\overline{\rho}\|_s^2+\|\overline{u}\|_s^2\right)+C+C\|(v\cdot\nabla) u_0\|_s^2+C\|\nabla u_0\|_s^2
    \\&&\nonumber
    +C\sum_{|\alpha_1|\leq s}\left\|\nabla^{\alpha_1}\left(\frac{1}{\xi}\Delta u_0\right)-\frac{1}{\xi}\nabla^{\alpha_1}\Delta u_0\right\|^2
    \\&\leq&
    C\left(K_1+K_2\right)\left(\|\lambda\overline{\rho}\|_s^2+\|\overline{u}\|_s^2\right)+CK_1.
\eeq

Then it follows from the Gronwall inequality, (\ref{l45}) and (\ref{l48}) that
\beq
    \label{s86}
    \|\lambda\overline{\rho}\|_s^2+\|\overline{u}\|_s^2\leq C\exp{\left(C(K_1+K_2)T_0\right)}\left(\lambda^{-2}\delta_0^2+K_1T_0\right)<c_0^{-1}\delta,
\eeq
where we have chosen $T_0(<T_3)$ and $\lambda^{-1}$ small enough such that (\ref{s83}) and (\ref{s86}) hold.

As a conclusion, we have the following lemma.

\begin{Lemma}\label{le:2.2} Suppose that $B_{T_0}^\lambda(U_0)$ is
defined by (\ref{s1}) and $\Lambda:V\rightarrow U$ is defined by the
system (\ref{s2}). Then, under the assumptions in Theorem \ref{th.1.1}, there exist constants $T_0$, $\delta$,
$K_1$ and $K_2$ independent
of $\lambda$ such that $\Lambda$ maps $B_{T_0}^\lambda(U_0)$ into
itself.
\end{Lemma}

Now we plan to show that $\Lambda$ is a contractive map. In the proof of this part and the following lemmas in this section, denote by $C$ the constant
depending on the initial data, the domain, $N$, $s$ and the viscosity coefficients $\mu$, $\kappa$, $\nu$, and $\theta$, but independent of $\lambda$.

Let $U=\Lambda(V)$ and
$\widehat{U}=\Lambda(\widehat{V})$, where $V,\ \widehat{V}\in
B^\lambda_{T_0}(U_0)$. Then by the definition, we have \beq \label{s87}
\begin{cases}
         (\rho-\widehat{\rho})_t+(v\cdot\nabla)(\rho-\widehat{\rho})+\left((v-\widehat{v})\cdot\nabla\right)\widehat{\rho} +\xi\nabla\cdot(u-\widehat{u})+(\xi-\widehat{\xi})\nabla\cdot \widehat{u}=0,\\
         (u-\widehat{u})_t+(v\cdot \nabla) (u-\widehat{u})+\left((v-\widehat{v})\cdot\nabla\right)\widehat{u}
         +\lambda^2\frac{P'(\xi)}{\xi}\nabla(\rho-\widehat{\rho})
         +\lambda^2\left(\frac{P'(\xi)}{\xi}-\frac{P'(\widehat{\xi})}{\widehat{\xi}}\right)\nabla\widehat{\rho}
         \\=
         \frac{\mu}{\xi}
         \Delta(u-\widehat{u})+\left(\frac{\mu}{\xi}-\frac{\mu}{\widehat{\xi}}\right)\Delta
         \widehat{u}+\frac{\kappa+\mu}{\xi}\nabla
         \left(\nabla\cdot
         (u-\widehat{u})\right)+\left(\frac{\kappa+\mu}{\xi}-\frac{\kappa+\mu}{\widehat{\xi}}\right)\nabla(\nabla\cdot
         \widehat{u})
         \\ \ \ -
         \frac{\nu}{\xi}\left(\Delta (n-\widehat{n})\cdot\nabla n +\Delta \widehat{n}\cdot\nabla(n-\widehat{n})\right)-\left(\frac{\nu}{\xi}
        -\frac{\nu}{\widehat{\xi}}\right)\Delta\widehat{n}\cdot\nabla\widehat{n},\\
        (n-\widehat{n})_t-\theta\Delta(n-\widehat{n})
        =-({v}\cdot\nabla)(m-\widehat{m})-\left((v-\widehat{v})\cdot\nabla\right)\widehat{m}+\theta|\nabla m|^2(m-\widehat{m})
        \\
        \quad\quad\quad\quad\quad\quad\quad\quad\quad\quad\ +\theta\left((\nabla m-\nabla\widehat{m}):(\nabla m+\nabla\widehat{m})\right)\widehat{m}.%\\[2mm]
\end{cases}\eeq

Firstly, multiplying (\ref{s87})$_3$ by $(n-\widehat{n})$ and
$\Delta(n-\widehat{n})$ respectively and then integrating by parts, one obtains
\beq
    \label{s88}
    \nonumber&&\frac12\frac{{\rm d}}{{\rm d}t}\int_{T^N}|n-\widehat{n}|^2+\theta\int_{T^N}|\nabla(n-\widehat{n})|^2
    \\\nonumber&=&
    -\int_{T^N}(v\cdot\nabla)(m-\widehat{m})\cdot(n-\widehat{n})-\int_{T^N}\left((v-\widehat{v})\cdot\nabla\right)\widehat{m}\cdot(n-\widehat{n})
    \\&&
    +\theta\int_{T^N}|\nabla m|^2(m-\widehat{m})\cdot(n-\widehat{n})+\theta\int_{T^N}\left((\nabla m+\nabla \widehat{m}):(\nabla m-\nabla \widehat{m})\right)\widehat{m}\cdot(n-\widehat{n}),\quad\quad\quad
\eeq
and
\beq
    \label{s89}
    \nonumber &&
    \frac12\frac{{\rm d}}{{\rm d}t}\int_{T^N}|\nabla(n-\widehat{n})|^2+\theta\int_{T^N}|\Delta(n-\widehat{n})|^2
    \\\nonumber&=&
    \int_{T^N}(v\cdot\nabla)(m-\widehat{m})\cdot\Delta(n-\widehat{n})+\int_{T^N}\left((v-\widehat{v})\cdot\nabla\right)\widehat{m}\cdot\Delta(n-\widehat{n}) \\&& \nonumber
    -\theta\int_{T^N}|\nabla m|^2(m-\widehat{m})\cdot\Delta(n-\widehat{n})-\theta\int_{T^N}\left((\nabla m-\nabla \widehat{m}):(\nabla m+\nabla \widehat{m})\right)\widehat{m}\cdot\Delta(n-\widehat{n}).\\
\eeq

Secondly, multiplying (\ref{s87})$_1$ and (\ref{s87})$_2$ by $\lambda^2\frac{P'(\xi)}{\xi}(\rho-\widehat{\rho})$ and $\xi(u-\widehat{u})$ respectively, we have
\beq
    \label{s90}
    && \nonumber
    \frac{\lambda^2}{2}\frac{{\rm d}}{{\rm d}t}\int_{T^N}\frac{P'(\xi)}{\xi}|\rho-\widehat{\rho}|^2
    -\frac{\lambda^2}{2}\int_{T^N}|\rho-\widehat{\rho}|^2\left(\frac{P'(\xi)}{\xi}\right)_t
    \\&&\nonumber
    +\lambda^2\int_{T^N}\frac{P'(\xi)}{\xi}(\rho-\widehat{\rho})(v\cdot\nabla)(\rho-\widehat{\rho})+
    \lambda^2\int_{T^N}\frac{P'(\xi)}{\xi}(\rho-\widehat{\rho})\left((v-\widehat{v})\cdot\nabla\right)\widehat{\rho}
    \\&&
    +\lambda^2\int_{T^N}
    P'(\xi)(\rho-\widehat{\rho})\nabla\cdot(u-\widehat{u})+\lambda^2\int_{T^N}\frac{P'(\xi)}{\xi}(\rho-\widehat{\rho})(\xi-\widehat{\xi})\nabla\cdot\widehat{u}=0,
\eeq
and
\beq
    \label{s91}
    &&\nonumber
    \frac12\frac{{\rm d}}{{\rm d}t}\int_{T^N}\xi|u-\widehat{u}|^2-\frac12\int_{T^N}\xi_t|u-\widehat{u}|^2
    +\int_{T^N}\xi(v\cdot\nabla)(u-\widehat{u})\cdot(u-\widehat{u})
    \\&&+\nonumber
    \int_{T^N}\xi\left((v-\widehat{v})\cdot\nabla\right)\widehat{u}\cdot(u-\widehat{u})+\lambda^2\int_{T^N} P'(\xi)\left((u-\widehat{u})\cdot\nabla\right)(\rho-\widehat{\rho})
    \\&&\nonumber
    +\lambda^2\int_{T^N}\xi\left(\frac{P'(\xi)}{\xi} -\frac{P'(\widehat{\xi})}{\widehat{\xi}}\right)\left((u-\widehat{u})\cdot\nabla\right)\widehat{\rho}
    \\&=&\nonumber
    -\mu\int_{T^N}|\nabla(u-\widehat{u})|^2
    +\int_{T^N}\xi\left(\frac{\mu}{\xi}-\frac{\mu}{\widehat{\xi}}\right)\Delta\widehat{u}\cdot(u-\widehat{u})
    -(\kappa+\mu)\int_{T^N}|\nabla\cdot(u-\widehat{u})|^2
    \\&&\nonumber
    +\int_{T^N}\xi\left(\frac{\kappa+\mu}{\xi}-\frac{\kappa+\mu}{\widehat{\xi}}\right)\left((u-\widehat{u})\cdot\nabla\right)(\nabla\cdot\widehat{u})
    -\int_{T^N}\xi\left(\frac{\nu}{\xi}-\frac{\nu}{\widehat{\xi}}\right)\left((u-\widehat{u})\cdot\nabla\right)\widehat{n}\cdot\Delta\widehat{n}
    \\&&
    -\nu\int_{T^N}\left(\Delta(n-\widehat{n})\cdot\nabla n+\Delta\widehat{n}\cdot\nabla(n-\widehat{n})\right)\cdot(u-\widehat{u}).
\eeq

In conclusion, putting (\ref{s88})--(\ref{s91}) together, processing as before and using the Cauchy inequality, we obtain
\beq
    \label{s92}
    \nonumber
    &&\frac12\frac{{\rm d}}{{\rm d}t}\int_{T^N}\left(\frac{P'(\xi)}{\xi}|\lambda(\rho-\widehat{\rho})|^2
    +\xi|u-\widehat{u}|^2+|n-\widehat{n}|^2+|\nabla (n-\widehat{n})|^2\right)+\mu\int_{T^N}|\nabla(u-\widehat{u})|^2
    \\&&\nonumber
    +(\kappa+\mu)\int_{T^N}|\nabla\cdot(u-\widehat{u})|^2
    +\theta\int_{T^N}|\nabla(n-\widehat{n})|^2+\theta\int_{T^N}|\Delta(n-\widehat{n})|^2
    \\&\leq&
    C\left(\|\lambda(\rho-\widehat{\rho})\|^2+\|u-\widehat{u}\|^2+\|n-\widehat{n}\|^2_1+\|\lambda(\xi-\widehat{\xi})\|^2+\|v-\widehat{v}\|^2
    +\|m-\widehat{m}\|_1^2\right),
\eeq
where we have used the following estimate different from before
\beq
    \label{s93}
    \nonumber
    &&\int_{T^N}\xi\left(\frac{\mu}{\xi}-\frac{\mu}{\widehat{\xi}}\right)\Delta\widehat{u}\cdot(u-\widehat{u})
    +\int_{T^N}\xi\left(\frac{\kappa+\mu}{\xi}-\frac{\kappa+\mu}{\widehat{\xi}}\right)\left((u-\widehat{u})\cdot\nabla\right)(\nabla\cdot\widehat{u})
    \\&\leq& \nonumber
    C\|\xi-\widehat{\xi}\|\left(\|\Delta\widehat{u}\|_{L^4}+\|\nabla(\nabla\cdot\widehat{u})\|_{L^4}\right)\|u-\widehat{u}\|_{L^4}
    \\&\leq& \nonumber
    C\|\xi-\widehat{\xi}\|\|\widehat{u}\|_3(\|u-\widehat{u}\|+\|\nabla (u-\widehat{u})\|)
    \\&\leq&
    \frac{\mu}{2}\|\nabla(u-\widehat{u})\|^2+C(\lambda^{-1}\|\lambda(\xi-\widehat{\xi})\|^2+\|u-\widehat{u}\|^2).\quad
\eeq
%and
%\beq
%    \label{s94}
%    \nonumber
%    \left|\int_{T^N}\xi(\frac{\kappa+\mu}{\xi}
%    -\frac{\kappa+\mu}{\widehat{\xi}})\nabla(\nabla\cdot\widehat{u})\cdot(u-\widehat{u})\right|\leq\frac{\mu}{8}\|\nabla(u-\widehat{u})\|^2+
%    C(\|\lambda(\xi-\widehat{\xi})\|^2+\|u-\widehat{u}\|^2).
%    \\
%\eeq

Then noting that $(U-\widehat{U})(0)=0$, the Gronwall inequality and (\ref{s92}), we obtain
\beq
    \label{s95}
    && \nonumber \sup\limits_{0\leq t\leq T_0}\left(\|\lambda(\rho-\widehat{\rho})\|^2+\|u-\widehat{u}\|^2+\|n-\widehat{n}\|^2_1\right)
    \\&\leq&
    CT_0e^{CT_0}\sup\limits_{0\leq t\leq T_0}\left(\|\lambda(\xi-\widehat{\xi})\|^2+\|v-\widehat{v}\|^2+\|m-\widehat{m}\|^2_1\right).
\eeq
Moreover, we integrate (\ref{s92}) over $[0,t]\subseteq [0,T_0]$ to get
\beq
    \label{s96}
    &&\nonumber
    \int_0^t\left(\|u-\widehat{u}\|_1^2+\|n-\widehat{n}\|_2^2\right)
    \\&\leq&
    C\left(T_0^2e^{CT_0}+T_0\right)\sup\limits_{0\leq t\leq T_0}\left(\|\lambda(\xi-\widehat{\xi})\|^2+\|v-\widehat{v}\|^2+\|m-\widehat{m}\|^2_1\right).
\eeq

If we take $T_0(<T_3)$ small enough such that (\ref{s83}), (\ref{s86}) and
$C\left(T_0^2e^{CT_0}+T_0e^{CT_0}+T_0\right)<1$ are valid, where $C$ is the uniform constants mentioned in (\ref{s95})--(\ref{s96}), then we can prove the contraction.

In conclusion, we have the following lemma.
\begin{Lemma}\label{le:2.3} Under the assumptions in Theorem \ref{th.1.1}, the maps $\Lambda:V\rightarrow U$ is a contraction in
the sense that
\beq
    \label{s97}
    \nonumber
    &&\sup\limits_{0\leq t\leq T_0}\left(\|\lambda(\rho-\widehat{\rho})\|^2+\|u-\widehat{u}\|^2+\|n-\widehat{n}\|^2_1\right) +\int_0^t\left(\|u-\widehat{u}\|_1^2+\|n-\widehat{n}\|_2^2\right)
    \\&\leq&
    \eta\sup\limits_{0\leq t\leq T_0}\left(\|\lambda(\xi-\widehat{\xi})\|^2+\|v-\widehat{v}\|^2+\|m-\widehat{m}\|^2_1\right)
\eeq
for some $0<\eta<1$, provided that $T_0$ is small enough.
\end{Lemma}

In order to prove Theorem \ref{th.1.1}, we also need the following
lemma.

\begin{Lemma}\label{le:2.4} Consider the incompressible system of
liquid crystals (\ref{l40}) with the initial condition (\ref{l47})
for $s\geq 2$. Then there exists at most one strong solution $(u,n)\in \{(u,n):(u,\nabla n)\in L^\infty(0,T;H^s)\}$ for $0<T\leq\infty$.
\end{Lemma}
\begin{Proof} Assume that $(u_1,n_1)$ and $(u_2,n_2)$ are two strong solutions of (\ref{l40}) with the same initial data (\ref{l47}). Then we have
\beq \label{s98}
\begin{cases}\nabla\cdot(u_1-u_2)=0,\\
         \partial_t(u_1-u_2)+(u_1\cdot \nabla) (u_1-u_2)+\left((u_1-u_2)\cdot\nabla\right) u_2+\nabla(p_1-p_2)
         \\=\mu\Delta(u_1-u_2)-\nu\nabla\cdot\left(\nabla(n_1-n_2)\odot\nabla n_1\right)-\nu\nabla\cdot\left(\nabla n_2\odot\nabla (n_1-n_2)\right),\\
        \partial_t(n_1-n_2)+({u_1}\cdot\nabla)(n_1-n_2)+\left((u_1-u_2)\cdot\nabla\right) n_2
        \\=\theta\Delta(n_1-n_2)+\theta|\nabla n_1|^2(n_1-n_2)+\theta\left((\nabla n_1-\nabla n_2):(\nabla n_1+\nabla n_2)\right)n_2.%\\[2mm]
\end{cases}\eeq

Multiplying $(\ref{s98})_2$ by $(u_1-u_2)$ and multiplying $(\ref{s98})_3$ by $(n_1-n_2)$ and
$\Delta(n_1-n_2)$ respectively, and then using integration by parts and the Cauchy inequality, we have
\beq
    \label{s99}
    \nonumber&&
    \frac12\frac{{\rm d}}{{\rm d}t}\left(\|u_1-u_2\|^2+\|(n_1-n_2)\|_1^2\right)
    \\&& \nonumber
    +\mu\|\nabla (u_1-u_2)\|^2+\theta\|\nabla(n_1-n_2)\|^2+\theta\|\Delta(n_1-n_2)\|^2
    \\&\leq&
    M\left(\|u_1-u_2\|^2+\|n_1-n_2\|_1^2\right)+\frac{\mu}{2}\|\nabla(u_1-u_2)\|^2+\frac{\theta}{2}\|\Delta (n_1-n_2)\|^2.
\eeq
where $M\geq \sup_{t\in[0,T]}\left(\|u_i\|_2(t)+\|\nabla n_i\|_2(t)\right)$ for $i=1,2$.

Note that $\|u_1(x,0)-u_2(x,0)\|^2+\|n_1(x,0)-n_2(x,0)\|_1^2=0$. By the Gronwall inequality, we have
\beq
    \label{s100}
    \|u_1-u_2\|^2+\|n_1-n_2\|_1^2=0,
\eeq
which completes the proof of Lemma \ref{le:2.4}.
\end{Proof}

Before proving Theorem \ref{th.1.1}, we give the following lemma.
\begin{Lemma}\label{le:2.5}(\cite{Simon}) Assume that $X\subset
E\subset Y$ are Banach spaces and $X\hookrightarrow\hookrightarrow
E$. Then the following embeddings are compact: $$\D(i)\ \
\left\{\varphi:\varphi\in L^q(0,T; X),
\frac{\partial\varphi}{\partial t}\in
L^1(0,T;Y)\right\}\hookrightarrow\hookrightarrow L^q(0,T; E),\ \
{\rm if}\ \  1\leq q\leq\infty; $$ $$(ii)\ \
\left\{\varphi:\varphi\in L^\infty(0,T; X),
\frac{\partial\varphi}{\partial t}\in
L^r(0,T;Y)\right\}\hookrightarrow\hookrightarrow C([0,T]; E),\ \
{\rm if}\ \   1< r\leq\infty.
$$
\end{Lemma}

Now we are in a position to prove Theorem \ref{th.1.1}.

\begin{Proof of Theorem $2.1$} For any fixed $\lambda$, the standard procedure produces a sequence $\{(\rho_i,u_i,n_i)\}_{i=0}^\infty$ satisfying
\beqno \label{s101}
\begin{cases}
         \partial_t\rho^\lambda_{i+1}+\left(u^\lambda_i\cdot\nabla\right)\rho^\lambda_{i+1}+\rho^\lambda_i\nabla\cdot u^\lambda_{i+1}=0,\\
         \partial_tu^\lambda_{i+1}+\left(u^\lambda_i\cdot \nabla\right) u^\lambda_{i+1}+\lambda^2\frac{P'\left(\rho^\lambda_i\right)}{\rho^\lambda_i}\nabla\rho^\lambda_{i+1}
         =\frac{\mu}{\rho^\lambda_i}
         \Delta{u^\lambda_{i+1}}+\frac{(\kappa+\mu)}{\rho^\lambda_i}\nabla
         (\nabla\cdot u^\lambda_{i+1})-\frac{\nu}{\rho^\lambda_i}\left(\Delta n^\lambda_{i+1}\cdot \nabla n^\lambda_{i+1}\right),\\
        \partial_t{n}^\lambda_{i+1}-\theta\Delta{n}^\lambda_{i+1}=-\left({u^\lambda_i}\cdot\nabla\right){n}^\lambda_{i}
        +\theta|\nabla n^\lambda_i|^2{n}^\lambda_{i},%\\[2mm]
\end{cases}
\eeqno
as well as the uniform estimates
\beq \label{s102}
\begin{cases}E_s\left(U^\lambda_i(t)\right)+\|n^\lambda_i\|^2+\displaystyle\int_0^t\left(\mu\|\nabla
             u^\lambda_i\|_s^2+(\kappa+\mu)\|\nabla\cdot u^\lambda_i\|_s^2+\theta \|\nabla n^\lambda_i\|_{s+1}^2\right)\leq
             C,
             \\E_{s-1}\left(\partial_tU^\lambda_i(t)\right)+\|\partial_tn^\lambda_i\|^2+\displaystyle\int_0^t\left(\mu\|\nabla
             \partial_t
             u^\lambda_i\|_{s-1}^2+(\kappa+\mu)\|\nabla\cdot \partial_tu^\lambda_i\|_{s-1}^2
             +\theta \|\nabla \partial_tn^\lambda_i\|_{s}^2\right)\leq
             C.
             %\\[2mm]
\end{cases}
\eeq

Let
$\widehat{\rho}^\lambda_{i+1}=\rho^\lambda_{i+1}-\rho^\lambda_i$,
$\widehat{u}^\lambda_{i+1}=u^\lambda_{i+1}-u^\lambda_i$ and
$\widehat{n}^\lambda_{i+1}=n^\lambda_{i+1}-n^\lambda_i$. In view of
Lemma \ref{le:2.3}, we have \beq
    \label{s103}
    \sum\limits_{i=2}^\infty\|\widehat{\rho}^\lambda_i\|<\infty,\ \
    \sum\limits_{i=2}^\infty\left(\|\widehat{u}^\lambda_i\|+\int_0^{T_0}\|\widehat{u}^\lambda_i\|_1^2\right)<\infty,\ \
    \sum\limits_{i=2}^\infty\left(\|\widehat{n}^\lambda_i\|_1+\int_0^{T_0}\|\widehat{n}^\lambda_i\|_2^2\right)<\infty.\ \
\eeq

Let $\rho^\lambda=\rho^\lambda_1+\sum\limits_{i=2}^\infty\widehat{\rho}^\lambda_i$,
$u^\lambda=u^\lambda_1+\sum\limits_{i=2}^\infty\widehat{u}^\lambda_i$
and
$n^\lambda=n^\lambda_1+\sum\limits_{i=2}^\infty\widehat{n}^\lambda_i$,
then we have
$$\rho^\lambda_i\rightarrow \rho^\lambda,\ \ {\rm in} \ \ L^\infty([0,T_0];L^2),$$
$$u^\lambda_i\rightarrow u^\lambda,\ \ {\rm in} \ \ L^\infty([0,T_0];L^2)\cap L^2([0,T_0];H^1),$$
$$n^\lambda_i\rightarrow n^\lambda,\ \ {\rm in} \ \ L^\infty([0,T_0];H^1)\cap
L^2([0,T_0];H^2).$$
It follows obviously that
\beqno
      \rho^\lambda,u^\lambda, \nabla n^\lambda\in L^\infty([0,T_0];H^s)\cap {\rm Lip}([0,T_0];H^{s-1})
\eeqno
satisfy the estimates (\ref{s102}) according to the lower semi-continuity.

For any $s'\in[0,s)$, by the Sobolev interpolation inequalities, we have
\beq
    \label{s104}
    \nonumber
    &&\|(\rho_i^\lambda,u_i^\lambda,\nabla n_i^\lambda)-(\rho^\lambda,u^\lambda,\nabla n^\lambda)\|_{s'}
    \\&\leq& \nonumber
    C\|(\rho_i^\lambda,u_i^\lambda,\nabla n_i^\lambda)-(\rho^\lambda,u^\lambda,\nabla n^\lambda)\|^\theta(\|(\rho_i^\lambda,u_i^\lambda,\nabla n_i^\lambda)\|_s +\|(\rho^\lambda,u^\lambda,\nabla n^\lambda)\|_s)^{1-\theta}
    \rightarrow 0 \quad
\eeq
as $i\rightarrow \infty$ for some $\theta\in(0,1)$, where we have used lemma \ref{le:2.3} to get
\beq
    \label{s105}
    \|(\rho_i^\lambda,u_i^\lambda,\nabla n_i^\lambda)-(\rho^\lambda,u^\lambda,\nabla n^\lambda)\|\leq \sum\limits_{j=i+1}^\infty\|(\rho_j^\lambda,u_j^\lambda,\nabla n_j^\lambda)-(\rho_{j-1}^\lambda,u_{j-1}^\lambda,\nabla n_{j-1}^\lambda)\|\leq \frac{C\tau^i}{1-\tau}.\quad
\eeq
Hence, $(\rho^\lambda,u^\lambda,\nabla n^\lambda)\in C([0,T_0];H^{s'})$. In addition, by Lemma \ref{le:2.5} and (\ref{s102}), one deduces easily that
$(\rho^\lambda,u^\lambda,n^\lambda)$ is a strong solution of compressible liquid crystal model (\ref{l38}).

Finally, multiplying (\ref{l38})$_3$ by $n^\lambda$, we get an equation for $\left(|n^\lambda|^2-1\right)$ as
\beq
    \label{s106}
    \frac12\left(|n^\lambda|^2-1\right)_t+\frac12\left(u^\lambda\cdot\nabla\right)\left(|n^\lambda|^2-1\right) =\frac{\theta}{2}\Delta\left(|n^\lambda|^2-1\right)+\theta|\nabla n^\lambda|^2\left(|n^\lambda|^2-1\right).
\eeq
Multiplying (\ref{s106}) by $\left(|n^\lambda|^2-1\right)$, integrating over
${T^N}$, and then using the Gronwall inequality and the
assumption that $|n^\lambda(x,0)|^2=1$, we have
\beq
    \label{s107}
    \int_{T^N}\left||n^\lambda|^2-1\right|^2=0, \  t\in [0,T_0],
\eeq
Then, from the regularity of $n^\lambda$, we conclude that
\beq
    \label{s108}
    |n^\lambda|=1,\ \ \textrm{in} \ \ Q_{T_0}.
\eeq

The uniqueness can be proved by a similar argument as Lemma \ref{le:2.3}, which completes the proof of the uniform stability part of Theorem \ref{th.1.1}.

Next, we show that $(\rho^\lambda, u^\lambda, n^\lambda)$ converges
to the unique strong solution to the corresponding incompressible
system (\ref{l40}) as $\lambda\rightarrow\infty$. To see this, note
first that (\ref{l49}) implies that $\rho^\lambda\rightarrow 1$ in
$L^\infty([0,T_0];H^s)\cap {\rm Lip}([0,T_0];H^{s-1})$, and there exists a
subsequence $\{(u^{\lambda_j}, n^{\lambda_j})\}_j$ of
$\{(u^{\lambda}, n^{\lambda})\}_\lambda$ with a limit $u$ and $n$
such that
\beq \label{s109}
\begin{cases}u^{\lambda_j}\rightharpoonup u  \quad {\rm weakly}^* \ \ {\rm in}\ \
             L^\infty([0,T_0];H^s)\cap {\rm Lip}([0,T_0];H^{s-1}),
             \\
             n^{\lambda_j}\rightharpoonup n  \quad {\rm weakly}^* \ \ {\rm in}\ \
             L^\infty([0,T_0];H^{s+1})\cap {\rm Lip}([0,T_0];H^{s}),
             \\
             u^{\lambda_j}\rightarrow u \quad {\rm in}\ \ C([0,T_0];H^{s'}),
             \\
             n^{\lambda_j}\rightarrow n \quad {\rm in}\ \ C([0,T_0];H^{s'+1})
             %\\[2mm]
\end{cases}\eeq
for any $0\leq s'<s$, where we have used the fact that the embedding
$H^s\hookrightarrow H^{s'}$ is compact and Lemma \ref{le:2.5}.

Now we are to let $j\rightarrow\infty$ $({\lambda_j}\rightarrow
\infty)$ in (\ref{l38}).

First of all, multiplying (\ref{l38})$_1$ and (\ref{l38})$_3$ by two smooth test functions $\psi_1(x,t)$ and
$\psi_3(x,t)$ with compact supports in $t\in [0,T_0]$ respectively, we have
\beq
    \label{s110}&& \int_0^{T_0}\int_{T^N}\nabla\cdot u^{\lambda_j}\psi_1 =\int_0^{T_0}\int_{T^N}\left(\rho^{\lambda_j}_t+\left(u^{\lambda_j}\cdot\nabla\right)\rho^{\lambda_j} +\left(\rho^{\lambda_j}-1\right)\nabla\cdot u^{\lambda_j}\right)\psi_1,
    \\&& \label{s111}
    \int_0^{T_0}\int_{T^N}\left(n^{\lambda_j}_t+\left(u^{\lambda_j}\cdot\nabla\right)n^{\lambda_j}-\theta\left(\Delta n^{\lambda_j}+|\nabla n^{\lambda_j}|^2n^{\lambda_j}\right)\right)\psi_3=0.
\eeq
Then $(u,n)$ satisfies (\ref{l40})$_1$ and (\ref{l40})$_3$ by sending $j\rightarrow\infty$.

Let $\psi_2(x,t)$ be a smooth test function of (\ref{l38})$_2$
with compact supports in $t\in [0,T_0]$ and the divergence free
condition $\nabla\cdot\psi_2=0$. Then we have
\beq
    \label{s112}
    && \nonumber
    \int_0^{T_0}\int_{T^N}\left(u^{\lambda_j}_t+\left(u^{\lambda_j}\cdot\nabla\right)u^{\lambda_j}-\frac{\mu}{\rho^{\lambda_j}}\Delta u^{\lambda_j}-\frac{\kappa+\mu}{\rho^{\lambda_j}}\nabla\left(\nabla\cdot u^{\lambda_j}\right)\right.
    \\&&\nonumber\quad\quad\quad\quad
    \left.+\frac{\nu}{\rho^{\lambda_j}}\nabla\cdot\left(\nabla n^{\lambda_j}\odot\nabla n^{\lambda_j}\right)-\frac{\nu}{\rho^{\lambda_j}}\nabla\left(\frac{|\nabla n^{\lambda_j}|^2}{2}\right)\right)\psi_2
    \\&=&
    -\lambda_j^2\int_0^{T_0}\int_{T^N}\psi_2\nabla\left(\int_1^{\rho^{\lambda_j}(x,t)}\frac{P'(\xi)}{\xi}{\rm d}\xi\right)=0.
\eeq
Then, let $j\rightarrow \infty$ and get
\beq
    \label{s113}
    \mathbf{P}\left(u_t+(u\cdot\nabla) u-\mu\Delta u+\nu\nabla\cdot(\nabla n\odot\nabla n)-\nu\nabla\left(\frac{|\nabla n|^2}{2}\right)\right)=0,\eeq
where $\mathbf{P}$ is the $L^2$-projection on the divergence free vector fields.

If
\beq
    \label{s114}
    u_t+(u\cdot\nabla) u-\mu\Delta u+\nu\nabla\cdot(\nabla n\odot\nabla n)-\nu\nabla\left(\frac{|\nabla n|^2}{2}\right)=-\nabla\widehat{p}
\eeq
for some $\widehat{p}\in L^\infty([0,T_0];H^{s-1})\cap L^2([0,T_0];H^{s})$, then we have
\beq
    \label{s115}
    \frac{\lambda^2_j}{\rho^{\lambda_j}}\nabla\left(P(\rho^{\lambda_j})\right)\rightarrow\nabla \widehat{p}, \quad {\rm weakly}^*\ \ {\rm in} \ \ L^\infty([0,T_0];H^{s-2})\cap L^2([0,T_0];H^{s-1}).
\eeq

Taking $p=\widehat{p}-\nu\frac{|\nabla n|^2}{2}$, one sees that
(\ref{l40})$_2$ follows from (\ref{s114}) directly. Actually, Lemma
\ref{le:2.4} ensures that the convergence is in fact valid for the sequence
$(\rho^\lambda,u^\lambda,n^\lambda)$ themselves.

Moreover, similarly as in (\ref{s106})--(\ref{s108}), we have $|n|=1$ in
$Q_{T_0}$, which completes the proof of Theorem \ref{th.1.1}.
\end{Proof of Theorem $2.1$}

\begin{Remark} \label{re:2.1} It follows from (\ref{l49}) and (\ref{s115}) that
\beq
    \label{s116}
    \|\lambda^2\nabla \rho^\lambda\|_{s-2}+\int_0^t\|\lambda^2\nabla \rho^\lambda\|_{s-1}^2\leq C, \  t\in [0,T_0].
\eeq
\end{Remark} \begin{Remark} \label{re:2.2} It follows from (\ref{l38}) and (\ref{l49}) that
\beq
    \label{s117}
    \|\lambda\nabla\cdot u^\lambda\|_{s-1}\leq \|\lambda\rho^\lambda_t\|_{s-1}+\|\lambda(u^\lambda\cdot\nabla)\rho^\lambda\|_{s-1} +\|\lambda(\rho^\lambda-1)\nabla\cdot u^\lambda\|_{s-1}\leq C, \  t\in [0,T_0].\quad
\eeq
\end{Remark}

\setcounter{section}{3} \setcounter{equation}{0}
\section{Dispersive energy estimates and proof of Theorem \ref{th.1.2}}
In this section, we devote ourselves to getting {\it a} {\it priori} energy
estimates (\ref{l53}) and (\ref{l54}) with small
initial displacements and small initial data. Then, together with Theorem \ref{th.1.1}, Theorem \ref{th.1.2} can be proved by a standard procedure.

First of all, we assume that $|n|=1$ for all $(x,t)\in Q_{T^\lambda}$, and
\beq
    \label{r1}
    E_s(U(t))+\displaystyle\int_0^t\left(\mu\|\nabla u\|_s^2+(\kappa+\mu)\|\nabla\cdot u\|_s^2+\theta \|\nabla^2 n\|_{s}^2\right)\leq 4\left(\varepsilon_0+\lambda^{-2}\delta_0^2\right)
\eeq
for $t\in [0,T^\lambda]$, and then what we need to do is to prove the following desired estimates
\beq
    \label{r2}
    E_s(U(t))+\int_0^t\left(\mu\|\nabla u\|_s^2+(\kappa+\mu)\|\nabla\cdot u\|_s^2+\theta \|\nabla^2 n\|_{s}^2\right)\leq 3\left(\varepsilon_0+\lambda^{-2}\delta_0^2\right)
\eeq
for $t\in [0,T^\lambda]$. Then (\ref{l53}) follows by the standard continuity argument and the fact $E_s(U(0))< 4(\varepsilon_0+\lambda^{-2}\delta_0^2)$
(see also \cite{Lei,Lei2}).

Firstly, we plan to give (\ref{l54}) under the assumptions that
$|n|=1$ and (\ref{r1}). We go back to (\ref{s11})--(\ref{s12}),
(\ref{xin26})--(\ref{xin28}) (\ref{xin30})--(\ref{xin36}) and (\ref{s57})--(\ref{s71}), replace $(\xi,v,m)$ by $(\rho,u,n)$, make use of (\ref{r1})
and the Cauchy inequality, and then make a similar (just a little different) argument to give the following facts:
\beq
    \label{r6}
    \nonumber
    &&\frac{{\rm d}}{{\rm d}t}\sum_{|\beta|\leq s-1, |\alpha_1|\leq s}\int_{T^N}\left(\frac{P'(\rho)}{\rho}|\lambda\nabla^\beta\rho_t|^2+\rho|\nabla^\beta u_t|^2+|\nabla^{\alpha_1} n_t|^2\right)
    \\&& \nonumber
    +\sum_{|\beta|\leq s-1}\mu\int_{T^N}|\nabla^\beta\nabla u_t|^2+\sum_{|\beta|\leq s-1}(\kappa+\mu)\int_{T^N}|\nabla^\beta(\nabla\cdot u_t)|^2+\sum_{|\alpha_1|\leq s}\theta\int_{T^N}|\nabla^{\alpha_1}\nabla n_t|^2
    \\&\leq&
    C\left(\|\nabla u\|_s^2+\|\nabla\cdot u\|_s^2+1\right)\left(\|\lambda\rho_t\|_{s-1}^2+\|u_t\|_{s-1}^2+\|n_t\|_{s}^2\right),
\eeq
provided that $\varepsilon_0$ and $\lambda^{-1}$ are both small enough. Here we have used the fact $\|\rho_t\|_\infty\leq \|\nabla\cdot(\rho u)\|_\infty\leq C$
inferred from (\ref{r1}) to estimate $I_1$. Then by the Gronwall inequality, (\ref{l45}), (\ref{s78}) and (\ref{r1}), we get
\beq
    \label{r7}
    \nonumber
    &&\|\lambda\rho_t\|_{s-1}^2+\|u_t\|_{s-1}^2+\|n_t\|_{s}^2
    +\displaystyle\int_0^t\left(\mu\|\nabla u_t\|_{s-1}^2+(\kappa+\mu)\|\nabla\cdot u_t\|_{s-1}^2+\theta\|\nabla n_t\|_{s}^2\right)
    \\&\leq&
    C(1+t)\exp{Ct}
    \leq C\exp{Ct}
\eeq
for $t\in [0,T^\lambda]$, and thus we conclude that Remark \ref{re:2.1} and \ref{re:2.2} hold for $t\in [0,T^\lambda]$ by the same discussion in Section $3$.

Secondly, we shall begin to prove (\ref{r2}). We go back to (\ref{s4}), (\ref{s9}) and (\ref{s10}), replace $(v,m)$ by $(u,n)$,
let $D^{\alpha_2}=\nabla_i\nabla^{\alpha_1}$ with $|\alpha_1|\leq s$, use integration by parts and then make a different argument by noting the
fact that $n(x,t)\in S^2$ as follows:
\beq
    \label{r8}
    \nonumber &&\frac12\sum_i\frac{{\rm d}}{{\rm d}t}\int_{T^N} |\nabla_i\nabla^{\alpha_1} n|^2+\theta\sum_{i}\int_{T^N}|\nabla\nabla_i\nabla^{\alpha_1} n|^2
    \\&=&\nonumber
    \int_{T^N}\nabla^{\alpha_1} (u\cdot\nabla n)\cdot \Delta\nabla^{\alpha_1} n-\theta\int_{T^N}\nabla^{\alpha_1}\left(|\nabla n|^2n\right)\cdot \Delta\nabla^{\alpha_1} n
    \\&\leq&\nonumber
    C\|\Delta\nabla^{\alpha_1} n\|\left(\|\nabla^{\alpha_1}(u\cdot\nabla n)\|+\|\nabla^{\alpha_1}\left(|\nabla n|^2n\right)\|\right)
    \\&\leq& \nonumber
    \tau\|\Delta\nabla^{\alpha_1} n\|^2+C(\tau)\left(\|u\|_\infty^2\|\nabla n\|_s^2+\|\nabla n\|_\infty^2\|u\|_s^2\right)
    \\&& \nonumber
    +C\left(\|\nabla n\|_\infty^4\|n\|_{s}^2+\|n\|_\infty^2\|\nabla n\|_\infty^2\|\nabla n\|_s^2\right)
    \\&\leq&
    \tau\|\Delta\nabla^\alpha n\|^2+C(\tau)\left(\varepsilon_0+\lambda^{-2}\delta_0^2\right)\|\nabla^{s+2} n\|^2,
\eeq
where we have used the Poincar${\rm\acute{e}}$ inequality in the last step since
\beq
    \label{r9}
    \int_{T^N}\nabla^k n=0,\ \ k\geq 1.
\eeq

Summing over $\alpha_1$ and then choosing $\varepsilon_0$ and $\lambda^{-1}$ small enough, we conclude from (\ref{r8}) that
\beq
    \label{r10}
    \sum_{|\alpha_1|\leq s}\frac{{\rm d}}{{\rm d}t}\int_{T^N} |\nabla\nabla^{\alpha_1} n|^2
    +\theta\sum_{|\alpha_1|\leq s}\int_{T^N}|\nabla^2\nabla^{\alpha_1} n|^2\leq 0.
\eeq
Then by integrating (\ref{r10}) over $[0,t]$, we have
\beq
    \label{r11}
    \|\nabla n\|_{s}^2(t)+\theta\int_0^t\|\nabla^2 n\|_{s}^2\leq \|\nabla n_0\|_{s}^2+\lambda^{-2}\delta_0^2\leq \varepsilon_0+\lambda^{-2}\delta_0^2
\eeq for $t\in [0,T^\lambda]$.

On the other hand, recalling (\ref{s23}) and replacing $(\xi,v)$ by $(\rho,u)$, one has
\beq
    \label{r12}
    \nonumber
    &&\frac12\frac{{\rm d}}{{\rm d}t}\int_{T^N}\left(\lambda^2\frac{P'(\rho)}{\rho}|D^{\alpha_1}(\rho-1)|^2+\rho|D^{\alpha_1} u|^2\right)
    \\&&
    +\mu\int_{T^N}|\nabla D^{\alpha_1} u|^2+(\kappa+\mu)\int_{T^N}|\nabla\cdot D^{\alpha_1} u|^2=\sum_{k=1}^8 J_k,
\eeq
where there are some slight changes from $I_j$ as follows
\beqno
    J_1&=&\frac{\lambda^2}{2}\int_{T^N}\frac{2P'(\rho)-P''(\rho)\rho}{\rho}\nabla\cdot u|D^{\alpha_1}(\rho-1)|^2, \label{r13}
    \\J_2&=&\lambda^2\int_{T^N} P''(\rho)
    D^{\alpha_1}(\rho-1)D^{\alpha_1} u\cdot\nabla\rho,\label{r14}
    \\
    J_3&=&-\lambda^2\int_{T^N}\frac{P'(\rho)}{\rho}D^{\alpha_1}(\rho-1)\{[D^{\alpha_1}(u\cdot\nabla\rho)-u\cdot\nabla D^{\alpha_1}\rho]+\left[D^{\alpha_1}(\rho\nabla\cdot u)-\rho\nabla \cdot D^{\alpha_1} u\right]\},\quad\quad \label{r15}
    \\
    J_4&=&-\nu\int_{T^N}D^{\alpha_1}(\Delta n\cdot\nabla n)\cdot D^{\alpha_1} u,\label{r16}
    \\
    J_5&=&-\int_{T^N}\rho\left[D^{\alpha_1}(u\cdot\nabla u)-u\cdot\nabla D^{\alpha_1} u\right]\cdot D^{\alpha_1} u,\label{r17}
    \\
    J_6&=&-\lambda^2\int_{T^N}\rho\left[D^{\alpha_1}\left(\frac{P'(\rho)}{\rho}\nabla\rho\right)-\frac{P'(\rho)}{\rho}\nabla D^{\alpha_1} \rho\right]\cdot D^{\alpha_1} u, \label{r18}
    \\ \nonumber
    J_7&=&\int_{T^N}\rho\left[D^{\alpha_1}\left(\frac{\mu}{\rho}\Delta u\right)-\frac{\mu}{\rho}\Delta D^{\alpha_1} u\right]\cdot D^{\alpha_1} u
    \\&&
    +\int_{T^N}\rho\left\{D^{\alpha_1}\left[\frac{\kappa+\mu}{\rho}\nabla(\nabla\cdot u)\right]-\frac{\kappa+\mu}{\rho}\nabla D^{\alpha_1}(\nabla\cdot u)\right\}\cdot D^{\alpha_1} u,\label{r19}
    \\
    J_8&=&-\int_{T^N}\rho\left\{D^{\alpha_1}\left[\frac{\nu}{\rho}(\Delta n\cdot\nabla n)\right]-\frac{\nu}{\rho}D^{\alpha_1}(\Delta n\cdot\nabla n)\right\}\cdot D^{\alpha_1} u,
    \quad\quad\quad\quad\quad\quad\quad\quad\quad\quad \label{r20}
\eeqno
since we have used (\ref{l38})$_1$ to give
\beq
    \label{r21}
    \nonumber
    &&\frac{\lambda^2}{2}\int_{T^N}|D^{\alpha_1}(\rho-1)|^2\partial_t\left(\frac{P'(\rho)}{\rho}\right) +\frac{\lambda^2}{2}\int_{T^N}|D^{\alpha_1}(\rho-1)|^2\nabla\cdot\left(\frac{P'(\rho)}{\rho}u\right)
    \\&=&  \nonumber
    \frac{\lambda^2}{2}\int_{T^N}\left(\left(\frac{P'(\rho)}{\rho}\right)'\rho_t +\nabla\cdot\left(\frac{P'(\rho)}{\rho}u\right)\right)|D^{\alpha_1}(\rho-1)|^2
    \\&=&  \nonumber
    \frac{\lambda^2}{2}\int_{T^N}\left(\left(\frac{P'(\rho)}{\rho}\right)'\rho_t+\left(\frac{P'(\rho)}{\rho}\right)'\nabla\rho\cdot u+\left(\frac{P'(\rho)}{\rho}\right)'\rho \nabla\cdot u\right)|D^{\alpha_1}(\rho-1)|^2
    \\&& \nonumber
    +\frac{\lambda^2}{2}\int_{T^N}\left(\frac{P'(\rho)}{\rho}-\left(\frac{P'(\rho)}{\rho}\right)'\rho\right)\nabla\cdot u|D^{\alpha_1}(\rho-1)|^2
    \\&=&  \nonumber
    \frac{\lambda^2}{2}\int_{T^N}\left(\frac{P'(\rho)}{\rho}-\left(\frac{P'(\rho)}{\rho}\right)'\rho\right)\nabla\cdot u|D^{\alpha_1}(\rho-1)|^2
    \\&=&
    \frac{\lambda^2}{2}\int_{T^N}\frac{2P'(\rho)-P''(\rho)\rho}{\rho}\nabla\cdot u|D^{\alpha_1}(\rho-1)|^2.
\eeq

Now we give the dispersive estimates about $J_k$ for $k=1,2,\ldots,8$, when $D^{\alpha_1}=\nabla^{\alpha_1}$ with $|\alpha_1|\leq s$, $s\geq 3$.

Firstly, Remark \ref{re:2.2} implies that
\beq
    \label{r22}
    |J_1|&\leq& C\lambda^{-1}\|\lambda\nabla\cdot u\|_\infty\|\lambda \nabla^{\alpha_1}(\rho-1)\|^2
    \leq C\lambda^{-1}\|\lambda\nabla^{\alpha_1}(\rho-1)\|^2.
\eeq

Furthermore, it follows from the H\"older inequality, the Sobolev embedding and Remark \ref{re:2.1} that
\beq
    \label{r23}
    \nonumber
    |J_2|&\leq& C\lambda\|\nabla \rho\|_{L^4}\|\nabla^{\alpha_1} u\|_{L^4}\|\lambda \nabla^{\alpha_1}(\rho-1)\|
%    \\&\leq& \nonumber
%    C\lambda\|\nabla \rho\|_1\left(\|\nabla^{\alpha_1} u\|+\|\nabla\nabla^{\alpha_1} u\|\right)\|\lambda \nabla^{\alpha_1}(\rho-1)\|
    \\&\leq& \nonumber
    C\lambda^{-1}\|\lambda^2\nabla\rho\|_{s-2}\left(\|\nabla^{\alpha_1} u\|+\|\nabla\nabla^{\alpha_1} u\|\right)\|\lambda \nabla^{\alpha_1}(\rho-1)\|
    \\&\leq&
    \tau\|\nabla\nabla^{\alpha_1} u\|^2+C(\tau)\lambda^{-1}(\|\nabla^{\alpha_1} u\|^2+\|\lambda\nabla^{\alpha_1}(\rho-1)\|^2).
\eeq

Then we use the H\"older inequality, Lemma \ref{le:2.1}, the Cauchy inequality and the Poincar${\rm\acute{e}}$ inequality to give
\beq
    \nonumber
    |J_3|&\leq&
    C\lambda \|\lambda \nabla^{\alpha_1}(\rho-1)\|\left(\|\nabla u\|_\infty\|\nabla \rho\|_{s-1}+\|\nabla\rho\|_\infty\|u\|_s+\|\nabla\rho\|_\infty\|\nabla\cdot u\|_{s-1}\right.
    \\&& \nonumber
    \left.+\|\nabla\cdot u\|_\infty\|\nabla\rho\|_{s-1}\right)
    \\&\leq& \nonumber
    C\lambda^{-1}\|\lambda\nabla^{\alpha_1}(\rho-1)\|\|u\|_{s}\|\lambda^2\nabla\rho\|_{s-1}
    \\&\leq&
    C\lambda^{-1}\|\lambda(\rho-1)\|_s^2+C\lambda^{-1}\|\lambda^2\nabla\rho\|_{s-1}^2,\label{r24}
    \\ \nonumber
    |J_4|&=&\left|\int_{T^N}\nabla^{\alpha_1}\left(\nabla n\odot\nabla n-\frac{|\nabla n|^2}{2}I_N\right)\cdot\nabla\nabla^{\alpha_1} u\right|
    \\&\leq&\nonumber
    \tau\|\nabla\nabla^{\alpha_1} u\|^2+C(\tau)\|\nabla n\|_\infty^2\|\nabla n\|_s^2
    \\&\leq&
    \tau\|\nabla\nabla^{\alpha_1} u\|^2+C(\tau)\left(\varepsilon_0^2+\lambda^{-4}\delta_0^2\right)\|\nabla^2 n\|_s^2,
    \\
    |J_5|&\leq&
    C\|\nabla^{\alpha_1} u\|\|\nabla u\|_\infty\|\nabla u\|_{s-1}\leq C\left(\varepsilon_0^{\frac12}+\lambda^{-1}\delta_0\right)\|\nabla u\|_{s}^2,\label{r25}
    \\ \nonumber
    |J_6|&\leq&
    C\lambda^2\|\nabla^{\alpha_1} u\|\left(\left\|\nabla\rho\right\|_\infty\|\nabla\rho\|_{s-1}+\|\nabla\rho\|_\infty\left\|\nabla
    \rho\right\|_{s-1}\right)
    \\&\leq&  \nonumber
    C\lambda^{2}\|\nabla^{\alpha_1} u\|\|\nabla\rho\|_2\|\nabla\rho\|_{s-1}
    \\&\leq&  \nonumber
    C\lambda^{-1}\|\nabla^{\alpha_1} u\|\|\lambda^2\nabla\rho\|_{s-1}\|\lambda\nabla\rho\|_{s-1}
    \\&\leq&
    C\lambda^{-1}\|u\|_s^2+C\lambda^{-1}\|\lambda^2\nabla\rho\|_{s-1}^2,\label{r26}
    \\ \nonumber
    |J_7|&\leq&
    C\|\nabla^{\alpha_1} u\|\left.(\left\|\nabla\rho\right\|_\infty\|\Delta u\|_{s-1}+\|\Delta u\|_\infty\left\|\nabla\rho\right\|_{s-1}+\left\|\nabla\rho\right\|_\infty
    \|\nabla(\nabla\cdot u)\|_{s-1}\right.
    \\&& \nonumber
    \left.+\left\|\nabla(\nabla\cdot u)\|_\infty\|\nabla\rho\right\|_{s-1}\right)
    \\&\leq& \nonumber
    C\lambda^{-1}\|u\|_s\|\lambda\nabla\rho\|_{s-1}\|\nabla u\|_{s}
    \\&\leq&
    C\lambda^{-1}\left(\|u\|_s^2+\|\nabla u\|_{s}^2\right),\label{r27}
    \\ \nonumber
    |J_8|&\leq&
    C\|\nabla^{\alpha_1} u\|\left(\left\|\nabla\rho\right\|_\infty\|\Delta n\cdot\nabla n\|_{s-1}
    +\|\Delta n\|_\infty\|\nabla n\|_\infty\left\|\nabla\rho\right\|_{s-1}\right)
    \\&\leq& \nonumber
    C\|\nabla^{\alpha_1} u\|\left(\|\nabla\rho\|_2(\|\Delta n\|_\infty\|\nabla n\|_{s-1}+\|\nabla n\|_\infty\|\Delta n\|_{s-1})+\|\Delta n\|_2\|\nabla n\|_2\|\nabla\rho\|_{s-1}\right)
    \\&\leq&
    C\lambda^{-1}(\|u\|_s^2+\|\lambda\nabla\rho\|_{s-1}^2).\label{r28}
\eeq

In conclusion, we put (\ref{r22})--(\ref{r28}) together, sum over $\alpha_1$, and then choose the parameter $\tau$, $\varepsilon_0$ and $\lambda^{-1}$ small enough
to give the following dispersive energy estimate
\beq
    \label{r31}
    \nonumber
    &&\frac12\frac{{\rm d}}{{\rm d}t}\sum_{|{\alpha_1}|\leq s}\int_{T^N}\left(\lambda^2\frac{P'(\rho)}{\rho}|\nabla^{\alpha_1}(\rho-1)|^2+\rho|\nabla^{\alpha_1}
    u|^2\right)
    \\&& \nonumber
    +\mu\sum_{|{\alpha_1}|\leq s}\int_{T^N}|\nabla^{\alpha_1}\nabla u|^2+(\kappa+\mu)\sum_{|\alpha_1|\leq s}\int_{T^N}|\nabla^{\alpha_1}(\nabla\cdot u)|^2
    \\&\leq&
    C\lambda^{-1}\left(\|u\|_s^2+\|\lambda(\rho-1)\|_s^2\right)+C\lambda^{-1}\|\lambda^2\nabla\rho\|_{s-1}^2 +C\left(\varepsilon_0+\lambda^{-2}\right)\|\nabla^2 n\|_s^2.
\eeq

Then we can use the Gronwall inequality, Remark \ref{re:2.1} and (\ref{r11}) to give
(\ref{r2}) for $0\leq t\leq T^\lambda$ with $T^\lambda=\lambda^{1-\delta}$ ($\delta<1$ is a small positive constant), provided that $\lambda^{-1}$ and $\varepsilon_0$
are both small enough.

Moreover, $|n|=1$ in $Q_{T^\lambda}$ can be proved by repeating the procedure shown in (\ref{s106})--(\ref{s108}).

Finally, by (\ref{l53}) and (\ref{l54}), we can proceed as in Section 3 by taking smooth test functions to prove that
the limiting function $(1, u, n)$ and $\nabla p$ satisfy the incompressible system (\ref{l40})
in the time interval $[0,T]$ with the initial data (\ref{l51}) satisfying the constraints (\ref{l47}) and (\ref{l52}).
Since $T$ is an arbitrary positive constant, we have in fact obtained a unique global strong solution of (\ref{l40}) and
proved Theorem \ref{th.1.2}.

\setcounter{section}{4} \setcounter{equation}{0}
\section{The convergence rates about $u^\lambda$ and $n^\lambda$ when $\lambda\rightarrow \infty$}
In this section, we will prove Theorem \ref{th.1.3} by the modulated energy
method with the help of uniform estimates (\ref{l49}).

\begin{Proof of Theorem $2.3$}
First of all, let us rewrite (\ref{l38}) as follows
\beq
    \label{4.1}
    \begin{cases}\rho^\lambda_t+{\rm{div}}(\rho^\lambda u^\lambda)=0,\\
    (\rho^\lambda u^\lambda)_t+\nabla\cdot(\rho^\lambda u^\lambda\otimes u^\lambda)+\lambda^2\nabla P(\rho^\lambda)=\mu
    \Delta{u^\lambda}+(\kappa+\mu)\nabla(\nabla\cdot u^\lambda)-\nu(\Delta n^\lambda\cdot\nabla n^\lambda),\\
    {n^\lambda_t}+({u^\lambda}\cdot\nabla){n^\lambda}=\theta(\Delta{n^\lambda}+|\nabla n^\lambda|^2{n^\lambda}).%\\[2mm]
\end{cases}
\eeq

Multiplying (\ref{4.1})$_2$ by $u^\lambda$, and then integrating over $T^N$ and using integration by parts, we have
\beq
    \label{4.2}\nonumber
    &&\frac{{\rm d}}{{\rm d}t}\int_{T^N} \left(\frac12\rho^\lambda|u^\lambda|^2+\lambda^2 \omega(\rho^\lambda)\right)+\mu\int_{T^N}|\nabla u^\lambda|^2+(\kappa+\mu)\int_{T^N}|\nabla\cdot u^\lambda|^2
    \\&&=
    -\nu\int_{T^N}(u^\lambda\cdot\nabla)n^\lambda\cdot\Delta n^\lambda,
\eeq
where $\omega(\rho)=\rho\int_1^\rho\frac{P(z)}{z^2}dz$.

Multiplying (\ref{4.1})$_3$ by $\nu(\Delta n^\lambda+|\nabla
n^\lambda|^2n^\lambda)$, and then integrating over $T^N$ and noting that $|n^\lambda|=1$,
we have
\beq
    \label{4.3}
    \frac{{\rm d}}{{\rm d}t}\int_{T^N}\frac{\nu}{2}|\nabla n^\lambda|^2+\nu\theta\int_{T^N}|\Delta n^\lambda+|\nabla n^\lambda|^2n^\lambda|^2=\nu\int_{T^N}(u^\lambda\cdot\nabla)n^\lambda\cdot\Delta n^\lambda.
\eeq

Putting (\ref{4.2}) and (\ref{4.3}) together, we obtain
\beq
    \label{4.4}
    &&\nonumber
    \frac{{\rm d}}{{\rm d}t}\int_{T^N} \left(\frac12\rho^\lambda|u^\lambda|^2+\frac{\nu}{2}|\nabla n^\lambda|^2+\lambda^2\omega(\rho^\lambda)\right)+\mu\int_{T^N}|\nabla u^\lambda|^2
    \\&&
    +(\kappa+\mu)\int_{T^N}|\nabla\cdot u^\lambda|^2+\nu\theta\int_{T^N}\big|\Delta n^\lambda+|\nabla n^\lambda|^2n^\lambda\big|^2=0.
\eeq

Let
\beq
    \label{4.5}
    \Pi^\lambda(x,t)=\lambda^2\left(\omega(\rho^\lambda)-P(1)(\rho^\lambda-1)\right),
\eeq
and then we have $\displaystyle\int_{T^N}\Pi^\lambda(x,0)\leq C\lambda^{-2}$
from the Taylor series and (\ref{l48}).

Integrating (\ref{4.4}) over $[0,t]$ and using the law of conservation of mass, we get the basic energy estimate
\beq
    \label{4.6}
    && \nonumber
    \int_{T^N}\left(\frac12\rho^\lambda|u^\lambda|^2+\frac{\nu}{2}|\nabla n^\lambda|^2+\Pi^\lambda(x,t)\right)+\mu\int_0^t\int_{T^N}|\nabla u^\lambda|^2+(\kappa+\mu)\int_0^t\int_{T^N}|\nabla\cdot u^\lambda|^2
    \\&&+
    \nu\theta\int_0^t\int_{T^N}|\Delta n^\lambda+|\nabla n^\lambda|^2n^\lambda|^2=\int_{T^N} \left(\frac12\rho_0^\lambda|u_0^\lambda|^2+\frac{\nu}{2}|\nabla n_0^\lambda|^2+\Pi^\lambda(x,0)\right).
\eeq

Similarly, from (\ref{l40}), we get the basic energy law (see \cite{llw} for example) for the incompressible system as
\beq
    \label{4.7}
    && \nonumber \int_{T^N} \left(\frac12|u|^2+\frac{\nu}{2}|\nabla n|^2\right)+\mu\int_0^t\int_{T^N}|\nabla u|^2+\nu\theta\int_0^t\int_{T^N}|\Delta n+|\nabla n|^2n|^2
    \\&=&
    \int_{T^N}\left(\frac12|u_0|^2+\frac{\nu}{2}|\nabla n_0|^2\right).
\eeq

Secondly, multiplying (\ref{4.1})$_2$ by $u$ and then integrating over $Q_t$, we have
\beq
    \label{4.8}
    && \nonumber
    \int_{T^N}\rho^\lambda u^\lambda\cdot u+\int_0^t\int_{T^N}\rho^\lambda u^\lambda\cdot\left((u\cdot\nabla)u+\nabla p-\mu\Delta u+\nu\nabla\cdot(\nabla n\odot\nabla n)\right)
    \\&&\nonumber
    -\int_0^t\int_{T^N}(\rho^\lambda u^\lambda\otimes u^\lambda)\cdot\nabla u+\mu\int_0^t\int_{T^N} \nabla u^\lambda\cdot\nabla u+\nu\int_0^t\int_{T^N}(u\cdot\nabla)n^\lambda\cdot\Delta n^\lambda
    \\&=&
    \int_{T^N}\rho_0^\lambda u_0^\lambda\cdot u_0.
\eeq

On the other hand, multiplying (\ref{4.1})$_3$ by $\nu n$ and then integrating over $Q_t$, we obtain
\beq
    \label{4.9}
    &&\nonumber
    \nu\int_{T^N} n^\lambda\cdot n-\nu\int_{T^N} n_0^\lambda\cdot n_0+\nu\int_0^t\int_{T^N}n^\lambda\cdot\left((u\cdot\nabla)n-\theta(\Delta n+|\nabla n|^2n)\right)
    \\&&+
    \nu\int_0^t\int_{T^N} (u^\lambda\cdot\nabla)n^\lambda\cdot n=\nu\theta\int_0^t\int_{T^N}(\Delta n^\lambda+|\nabla n^\lambda|^2n^\lambda)\cdot n.
\eeq

Similarly, multiplying (\ref{4.1})$_3$ by $\nu\Delta n$ and then integrating over $Q_t$, we have
\beq
    \label{4.10}
    && \nonumber
    -\nu\int_{T^N}\nabla n^\lambda\cdot\nabla n+\nu\int_0^t\int_{T^N}\nabla n^\lambda\cdot\left(-\nabla\left((u\cdot\nabla)n\right)+\theta\left(\nabla\Delta n+\nabla\left(|\nabla n|^2n\right)\right)\right)
    \\&& \nonumber
    +\nu\int_0^t\int_{T^N}(u^\lambda\cdot\nabla)n^\lambda\cdot\Delta n-\nu\theta\int_0^t\int_{T^N}(\Delta n^\lambda+|\nabla n^\lambda|^2n^\lambda)\cdot\Delta n
    \\&=&
    -\nu\int_{T^N} \nabla n_0^\lambda\cdot\nabla n_0.
\eeq

In conclusion, using (\ref{4.6})--(\ref{4.10}), and noting the fact that
$|n^\lambda|=|n|=1$, we have
\beq
    \label{4.11}
    &&\nonumber
    \int_{T^N}\left(\frac12|\sqrt{\rho^\lambda}u^\lambda-u|^2+\frac{\nu}{2}|n^\lambda-n|^2+\frac{\nu}{2}|\nabla n^\lambda-\nabla n|^2+\Pi^\lambda(x,t)\right)
    \\&&\nonumber
    +\mu\int_0^t\int_{T^N}|\nabla u^\lambda-\nabla u|^2+(\kappa+\mu)\int_0^t\int_{T^N}|\nabla\cdot u^\lambda|^2
    \\&& \nonumber
    +\nu\theta\int_0^t\int_{T^N}\left|(\Delta n^\lambda+|\nabla n^\lambda|^2n^\lambda)-(\Delta n+|\nabla n|^2n)\right|^2
    \\&=& \nonumber
    \int_{T^N}\left(\frac12\left|\sqrt{\rho_0^\lambda}u_0^\lambda-u_0\right|^2+\frac{\nu}{2}|n_0^\lambda-n_0|^2+\frac{\nu}{2}|\nabla n_0^\lambda-\nabla n_0|^2+\Pi^\lambda(x,0)\right)+\sum\limits_{1\leq i\leq 8}R_1^\lambda(t),
    \\
\eeq
where
\beq
    \nonumber&&
    R_1^\lambda(t)=\int_{T^N}\sqrt{\rho^\lambda}(\sqrt{\rho^\lambda}-1)u^\lambda\cdot u-\int_{T^N}\sqrt{\rho_0^\lambda}\left(\sqrt{\rho_0^\lambda}-1\right)u_0^\lambda\cdot u_0,
    \\&&\nonumber
    R_2^\lambda(t)=\int_0^t\int_{T^N} \rho^\lambda u^\lambda\cdot\left((u\cdot\nabla)u\right)-\int_0^t\int_{T^N}(\rho^\lambda u^\lambda\otimes u^\lambda)\cdot\nabla u,
    \\&&\nonumber
    R_3^\lambda(t)=\int_0^t\int_{T^N}\rho^\lambda u^\lambda\cdot\nabla p,
    \\&&\nonumber
    R_4^\lambda(t)=-\mu\int_0^t\int_{T^N}\nabla u^\lambda\cdot\nabla u -\mu\int_0^t\int_{T^N}\rho^\lambda u^\lambda\cdot\Delta u,
    \\&&\nonumber
    R_5^\lambda(t)=\nu\int_0^t\int_{T^N} \rho^\lambda u^\lambda\nabla\cdot(\nabla n\odot \nabla n)+\nu\int_0^t\int_{T^N} (u\cdot\nabla)n^\lambda\cdot\Delta n^\lambda
    \\&&\nonumber \quad\quad\quad
    -\nu\int_0^t\int_{T^N}(u^\lambda\cdot\nabla)n^\lambda\cdot\Delta n-\nu\int_0^t\int_{T^N}(u\cdot\nabla)n\cdot\Delta n^\lambda,
    \\&&\nonumber
    R_6^\lambda(t)=\nu\int_0^t\int_{T^N}(u^\lambda\cdot\nabla)n^\lambda\cdot n+\nu\int_0^t\int_{T^N}(u\cdot\nabla)n\cdot n^\lambda,
    \\&&\nonumber
    R_7^\lambda(t)=-\nu\theta\int_0^t\int_{T^N}(\Delta n+|\nabla n|^2n)\cdot n^\lambda-\nu\theta\int_0^t\int_{T^N}(\Delta n^\lambda+|\nabla n^\lambda|^2n^\lambda)\cdot n,
    \\&&\nonumber
    R_8^\lambda(t)=-\nu\theta\int_0^t\int_{T^N}(\Delta n+|\nabla n|^2n)\cdot|\nabla n^\lambda|^2n^\lambda-\nu\theta\int_0^t\int_{T^N}(\Delta n^\lambda+|\nabla n^\lambda|^2n^\lambda)\cdot|\nabla n|^2n.
\eeq

Now we estimate every term above.

Firstly, it follows from the H\"older inequality, Lemma \ref{le:2.1}, (\ref{l48}) and (\ref{l49}) that
\beq
    \nonumber |R_1^\lambda(t)|
    &\leq& \nonumber
    \|u\|_\infty\left(\int_{T^N}\rho^\lambda|u^\lambda|^2\right)^{\frac12}\left(\int_{T^N}
    \left|\sqrt{\rho^\lambda}-1\right|^2\right)^{\frac12}
    \\&&\nonumber
    +\|u_0\|_\infty\left(\int_{T^N}\rho^\lambda|u_0^\lambda|^2\right)^{\frac12}\left(\int_{T^N}
    \left|\sqrt{\rho_0^\lambda}-1\right|^2\right)^{\frac12}
%    \\&\leq& \nonumber
%    \|u\|_2\left(\int_{T^N}\rho^\lambda|u^\lambda|^2\right)^{\frac12}\left(\int_{T^N}
%    \left|\frac{\rho^\lambda-1}{\sqrt{\rho^\lambda}+1}\right|^2\right)^{\frac12}
%    +\|u_0\|_2\left(\int_{T^N}\rho^\lambda|u_0^\lambda|^2\right)^{\frac12}\left(\int_{T^N}
%    \left|\frac{\rho^\lambda-1}{\sqrt{\rho_0^\lambda}+1}\right|^2\right)^{\frac12}
    \\&\leq& \nonumber
    C\left\|\frac{1}{\sqrt{\rho^\lambda}}+1\right\|_\infty\left(\int_{T^N}|\rho^\lambda-1|^2\right)^{\frac12}
    +C\left\|\frac{1}{\sqrt{\rho_0^\lambda}}+1\right\|_\infty\left(\int_{T^N}|\rho^\lambda_0-1|^2\right)^{\frac12}
    \\&\leq&
    C\lambda^{-1}.
    \\ \nonumber
    \label{R_2}
    |R_2^\lambda(t)|&=&
    \nonumber
    \bigg|-\int_0^t\int_{T^N}\left((\sqrt{\rho^\lambda}u^\lambda-u)\otimes
    (\sqrt{\rho^\lambda}u^\lambda-u)\right)\cdot\nabla u
    \\&& \nonumber
    +\int_0^t\int_{T^N}(\rho^\lambda-\sqrt{\rho^\lambda})u^\lambda\cdot\left((u\cdot\nabla)u\right) -\int_0^t\int_{T^N}(\sqrt{\rho^\lambda}u^\lambda-u)\cdot\nabla\left(\frac{|u|^2}{2}\right)\bigg|
    \\&\leq& \nonumber
    C\lambda^{-1}+C\int_0^t\int_{T^N}|\sqrt{\rho^\lambda}u^\lambda-u|^2
    +\left|\int_0^t\int_{T^N}(\sqrt{\rho^\lambda}u^\lambda-u)\cdot\nabla\left(\frac{|u|^2}{2}\right)\right|
    \\&\leq&
    C\lambda^{-1}+C\int_0^t\int_{T^N}|\sqrt{\rho^\lambda}u^\lambda-u|^2,
\eeq
where we have used in (\ref{R_2}) the following fact from the incompressibility $\nabla\cdot u=0$ that
\beq
    \nonumber
    &&\left|\int_0^t\int_{T^N}(\sqrt{\rho^\lambda}u^\lambda-u)\cdot\nabla\left(\frac{|u|^2}{2}\right)\right|
    \\&=& \nonumber
    \left|-\int_0^t\int_{T^N}\sqrt{\rho^\lambda}(\sqrt{\rho^\lambda}-1)u^\lambda\cdot\nabla\left(\frac{|u|^2}{2}\right)
    +\int_0^t\int_{T^N}\rho^\lambda
    u^\lambda\cdot\nabla\left(\frac{|u|^2}{2}\right)\right|
    \\&=& \nonumber
    \left|-\int_0^t\int_{T^N}\sqrt{\rho^\lambda}(\sqrt{\rho^\lambda}-1)u^\lambda\cdot\nabla\left(\frac{|u|^2}{2}\right)
    +\int_0^t\int_{T^N}\rho^\lambda_t\cdot\frac{|u|^2}{2}\right|
    \\&\leq&\nonumber
    C\left(\int_0^t\int_{T^N}|\rho^\lambda-1|^2\right)^{\frac12}
    +C\left(\int_0^t\int_{T^N}|\rho^\lambda_t|^2\right)^{\frac12}
    \\&\leq&
    C\lambda^{-1}.
\eeq

Then we estimate the third and forth terms in a similar way as
\beq
    |R_3^\lambda(t)|&=&\left|\int_0^t\int_{T^N}\rho^\lambda_t p\right|
    \leq
    C\left(|\int_0^t\int_{T^N}|\rho^\lambda_t|^2\right)^{\frac12}\left(\int_0^t\int_{T^N}
    p^2\right)^{\frac12}\leq C\lambda^{-1},
    \\|R_4^\lambda(t)| &\leq&
    \nonumber
    \mu\int_0^t\int_{T^N}|(1-\rho^\lambda)u^\lambda\cdot\Delta u|
    \\&\leq&
    C\left(\int_0^t\int_{T^N}|1-\rho^\lambda|^2\right)^{\frac12}\left(\int_0^t\int_{T^N}|u^\lambda|^2|\Delta u|^2\right)^{\frac12} \leq C\lambda^{-1}.
\eeq

Now we turn to estimate $R_5^\lambda(t)$.
\beq
    \nonumber R_5^\lambda(t)&=&
    \nonumber
    \nu\int_0^t\int_{T^N} \rho^\lambda(u^\lambda\cdot\nabla)n\cdot\Delta n+\nu\int_0^t\int_{T^N}
    \rho^\lambda u^\lambda\cdot\nabla\left(\frac{|\nabla
    n|^2}{2}\right)
    \\&& \nonumber
    +\nu\int_0^t\int_{T^N}
    (u\cdot\nabla)n^\lambda\cdot\Delta n^\lambda-\nu\int_0^t\int_{T^N}(u^\lambda\cdot\nabla)n^\lambda\cdot\Delta
    n-\nu\int_0^t\int_{T^N}(u\cdot\nabla)n\cdot\Delta n^\lambda
%    \\&=& \nonumber
%    \nu\int_0^t\int_{T^N} (\rho^\lambda-1)(u^\lambda\cdot\nabla)n\cdot\Delta n+\nu\int_0^t\int_{T^N}(u^\lambda\cdot\nabla)n\cdot\Delta n+\nu\int_0^t\int_{T^N}(\rho^\lambda -1)_t\frac{|\nabla n|^2}{2}
%    \\&& \nonumber
%    +\nu\int_0^t\int_{T^N} (u\cdot\nabla)n^\lambda\cdot\Delta n^\lambda-\nu\int_0^t\int_{T^N}(u^\lambda\cdot\nabla)n^\lambda\cdot\Delta n-\nu\int_0^t\int_{T^N}(u\cdot\nabla)n\cdot\Delta n^\lambda
    \\&=& \nonumber
    \nu\int_0^t\int_{T^N} (\rho^\lambda-1)(u^\lambda\cdot\nabla)n\cdot\Delta n+\nu\int_0^t\int_{T^N}\rho^\lambda_t\frac{|\nabla n|^2}{2}
    \\&& \nonumber
    -\nu\int_0^t\int_{T^N}((u^\lambda-u)\cdot\nabla)(n^\lambda-n)\cdot\Delta n+\nu\int_0^t\int_{T^N}(u\cdot\nabla)(n^\lambda-n)(\Delta n^\lambda-\Delta n)
    \\&=& \nonumber
    \nu\int_0^t\int_{T^N} (\rho^\lambda-1)(u^\lambda\cdot\nabla)n\cdot\Delta n+\nu\int_0^t\int_{T^N}\rho^\lambda_t\frac{|\nabla n|^2}{2}
    \\&& \nonumber
    -\nu\int_0^t\int_{T^N}\big((\sqrt{\rho^\lambda}u^\lambda-u)\cdot\nabla\big)(n^\lambda-n)\cdot\Delta n
    \\&& \nonumber
    -\int_0^t\int_{T^N}\big((1-\sqrt{\rho^\lambda})u^\lambda\cdot\nabla\big)(n^\lambda-n)\cdot\Delta n
    \\&&
    -\nu\int_0^t\int_{T^N}(\nabla u\cdot\nabla)(n^\lambda-n)\cdot(\nabla n^\lambda-\nabla n),
\eeq
where we have used integration by parts and the incompressibility $\nabla\cdot u=0$ in the last equality.
Then by the H\"older inequality and the uniform estimates (\ref{l49}), we get
\beq
    |R_5^\lambda(t)|&\leq&  \nonumber
    C\left(\int_0^t\int_{T^N}|\rho^\lambda-1|^2\right)^{\frac12}\left(\int_0^t\int_{T^N}|u^\lambda|^2|\nabla n|^2|\Delta n|^2\right)^{\frac12}
    \\&& \nonumber
    +C\left(\int_0^t\int_{T^N}|\rho^\lambda_t|^2\right)^{\frac12}\left(\int_0^t\int_{T^N}\frac{|\nabla n|^4}{4}\right)^{\frac12}
    \\&& \nonumber
    +C\|\Delta n\|_\infty\int_0^t\left(\int_{T^N}|\sqrt{\rho^\lambda}u^\lambda-u|^2\right)^{\frac12}\bigg(\int_{T^N}|\nabla n^\lambda-\nabla n|^2\bigg)^{\frac12}
    \\&& \nonumber
    +C\|u^\lambda\|_\infty\|\Delta n\|_\infty\int_0^t\left(\int_{T^N}|1-\sqrt{\rho^\lambda}|^2\right)^{\frac12}\left(\int_{T^N}|\nabla n^\lambda-\nabla n|^2\right)^{\frac12}
    \\&&\nonumber
    +C\|\nabla u\|_\infty\int_0^t\int_{T^N}|\nabla n^\lambda-\nabla n|^2
    \\&\leq&
    C(\lambda^{-1}+\lambda^{-2})+C\int_0^t\int_{T^N}\left(|\sqrt{\rho^\lambda}u^\lambda-u|^2+|\nabla n^\lambda-\nabla n|^2\right).\quad\quad\quad
\eeq

Note that $|n|=1$ and $|n^\lambda|=1$ in $Q_{T_0}$. Then we have
\beq
    \label{4.19}
    \nabla(|n^\lambda|^2)=\nabla(|n|^2)=0,\quad
    (\Delta n+|\nabla n|^2n)\cdot n=(\Delta n^\lambda+|\nabla n^\lambda|^2n^\lambda)\cdot n^\lambda=0.
\eeq

With the help of (\ref{4.19}), one obtains
\beq
    \nonumber
    R_6^\lambda(t)&=&\nu\int_0^t\int_{T^N}(u^\lambda\cdot\nabla)n^\lambda\cdot(n-n^\lambda)
    +\nu\int_0^t\int_{T^N}(u\cdot\nabla)n\cdot(n^\lambda-n)
    \\&=&
    \nu\int_0^t\int_{T^N}\big((u^\lambda-u)\cdot\nabla\big)n^\lambda\cdot(n-n^\lambda) +\nu\int_0^t\int_{T^N}(u\cdot\nabla)(n-n^\lambda)\cdot(n^\lambda-n).\quad\quad\quad
\eeq

Then it follows from the H\"older inequality and (\ref{l49}) that
\beq
    |R_6^\lambda(t)|&\leq& \nonumber C\|\nabla n^\lambda\|_\infty\left(\int_0^t\int_{T^N}|u^\lambda-u|^2\right)^{\frac12}\left(\int_0^t\int_{T^N}|n^\lambda-n|^2\right)^{\frac12}
    \\&&\nonumber
    +C\|u\|_\infty\left(\int_0^t\int_{T^N}|n^\lambda-n|^2\right)^{\frac12}\left(\int_0^t\int_{T^N}|\nabla n^\lambda-\nabla n|^2\right)^{\frac12}
    \\&\leq&
    C\int_0^t\int_{T^N}|u^\lambda-u|^2+C\int_0^t\int_{T^N}|n^\lambda-n|^2+C\int_0^t\int_{T^N}|\nabla n^\lambda-\nabla n|^2.\quad\quad
\eeq

Recalling (\ref{4.19}), we have
\beq
    \nonumber
    R_7^\lambda(t)&=&-\nu\theta\int_0^t\int_{T^N}(\Delta n+|\nabla n|^2n)\cdot(n^\lambda-n)-\nu\theta\int_0^t\int_{T^N}(\Delta n^\lambda+|\nabla n^\lambda|^2n^\lambda)\cdot (n-n^\lambda)
    \\&=&
    -\nu\theta\int_0^t\int_{T^N}\left((\Delta n+|\nabla n|^2n)-(\Delta n^\lambda+|\nabla n^\lambda|^2n^\lambda)\right)\cdot(n^\lambda-n),
\eeq
and then we can use the Cauchy inequality to get
\beq
    |R_7^\lambda(t)|\leq\frac{\nu\theta}{4}\int_0^t\int_{T^N}\big|(\Delta n+|\nabla n|^2n)-(\Delta n^\lambda+|\nabla n^\lambda|^2n^\lambda)\big|^2+C\int_0^t\int_{T^N}|n-n^\lambda|^2.\quad
\eeq

Similarly, with the help of (\ref{4.19}), we derive
\beq
    \nonumber
    R_8^\lambda(t)&=& \nonumber
    -\nu\theta\int_0^t\int_{T^N}(\Delta n+|\nabla n|^2n)\cdot|\nabla n^\lambda|^2(n^\lambda-n)
    \\&&\nonumber
    -\nu\theta\int_0^t\int_{T^N}(\Delta n^\lambda+|\nabla n^\lambda|^2n^\lambda)\cdot|\nabla n|^2(n-n^\lambda)
    \\&=& \nonumber
    \nu\theta\int_0^t\int_{T^N}\left((\Delta n^\lambda+|\nabla n^\lambda|^2n^\lambda)-(\Delta n+|\nabla n|^2n)\right)\cdot|\nabla n|^2(n^\lambda-n)
    \\&&
    -\nu\theta\int_0^t\int_{T^N}((\nabla n^\lambda-\nabla n)\cdot(\nabla n^\lambda+\nabla n))(\Delta n+|\nabla n|^2n)\cdot(n^\lambda-n),
\eeq
and then by the Cauchy inequality and (\ref{l49}), we have
\beq
    \nonumber |R_8^\lambda(t)|&\leq&\frac{\nu\theta}{4}\int_0^t\int_{T^N}\big|(\Delta n^\lambda+|\nabla n^\lambda|^2n^\lambda)-(\Delta n+|\nabla n|^2n)\big|^2+C\|\nabla n\|_\infty^4\int_0^t\int_{T^N}|n^\lambda-n|^2
    \\&&\nonumber
    +C\left(\|\nabla n^\lambda\|_\infty^2+\|\nabla n\|_\infty^2\right)\|\Delta n+|\nabla n|^2n\|_\infty^2\int_0^t\int_{T^N}|\nabla n^\lambda-\nabla n|^2
    \\&&   \nonumber
    +C\int_0^t\int_{T^N}|n^\lambda-n|^2
    \\&\leq& \nonumber
    \frac{\nu\theta}{4}\int_0^t\int_{T^N}|(\Delta n^\lambda+|\nabla n^\lambda|^2n^\lambda)-(\Delta n+|\nabla n|^2n)|^2
    \\&&
    +C\int_0^t\int_{T^N}\left(|n^\lambda-n|^2+|\nabla n-\nabla n^\lambda|^2\right).
\eeq

In conclusion, combining all the estimates of $R_i^\lambda(t)$ for $i=1,2,\ldots,8$, we have
\beq
    \label{4.27}
    &&\nonumber
    \int_{T^N}\left(\frac12|\sqrt{\rho^\lambda}u^\lambda-u|^2+\frac{\nu}{2}|n^\lambda-n|^2+\frac{\nu}{2}|\nabla n^\lambda-\nabla n|^2+\Pi^\lambda(x,t)\right)
    \\&&\nonumber
    +\mu\int_0^t\int_{T^N}|\nabla u^\lambda-\nabla u|^2+(\kappa+\mu)\int_0^t\int_{T^N}|\nabla\cdot u^\lambda|^2
    \\&& \nonumber
    +\frac{\nu\theta}{2}\int_0^t\int_{T^N}|(\Delta n^\lambda+|\nabla n^\lambda|^2n^\lambda)-(\Delta n+|\nabla n|^2n)|^2
    \\&\leq& \nonumber
    \int_{T^N}\left(\frac12\left|\sqrt{\rho_0^\lambda}u_0^\lambda-u_0\right|^2+\frac{\nu}{2}|\nabla n_0^\lambda-\nabla n_0|^2+\Pi^\lambda(x,0)\right)+C\lambda^{-1}
    \\&&
    +C\int_0^t\int_{T^N}\left(|\sqrt{\rho^\lambda}u^\lambda-u|^2+|n^\lambda-n|^2+|\nabla n^\lambda-\nabla n|^2\right).
\eeq

By the following facts indicated from the H\"older inequality, (\ref{l48}) and (\ref{l49}) that
\beq
    \label{4.28}
    \int_{T^N}\left|\sqrt{\rho_0^\lambda} u_0^\lambda-u_0\right|^2=\int_{T^N}\left|\sqrt{\rho_0^\lambda} u_0^\lambda-u_0^\lambda\right|^2+\int_{T^N}|u_0^\lambda-u_0|^2\leq C(\lambda^{-4}+\lambda^{-2}),
\eeq
and
\beq
    \label{4.29}
    \int_{T^N}|u^\lambda-u|^2 \leq\int_{T^N}|\sqrt{\rho^\lambda}u^\lambda-u^\lambda|^2+\int_{T^N}|\sqrt{\rho^\lambda}u^\lambda-u|^2
    \leq
    C\lambda^{-2}+\int_{T^N}|\sqrt{\rho^\lambda}u^\lambda-u|^2.\
\eeq

By (\ref{4.27})--(\ref{4.29}) and the Gronwall inequality, we have
\beq
    \label{4.30}
    \|u^\lambda-u\|^2+\nu\|\nabla n^\lambda-\nabla n\|^2\leq C\lambda^{-1} \ {\rm for} \ t\in[0,T_0].
\eeq
Then integrating (\ref{4.27}) over $[0,t]$, we get
\beq
    \label{4.31}
    \int_0^t(\mu\|\nabla u^\lambda-\nabla u\|^2+\nu\theta\|(\Delta n^\lambda+|\nabla n^\lambda|^2n^\lambda)-(\Delta n+|\nabla n|^2n)\|^2)\leq C\lambda^{-1}.
\eeq

Finally, it follows from the subtraction of (\ref{l38})$_3$ from (\ref{l40})$_3$ that
\beq
    \label{4.32}
    \nonumber
    (n^\lambda-n)_t-\theta\Delta(n^\lambda-n)&=&-\big((u^\lambda-u)\cdot\nabla\big)n^\lambda-(u\cdot\nabla)(n^\lambda-n)
    \\&&
    +\theta((\nabla n^\lambda-\nabla n):(\nabla n^\lambda+\nabla n))\cdot n^\lambda+\theta|\nabla n|^2(n^\lambda-n).\quad\quad
\eeq
Then the parabolic theory implies that
\beq
    \label{4.33}
    \nonumber
    \|n^\lambda-n\|_2^2&\leq&C\|\nabla n^\lambda\|_\infty^2\|u^\lambda-u\|^2+C\|u\|_\infty^2\|\nabla n^\lambda-\nabla n\|^2
    \\&&
    +C\left(\|\nabla n^\lambda\|_\infty^2+\|\nabla n\|_\infty^2\right)\|\nabla n^\lambda-\nabla n\|^2+C\|\nabla n\|_\infty^4\|n^\lambda-n\|^2
    \leq
    C\lambda^{-1} \quad\quad\
\eeq for $t\in[0,T_0]$.

Furthermore, applying $\nabla$ to (\ref{4.32}), we have
\beq
    \label{4.34}
    && \nonumber
    \nabla(n^\lambda-n)_t-\theta\Delta\nabla(n^\lambda-n)
    \\&=& \nonumber
    -(\nabla(u^\lambda-u)\cdot\nabla)n^\lambda-((u^\lambda-u)\cdot\nabla)\nabla n^\lambda-(\nabla u\cdot\nabla)(n^\lambda-n)-(u\cdot\nabla)\nabla(n^\lambda-n)
    \\&& \nonumber
    +\theta((\nabla^2n^\lambda-\nabla^2 n):(\nabla n^\lambda+\nabla n))\cdot n^\lambda+\theta((\nabla n^\lambda-\nabla n):(\nabla^2n^\lambda+\nabla^2 n))\cdot n^\lambda
    \\&&\nonumber
    +\theta((\nabla n^\lambda-\nabla n):(\nabla n^\lambda+\nabla n))\cdot\nabla n^\lambda +2\theta(\nabla n:\nabla^2n)\cdot(n^\lambda-n)+\theta|\nabla n|^2(\nabla n^\lambda-\nabla n),\\
\eeq
and then (\ref{4.30}), together with (\ref{4.31}) and (\ref{4.33}), implies that
\beq
    \int_0^t\|\nabla n^\lambda-\nabla n\|_2^2\leq C\lambda^{-1}.
\eeq

This completes the proof of Theorem \ref{th.1.3}.
\end{Proof of Theorem $2.3$}

\section*{Acknowledgment}
The authors would like to
thank Professor Chun Liu and Professor Zhen Lei for their suggestions and sincere help. This work is supported by the National
Basic Research Program of China (973 Program) (No. 2011CB808002), the National Natural
Science Foundations of China (No. 11071086 and No. 11128102), the University Special Research Foundation for Ph.D Program (No. 20104407110002), the Guangdong Province National Science Foundation (No. S2012010010408) and the China Postdoctoral Science Foundation funded project (No. 2012M521443).


\vskip 0.8cm


\section*{References}
\addcontentsline{toc}{section}{Bibliography}
\begin{thebibliography}{99}

%\bibitem{Caffarelli} Caffarelli, L., Kohn, R., Nirenberg, L.: Partial regularity of suitable weak solutions
%of the Navier-Stokes equations. {\em Comm. Pure Appl. Math.} {\bf 35}, 771--831 (1982).
%
%\bibitem{Chang} Chang, K.C., Ding, W.Y., Ye, R.: Finite-time blow-up of the heat flow of harmonic maps from surfaces. {\em J. Differential
%Geom.} {\bf 36}(2) 507--515 (1992).
%
%\bibitem{Cho} Cho, Y.G., Choe, H.J., Kim, H.: Unique solvability of the initial boundary value problems for compressible viscous fluids. {\em J. Math. Pures Appl.} {\bf 83}, 243--275 (2003).
%
%\bibitem{Choe} Choe H.J., Kim H.: Global existence of the radially symmetric solutions of the
%Navier-Stokes equations for the isentropic compressible fuids.
%{\em Math. Meth. Appl. Sci.} {\bf 28}, 1--28 (2005).
%
%\bibitem{Ding} Ding, S.J., Lin, J.Y., Wang, C.Y., Wen, H.Y.: Compressible hydrodynamic flow of liquid crystals in
%1-D. {\em Discrete Contin. Dyn. Syst.} {\bf 32}(2) 539--563 (2012).
%
%\bibitem{Ding2} Ding, S.J., Huang, J.R., Lin, J.Y.: Global existence of solutions to slightly compressible model of liquid crystals in two dimension. in prepare.
%
%\bibitem{Erickson1} Ericksen, J.: Conservation laws for liquid crystals, Trans. Soc. Rheol. {\bf 5}, 22--34 (1961).
%
%\bibitem{Erickson2} Ericksen, J.: Hydrostatic Theory of Liquid Crystal. {\em Arch. Rational Mech. Anal.}
%{\bf 9}, 371--378 (1962).
%
%\bibitem{Fan} Fan, J.S., Jiang, S., Ni, G.X.: Uniform boundedness of the radially symmetric solutions
%of the Navier-Stokes equations for isentropic compressible fluids.
%{\em Osaka J. Math.} {\bf 46}, 863--876 (2009).
%
%\bibitem{Genn} P.G. de Gennes: The Physics of Liquid Crystals. Oxford, 1974.
%
%\bibitem{Hagstrom} Hagstrom, T., Lorenz, J.: All-time existence of classical solutions for slightly compressible flows. {\em SIAM J. MATH. ANAL.} {\bf 29}, 652--672 (1998).
%
%\bibitem{Hagstrom2} Hagstrom, T., Lorenz, J.: All-time existence of smooth solutions to PDEs of mixed type and the invariant subspace of uniform states, {\em Adv. Appl. Math.} {\bf 16}, 219--257 (1995).
%
%\bibitem{Hoff} Hoff, D.: Spherically symmetric solutions of the Navier-Stokes
%equations for compressible, isothermal flow with large,
%discontinuous initial data. {\em Indiana Univ. Math. J.} {\bf 41}, 1225--1302 (1992).
%
%\bibitem{Huf}Hu, X.P., Wang, D.H.: Low Mach number limit of viscous compressible
%magnetohydrodynamic flows. {\em SIAM J. Math. Anal.} {\bf 41}, 1272--1294 (2009).
%.
%
%\bibitem{Huang-Li-Xin} Huang, X.D., Li, J., Xin, Z.P.: Global Well-Posedness of Classical solutions with Large oscillations and vacuum to the three-dimensional isentropic compressible Navier-Stokes equations, {\em Comm. Pure. Appl. Math.} to appear.
%
%\bibitem{JiangTan} Jiang, F., Tan, Z.: Global weak solution to the flow of liquid crystals system. {\em Math. Methods Appl. Sci.}
%{\bf 32}, 2243--2266 (2009).
%
%\bibitem{Jiang 1} Jiang, S., Ju, Q.C., Li, F.C.: Incompressible limit of the
%compressible magnetohydrodynamic equations with periodic boundary
%conditions. {\em Commun. Math. Phys.} {\bf 297}, 371--400 (2010).
%
%\bibitem{Jiang 2} Jiang, S., Zhang, P.: On spherically symmetric solutions of the
%compressible isentropic Navier- Stokes equations. {\em Comm. Math.
%Phys.} {\bf 215}, 559--581 (2001).
%
%\bibitem{Klainerman} Klainerman, S., Majda, A.: Sinaular limits
%of quasilinear hyperbolic system with large paramiters and the
%incompressible limit of compressible fluids. {\em Comm. Pure Appl. Math.}
%{\bf 34}, 481--524 (1981).
%
%\bibitem{Ladyzenskaja} Ladyzenskaja, O.A., Solonnikov, V.A., Ural'ceva, N.N.:{\em Linear and  quasilinear equations of parabolic type.} Amer. Math.
%Soc., Providence RI, 1968.
%
%\bibitem{Lei} Lei, Z.: Global existence of classical solution for some Oldroyd-B model via the incompressible limit.
%{\em Chin. Ann. Math.} {\bf 27B}(5) 565--580 (2006).
%
%\bibitem{Lei2} Lei, Z., Zhou, Y.: Global existence of classical solution for the two-dimensional Oldroyd model via the incompressible limit. {\em SIAM J.MATH.ANAL.} {\bf 37}(3) 797--814 (2005).
%
%\bibitem{leslie} Leslie, F.: Some Constitute Equations for Anisotropic Fluids.
%{\em Q. J. Mech. Appl. Math.} {\bf 19}, 357--370 (1966).
%
%\bibitem{1} Lin, F.H.: Nonlinear theory of defects in nematic liquid
%crystal: phase transition and flow phenomena. {\em Comm. Pure Appl.
%Math.} {\bf 42}, 789--814 (1989).
%
%\bibitem{2} Lin, F.H., Liu, C.: Nonparabolic dissipative systems
%modeling the flow of liquid crystals. {\em Comm. Pure Appl. Math.} Vol.
%{\bf XLVIII}, 501--537 (1995).
%
%\bibitem{3} Lin, F.H., Liu, C.: Partial regularities of the nonlinear
%dissipative systems modeling the flow of liquid crystals. {\em Discrete Contin. Dyn. Syst.}
%{\bf 2}, 1--23 (1996).
%
%\bibitem{4} Lin, F.H., Liu, C.: Existence of solutions for the
%Ericksen-Leslie system. {\em Arch. Rational Mech. Anal.} {\bf 154},
%135--156 (2000).
%
%\bibitem{5} Lin, F.H., Liu, C., Zhang, P.: On hydrodynamics of viscoelastic fluids. {\em Comm. Pure Appl. Math.} Vol. {\bf LVIII}, 1437--1471 (2005).
%
%\bibitem{llw} Lin, F.H., Lin, J.Y., Wang, C.Y.: Liquid crystal flows in dimensions two. {\em Arch. Rational Mech. Anal.} {\bf 197}, 297--336 (2010).
%
%\bibitem{zhang} Lin, P., Liu, C., Zhang, H.: An energy law preserving $C^0$ finite element scheme for
%simulating the kinematic effects in liquid crystal dynamics. {\em J. Comp. Phys.} {\bf 227}, 1411--1427 (2007).
%
%\bibitem{Liu C} Liu, C., Walkington, N.J.: Mixed methods for the approximation of
%liquid crystal flows. {\em Math. Modeling and Numer Anal.} {\bf 36}, 205--222 (2002).
%
%\bibitem{Liushen} Liu C., Shen J.: On liquid crystal flows with free-slip boundary conditions. {\em Discrete Contin. Dyn. Syst.}
%{\bf 7}, 307--318 (2001).
%
%\bibitem{Daniel} Coutand, D., Shkoller, S.: Well-posedness of the full Ericksen-Leslie model of nematic liquid crystals.
%{\em C. R. Acad. Sci. Paris.} S$\acute{e}$r. {\bf I333}, 919--924 (2001).
%
%\bibitem{Liu} Liu, X.G., Zhang, Z.Y.: Existence of the flow of liquid crystals
%system. {\em Chin. Ann. Math.} {\bf 30A}(1) 1--20 (2009).
%
%\bibitem{Sideris} Sideris, T.C., Thomases, B.: Global existence for 3D incompressible isotropic elastodynamics via the
%incompressible limit. {\em Comm. Pure Appl. Math.} {\bf 57}, 1--39 (2004).
%
%\bibitem{Simon} Simon, J.: Nonhomogeneous viscous incompressible fluids: existence of
%vecocity, density and pressure, {\em SIAM J. Math. Anal}. {\bf 21}(5) 1093--1117 (1990).
%
%
%\bibitem{Valli} Valli, A.: Periodic and stationary solutions for compressible Navier-Stokes equations via a stability
%method, {\em Ann. Scuola Norm. Sup. Pisa.} {\bf 10}, 607--647 (1983).
%
%\bibitem{Wang dh} Wang, D.H., Yu, C.: Incompressible limit for the compressible flow of liquid crystals, arXiv:1108.4941v1.
%
%\bibitem{Weigant} Weigant, V.A.: Example of non-existence in the large for the problem
%of the existence of solutions of Navier-Stokes equations for
%compressible viscous barotropic fluids. {\em Dokl. Akad. Nauk.} {\bf 339},
%155--156 (1994).






\bibitem{Blanca}
         B. Climent-Ezquerra, F. Guill${\rm \acute{e}}$n-Gonz${\rm \acute{a}}$lez, M. Rojas-Medar,
         Reproductivity for a nematic liquid crystal model,
         Z. angew. Math. Phys. 57 (2006) 984--998.

%\bibitem{chen}
%         Y.Z. Chen,
%         Parabolic equations of second order,
%         Beijing University Press, China, 2003

\bibitem{cnk}
         L. Caffarelli, R. Kohn, L. Nirenberg,
         Partial regularity of suitable weak solutions of the Navier-Stokes equations,
         Comm. Pure Appl. Math. 35 (1982) 771--831.

\bibitem{w1}
         S.J. Ding, J.Y. Lin, C.Y. Wang, H.Y. Wen,
         Compressible hydrodynamic flow of liquid crystals in 1D,
         Discret. Contin. Dynam. Systems 32 (2012) 539--563.

%\bibitem{Ding-Huang-Wen-Zi}
%         S. J. Ding, J. R. Huang, H. Y. Wen, R. H. Zi,
%         Incompressible Limit of the Compressible Hydrodynamic Flow of Liquid Crystals,
%         Preprint, 2011.

\bibitem{Ding2} S.J. Ding, J.R. Huang, J.Y. Lin, Global existence of solutions to slightly compressible hydrodynamic flow of liquid crystals in two dimension, Sci. China Math., in press.

\bibitem{w2}
         S.J. Ding, C.Y. Wang, H.Y. Wen,
         Weak solution to compressible hydrodynamic flow of liquid crystals in dimension one,
         Discret. Contin. Dynam. Systems 15 (2011) 357--371.

\bibitem{Erickson1} J. Ericksen, Conservation laws for liquid crystals, Trans. Soc. Rheol. 5 (1961) 22--34.

\bibitem{Erickson2}
         J. Ericksen,
         Hydrostatic theory of liquid crystal,
         Arch. Rational Mech. Anal. 9 (1962) 371--378.

%\bibitem{fang}
%         D.Y. Fang, T. Zhang,
%         Compressible Navier¨CStokes equations with vacuum state in the case of general pressure law,
%         Math.Methods Appl. Sci. 29 (2006) 1081--1106.


\bibitem{feireisl}
         E. Feireisl, E. Rocca, G. Schimperna,
         On a non-isothermal model for nematic liquid crystals,
         Nonlinearity 24(1) (2011) 243--257.

\bibitem{Genn} P.G. de Gennes, The Physics of Liquid Crystals, Oxford, 1974.

%\bibitem{guo1}
%         Z.H. Guo, S. Jiang, F. Xie,
%         Global weak solutions and asymptotic behavior to 1d compressible Navier-Stokes equations with degenerate
%         viscosity coeffcient,
%         Asymptot. Anal. 60 (2008) 101--123.
%
%\bibitem{guo3}
%         Z.H. Guo, C.J. Zhu,
%         Remarks on one-dimensional compressible Navier-Stokes equations with density-dependent viscosity and vacuum,
%         J. Diff. Eqs. 248 (2010) 2768--2799.

%\bibitem{Hagstrom} T. Hagstrom, J. Lorenz, All-time existence of classical solutions for slightly compressible flows, SIAM J. Math. Anal. 29 (1998) 652--672.

\bibitem{Hagstrom2} T. Hagstrom, J. Lorenz, All-time existence of smooth solutions to PDEs of mixed type and the invariant subspace of uniform states, Adv. Appl. Math. 16 (1995) 219--257.

\bibitem{Huang-Wang-Wen1}
         T. Huang, C.Y. Wang, H.Y. Wen,
         Strong solutions of the compressible nematic liquid crystal flow,
         J. Differential Equations 252 (2012) 2222--2265.

\bibitem{Huang-Wang-Wen2}
         T. Huang, C.Y. Wang, H.Y. Wen,
         Blow up criterion for compressible nematic liquid crystal flows in dimension three,
         Arch. Rational Mech. Anal. 204 (2012) 285--311.

\bibitem{hkl} Y. Hyon, D.Y. Kwak, C. liu, Energetic variational approach in complex fluids: maximum dissipation principle, Discret. Contin. Dynam. Systems 26 (2010) 1291--1304.


%\bibitem{j1}
%         S. Jiang,
%         Global smooth solutions of the equations of a viscous, heat-conducting one-dimensional gas with density dependent viscosity,
%         Math. Nachr. 190 (1998) 169--183.


\bibitem{jiang-tan}
         F. Jiang, Z. Tan,
         Global weak solution to the flow of liquid crystals system,
         Math. Methods Appl. Sci. 30 (2009) 2243--2266.

\bibitem{Jiang 1} S. Jiang, Q.C. Ju, F.C. Li, Incompressible limit of the compressible magnetohydrodynamic equations with periodic boundary conditions, Commun. Math. Phys. 297 (2010) 371--400.

%\bibitem{j2}
%         S. Jiang, Z.P. Xin and P. Zhang,
%         Global weak solutions to 1D compressible isentropy Navier¨CStokes with densitydependent viscosity,
%         Meth. Appl. Anal. 12 (2005) 239--252.

%\bibitem{k}
%         A.V. Kazhikhov,
%         Sur la solubilit$\acute{e}$ globale des probl$\acute{e}$me monodimensionnelle aux valeurs initialeslimit$\acute{e}$es pour les $\acute{e}$quations du gaz visqueux et calorif$\acute{e}$re,
%         C. R. Acad. Sci. Paris. 284 (1977) Ser. A, 317-320.

\bibitem{Klainerman}
         S. Klainerman, A. Majda,
         Singular limits of quasilinear hyperbolic system with large paramiters and the incompressible limit of compressible fluids, Comm. Pure Appl. Math. 34 (1981) 481--524.

%\bibitem{15}
%         O.A. Lady\u{z}enskaja, V.A. Solonnikov, N.N. Ural'ceva,
%         Linear and Quasilinear Equations of Parabolic Type,
%         Translations of Mathematical Monographs, 23, Amer. Math. Soc. Providence RI, 1967.

\bibitem{Lei} Z. Lei, Global existence of classical solution for some Oldroyd-B model via the incompressible limit,
Chin. Ann. Math. Ser. B 27(5) (2006) 565--580.

\bibitem{Lei2} Z. Lei, Y. Zhou, Global existence of classical solution for the two-dimensional Oldroyd model via the incompressible limit, SIAM J. Math. Anal. 37(3) (2005) 797--814.

\bibitem{leslie}
         F. Leslie, Some constitutive equations for liquid crystals, Arch. Rational Mech. Anal. 28 (1968) 265--283.

%\bibitem{Li-Li-Xin}
%        H. L. Li, J. Li, Z. P. Xin,
%        Vanishing of Vacuum States and Blow-up Phenomena of the Compressible Navier-Stokes Equations,
%        Commun. Math. Phys. 281(2008) 401--444.


\bibitem{Lijing}
        J. Li, Z.H. Xu, J.W. Zhang, Global well-posedness with large oscillations and vacuum to the three-dimensional equations of compressible nematic liquid crystal flows, preprint arxiv:1204.4966v1 (2012) 1--46.

\bibitem{Li-Wang}
         X.L. Li, D.H. Wang,
         Global strong solution to the density-dependent incompressible flow of liquid crystals,
         preprint arXiv:1202.1011v1 (2012) 1--36.

%\bibitem{Lian-Guo-Li}
%        R. X. Lian, Z. H. Guo, H. L. Li,
%        Dynamical behaviors for 1D compressible Navier-Stokes equations with density-dependent viscosity,
%        J. Differential Equations, 248 (2010) 1926--1954.

\bibitem{1}
         F.H. Lin,
         Nonlinear theory of defects in nematic liquid crystals: phase transition and flow phenomena,
         Comm. Pure Appl. Math. 42 (1989) 789--814.

\bibitem{2}
         F.H. Lin, C. Liu,
         Nonparabolic dissipative systems modeling the flow of liquid crystals,
         Comm. Pure Appl. Math. 48 (1995) 501--537.

\bibitem{4}
         F.H. Lin, C. Liu,
         Partial regularities of the nonlinear dissipative systems modeling the flow of liquid crystals,
         Discret. Contin. Dynam. Systems 2 (1996) 1--22.

\bibitem{3}
         F.H. Lin, C. Liu,
         Existence of solutions for the Ericksen-Leslie system,
         Arch. Rational Mech. Anal. 154 (2000) 135--156.

\bibitem{llw}
         F.H. Lin, J.Y. Lin, C.Y. Wang,
         Liquid crystal flows in two dimensions,
         Arch. Rational Mech. Anal. 197 (2010) 297--336.

\bibitem{zhang} P. Lin, C. Liu, H. Zhang, An energy law preserving $C^0$ finite element scheme for
simulating the kinematic effects in liquid crystal dynamics, J. Comp. Phys. 227 (2007) 1411--1427.

%\bibitem{llq}
%         C. Liu, X. Liu, J. Qing,
%         Existence of globally weak solutions to the Flow of Compressible Liquid Crystals System, 2011, preprint.????

\bibitem{liu-walkington}
         C. Liu, N.J. Walkington,
         Mixed methods for the approximation of liquid crystal flows,
         Math. Model. Numer. Anal. 36 (2002) 205--222.

\bibitem{liu-liu-hao1}
         X.G. Liu, L.M. Liu, Y.H. Hao,
         Existence of strong solutions for the compressible Ericksen-Leslie model,
         preprint arXiv:1106.6140v1 (2011) 1--32.

\bibitem{liu-liu-hao2}
         X.G. Liu, L.M. Liu, Y.H. Hao,
         A blow-up criterion of strong solutions to the compressible liquid crystals system,
         Chin. Ann. Math. Ser. A 32(4) (2011) 393--406.

\bibitem{liu-hao} X.G. Liu, Y.H. Hao, Incompressible limit of a compressible liquid crystals system, preprint arXiv:1201.5942 (2012) 1--17.



%\bibitem{t1}
%         T.P. Liu, Z. Xin, T. Yang,
%         Vacuum states of compressible flow,
%         Discrete Contin. Dynam. Systems, 4 (1998) 1--32.

\bibitem{liu-zhang}
         X. Liu, Z. Zhang,
         Existence of the flow of liquid crystals system,
         Chin. Ann. Math. Ser. A 30(1) (2009) 1--20.

%\bibitem{9}
%         T. Luo, Z. Xin, T. Yang,
%         Interface behavior of compressible Navier-Stokes equations with vacuum,
%         SIAM J. Math. Anal. 31 (2000) 1175--1191.
%
%\bibitem{m1}
%         M. Okada,
%         Free boundary value problems for the equation of one-dimensional motion of compressible viscous fluids,
%         Japan J. Appl. Math. 4 (1987) 219--235.
%
%\bibitem{m2}
%         M. Okada,
%         Free boundary value problems for the equation of one-dimensional motion of viscous gas,
%         Japan J. Appl. Math. 6 (1989) 161--177.

%\bibitem{m3}
%         M. Okada, T. Makino.
%         Free boundary problem for the equation of spherically symmetrical motion of viscous gas,
%         Japan J. Appl. Math. 10 (1993) 219--235.

%\bibitem{m4}
%         M. Okada, $\breve{S}$. Matu$\check{s}\dot{u}$-Ne$\breve{c}$asov$\acute{a}$, T. Makino,
%         Free boundary problem for the equation of one-dimensional motion of compressible gas with density-dependent viscosity,
%         Ann. Univ. Ferrara Sez. VII (N.S.) 48 (2002) 1-20.

\bibitem{Onsager1} L. Onsager, Reciprocal relations in irreversible processes I, Phys. Rev. 37 (1931) 405--426.

\bibitem{Onsager2} L. Onsager, Reciprocal relations in irreversible processes II, Phys. Rev. 38 (1931) 2265--2279.

\bibitem{Onsager3} L. Onsager, S. Machlup, Fluctuations and irreversible processes, Phys. Rev. 91 (1953) 1505--1512.


\bibitem{Sideris} T.C. Sideris, B. Thomases, Global existence for 3D incompressible isotropic elastodynamics via the
incompressible limit, Comm. Pure Appl. Math. 57 (2004) 1--39.

\bibitem{Simon}
         J. Simon,
         Nonhomogeneous viscous incompressible fluids: existence of vecocity, density, and pressure,
         SIAM J. Math. Anal. 21 (1990) 1093--1117.

\bibitem{sun} H. Sun, C. Liu, On energetic variational approaches in modeling the nematic liquid crystal flows, Discrete Contin. Dynam. Systems 23 (2009) 455--475.

%\bibitem{z3}
%         S.W. Vong, T. Yang, C.J. Zhu,
%         Compressible Navier-Stokes equations with degenerate viscosity coeffcient and vacuum II,
%         J. Diff. Eqs. 192 (2003) 475--501.

\bibitem{Wang-Yu}
         D.H. Wang, C. Yu,
         Incompressible limit for the compressible flow of liquid crystals,
         preprint arXiv:1108.4941 (2011) 1--17.

 \bibitem{Wang-Yu2012}
         D.H. Wang, C. Yu,
         Global weak solution and large-time behavior for the compressible flow of liquid crystals,
         Arch. Rational Mech. Anal. 204 (2012) 881--915.

\bibitem{Wen-Ding}
         H.Y. Wen, S.J. Ding,
         Solutions of incompressible hydrodynamic flow of liquid crystals,
         Nonlinear Anal.: Real World Appl. 12 (2011) 1510--1531.

\bibitem{wxl} H. Wu, X. Xu, C. Liu, On the general Ericksen-Leslie system: Parodi's relation, well-posedness and stablity, preprint arXiv:1105.2180v5 (2011) 1--37.

\bibitem{zlz} S.P. Zhang, C. Liu, H. Zhang, Numerical simulations of hydrodynamics of nematic liquid crystals: effects of kinematic transports, Commun. Comput. Phys. 9 (2011) 974--993.

%\bibitem{8}
%         Z. Xin,
%         Blow-up of smooth solutions to the compressible Navier-Stokes equations with compact density,
%         Comm. Pure Appl.Math. 351 (1998) 229--240.

%\bibitem{z1}
%         T. Yang, Z.A. Yao, C.J. Zhu,
%         Compressible Navier-Stokes equations with density-dependent viscosity and vacuum,
%         Comm.Partial Differential Equations, 26 (2001) 965--981.
%
%\bibitem{z2}
%         T. Yang, C.J. Zhu,
%         Compressible Navier-Stokes equations with degenerate viscosity coefcient and vacuum,
%         Comm. Math. Phys. 230 (2002) 329--363.



\end{thebibliography}
\end{document}